\documentclass[journal,twoside,web]{ieeecolor}
\usepackage[T1]{fontenc}
\usepackage{generic}
\usepackage{cite}
\usepackage{amsmath,amssymb,amsfonts}
\usepackage{mathtools}
\usepackage{graphicx}
\usepackage{algorithm,algpseudocode}
\usepackage{caption}
\usepackage{subcaption}
\usepackage{textcomp,gensymb}

\newtheorem{thm}{\protect\theoremname}
\newtheorem{defn}[thm]{\protect\definitionname}
\newtheorem{example}[thm]{\protect\examplename}
\newtheorem{rem}[thm]{\protect\remarkname}
\newtheorem{lem}[thm]{\protect\lemmaname}
\newtheorem{problem}[thm]{\protect\problemname}
\newtheorem{corr}[thm]{\protect\corautorefname}

\providecommand{\definitionname}{Definition}
\providecommand{\examplename}{Example}
\providecommand{\lemmaname}{Lemma}
\providecommand{\problemname}{Problem}
\providecommand{\remarkname}{Remark}
\providecommand{\theoremname}{Theorem}

\providecommand{\assumpname}{Assumption}

\providecommand{\corautorefname}{Corollary}

\providecommand{\assumpname}{Assumption}

\global\long\def\norm#1{\lVert#1\rVert}%
\global\long\def\trace#1{\text{trace}(#1)}%
\global\long\def\abs#1{\lvert#1\rvert}%
\global\long\def\mtl{\mathcal{L}}%
\global\long\def\mtx{\mathcal{X}}%
\global\long\def\mtp{\mathcal{P}}%
\global\long\def\mte{\mathcal{E}}%
\global\long\def\vecc#1{\text{vec}(#1)}%
\global\long\def\mtf{\mathcal{F}}%
\global\long\def\mtn{\mathcal{N}}%
\global\long\def\mtm{\mathcal{M}}%
\global\long\def\mtu{\mathcal{U}}%
\global\long\def\diag{\mathrm{diag}}%
\global\long\def\fdyn{\text{Direct}}%
\global\long\def\mat#1{\mathrm{mat}(#1)}%
\global\long\def\blkdiag{\mathrm{blkdiag}}%

\DeclareMathOperator*{\esssup}{ess\,sup}

\usepackage{comment}

\def\BibTeX{{\rm B\kern-.05em{\sc i\kern-.025em b}\kern-.08em
    T\kern-.1667em\lower.7ex\hbox{E}\kern-.125emX}}
\begin{document}
\title{Continuous-time Constrained Funnel Synthesis for Incrementally Quadratic Nonlinear Systems}
\author{Taewan Kim, Dayou Luo, and Beh{\c{c}}et A{\c{c}}{\i}kme{\c{s}}e
\thanks{Taewan Kim is with the Department of Aeronautics \& Astronautics, University of Washington, Seattle, WA 98115 USA (e-mail: twankim@uw.edu)}
\thanks{Dayou Luo is with the Department of Aeronautics \& Astronautics, University of Washington, Seattle, WA 98115 USA (e-mail: twankim@uw.edu)}
\thanks{Beh{\c{c}}et A{\c{c}}{\i}kme{\c{s}}e is with the Department of Aeronautics \& Astronautics, University of Washington, Seattle, WA 98115 USA (e-mail: twankim@uw.edu)}}

\maketitle

\begin{abstract}
This paper presents a convex optimization-based framework for synthesizing time-varying controlled invariant funnels and associated feedback control around a given nominal trajectory for nonlinear systems subject to bounded disturbances. Nonlinearities are modeled using incremental quadratic constraints, including Lipschitz, L-smooth, and sector-bounded nonlinearities. Funnel invariance is ensured via a DLMI. Together with pointwise-in-time LMIs for state and input constraints, we formulate a continuous-time funnel synthesis problem. To solve it using numerical optimal control techniques, the DLMI is reformulated into a differential matrix equality (DME) and an LMI, where the DME acts as a funnel dynamics equation. We explore different formulations of these funnel dynamics. Continuous-time constraint satisfaction is addressed through two convex methods: one based on intermediate constraint-checking points, and another using a successive convexification method with subgradients to handle nondifferentiable maximum eigenvalue functions. Theoretical justification is provided for the existence of a measurable and integrable subgradient for the latter. The method is demonstrated on two numerical examples: the control of a unicycle and a 6-degree-of-freedom quadrotor for obstacle avoidance.
\end{abstract}


\section{Introduction}
\label{sec:introduction}
\IEEEPARstart{T}{rajectory} generation, also known as motion planning or guidance, has been successfully applied to a variety of tasks \cite{tassa2012synthesis,malyuta2021advances}. One drawback of employing these trajectories is that, because they are merely open-loop state and input signals, they provide no robustness to uncertainty or guarantees that even a single feedback controller can reliably track them while satisfying safety constraints. This paper aims to strengthen the robustness of nominal trajectories by computing robust controlled invariant funnels and associated feedback control around them. Then, the computed funnels serves as a certificate that, for any state within them, there exists a control signal, provided by the associated controller, that ensures constraint satisfaction, and consequently, the system's safety, under the presence of disturbances.

The computed funnels and associated controllers can be used for various purposes. First, they enable constrained robust control for uncertain systems, as they guarantee the closed-loop system's constraint satisfaction under the presence of uncertainties \cite{accikmecse2011robust}. Second, the funnels can be used to generate feasible trajectories \cite{reynolds2021funnel} by numerically integrating the system dynamics with the synthesized controller. This eliminates the need to repeatedly solve the trajectory optimization problem when problem data, such as the initial condition, changes. This can address a common limitation of many nonconvex trajectory optimization methods, which often lack guarantees that the converged solution satisfies constraints \cite{malyuta2022convex}. Furthermore, the computed funnels can be used to ensure stability of recursive feasibility of robust model predictive control (MPC) \cite{BehcetMPCResolvability}. By enforcing a terminal constraint that requires the predicted final state to lie within the funnel, the method provides a natural mechanism to maintain feasibility across MPC iterations. In addition, the computed funnels can be used for feasible-trajectory data generation, which is useful for training neural network controllers under safety constraints \cite{kim2022guided}.

Due to these benefits of funnels, the problem of computing funnels has been actively studied in both the control \cite{mayne2005robust,accikmecse2011robust} and robotics \cite{burridge1999sequential,majumdar2017funnel} communities. Despite this progress, several key challenges remain. First, extending funnel synthesis to broader classes of nonlinear systems requires going beyond linear or polynomial dynamics by exploiting structural properties of nonlinearities. Second, ensuring invariance, which is typically derived through Lyapunov methods \cite{corless1993control}, Hamilton-Jacobi reachability analysis \cite{bansal2017hamilton}, or direct uncertainty propagation \cite{seo2021fast}, is often computationally expensive, particularly for high-dimensional or nonlinear systems. Third, achieving constraint satisfaction for all states and inputs within the associated state and input funnels is nontrivial and is often handled indirectly by adjusting the nominal trajectory, rather than explicitly enforcing constraints on the funnel itself. Finally, because funnels are inherently time-varying, they are commonly represented using basis function approximations, such as zero-order hold or first-order hold interpolations, which can limit expressiveness.

This paper presents a convex optimization-based method for computing the controlled invariant funnels and their associated feedback controllers using linear matrix inequalities (LMIs) \cite{boyd1994linear}. We consider nonlinear systems subject to external bounded disturbances, where nonlinearities satisfy incremental quadratic constraints ($\delta$QC) \cite{accikmecse2011observers,d2012incremental,xu2020observer} that can capture a broader class of nonlinearities, including Lipschitz, L-smooth, and sector-bounded nonlinearities. The state funnels are defined as time-varying ellipsoidal sets around a nominal trajectory, described by a quadratic Lyapunov function with a positive definite (PD) matrix-valued function, as decision variables. Alongside the nominal input, a linear time-varying (LTV) feedback controller is applied, with the LTV gain matrices serving as another decision variables. Funnel invariance is ensured by Lyapunov theory, resulting in a differential linear matrix inequality (DLMI). To ensure constraint satisfaction, all states and inputs within the funnel must satisfy state and input constraints. When these constraints are linear or locally linearized around the nominal trajectory, they can be encoded as pointwise-in-time LMIs.

The resulting funnel synthesis problem, involving a DLMI for invariance and pointwise in time LMIs for state and input constraints, is an infinite-dimensional problem due to the time-varying nature of the decision variables. To make the problem tractable, we reformulate the DLMI using a slack variable, transforming it into a differential matrix equality (DME) and a pointwise-in-time LMI. We refer to this DME as funnel dynamics, as it behaves analogously to system dynamics in an optimal control framework. This reformulation enables us to cast the funnel synthesis problem into a standard form of numerical optimal control so that we can apply numerical optimal control techniques such as multiple-shooting \cite{bock1984multiple} and control parameterization \cite{lin2014control}. The resulting funnel synthesis is a finite-dimensional convex semidefinite program (SDP), which can be efficiently solved using standard SDP solvers.

The proposed framework further incorporates continuous-time constraint satisfaction (CTCS) by extending enforcement beyond discrete node points. We introduce two methods. The first enforces pointwise-in-time LMIs at additional intermediate time points within each subinterval, improving constraint coverage without increasing the number of decision variables. The second approach follows a trajectory optimization technique commonly used in prior works \cite{lin2014control,ELANGO2025112464}, where continuous-time constraints are reformulated by integrating their violations using an exterior penalty function. This reformulation retains convexity but results in constraints that are not in a standard conic form, thus requiring iterative methods such as successive convexification (SCvx). In our setting, however, the constraint is expressed as an LMI involving the maximum eigenvalue function, which makes the constraints nondifferentiable. As a result, we explore the use of subgradient-based linearization within the SCvx framework to handle the lack of differentiability. 

\subsection{Related work}

The computation of funnels has been studied in both robotics and control communities. To the best of our knowledge, the notion of the funnel was first introduced in \cite{mason1985mechanics} where the term was used as a geometric metaphor to describe how mechanical interactions can guide a range of initial states toward a desired final configuration. The Lyapunov-theoretic interpretation where funnels represent invariant sets was later adopt in \cite{burridge1999sequential}. Following this direction, a significant body of work has focused on offline synthesis of funnels as time-varying controlled invariant sets \cite{tedrake2010lqr,tobenkin2011invariant,majumdar2017funnel}. Most existing approaches rely on sum-of-squares (SOS) programming, where the system dynamics are approximated by polynomials, and high-order polynomial Lyapunov functions are employed. 

While our work also addresses offline computation of time-varying funnels, the key distinctions lie in exactness, computational efficiency, and explicit constraint satisfaction. Specifically, we employ quadratic Lyapunov functions and locally approximate the nonlinear dynamics around the nominal trajectory. The resulting approximation residual is captured using $\delta$QC, allowing us to certify exact invariance with respect to the original nonlinear system. Unlike SOS-based methods, which require solving nonconvex formulations due to bilinear terms even when the controller is fixed, our approach yields a semidefinite program (SDP) with only a few scalar parameters that can be efficiently tuned via line search. In contrast to existing approaches that often fix the controller, we optimize both the funnel and the feedback controller simultaneously, while preserving convexity. Existing methods typically do not incorporate state and input constraints within the funnel, aside from input saturation, and instead adjust the nominal trajectory to indirectly achieve feasibility. Our approach provides handle such constraints directly within the synthesis process. 

In the control community, controlled invariant funnels, also referred to as tubes, have been studied extensively in the context of tube-based offline finite-horizon control and online MPC. Earlier work in this area focused on linear systems subject to disturbances and stability around an equilibrium point \cite{LANGSON2004125,mayne2005robust,Rakovic2025Invariant}. As the field evolved, approaches have extended toward trajectory tracking for nonlinear systems \cite{accikmecse2011robust}. In these problems, the main solution variables include the nominal trajectory, the invariant set (funnel), and the feedback controller, typically leading to nonconvex formulations. For tractability in MPC settings, some of these variables are fixed in advance and assumed to be time-invariant for simplicity \cite{villanueva2017robust,kohler2020computationally}. In contrast, several offline approaches optimize all three jointly, leveraging the additional flexibility to improve performance \cite{kim2022joint,leeman2023robust1}. In this paper, we focus on computing time-varying invariant funnels and associated feedback control with a given nominal trajectory. Our main distinction to the existing work lies in addressing how continuous-time constraints can be satisfied between discretization nodes. Most prior methods are developed for discrete-time systems without considering constraint violations between nodes, thereby overlooking potential violations that can arise in continuous-time settings.

A related line of work studies the analysis and synthesis of uncertain LTV systems \cite{seiler2014stability,seiler2019finite,buch2021finite,seiler2024trajectory} using integral quadratic constraints (IQC) . The approach in \cite{seiler2019finite} presents a computational framework for assessing finite-horizon robust performance of uncertain LTV systems with uncertainties described by IQCs. Building upon this, \cite{buch2021finite} extends the results to the synthesis of robust output-feedback controllers where Lyapunov parameters and control parameters are iteratively updated. Our key distinction is the focus on invariant set computation, which necessitates explicit consideration of state and input constraints. In addition, we employ state-feedback and formulate the problem as a convex semidefinite program (SDP) where both funnel and control parameters are computed together.

There exists a body of literature on funnel control \cite{Ilchmann_Ryan_Sangwin_2002,berger2023funnel}. Although our work and theirs share the term \textquotedbl funnel\textquotedbl , the approaches are quite different. (1) The definition of a funnel differs: in funnel control, it refers to a prescribed bound on the output tracking error over time, while in our work it denotes a computed controlled-invariant tube around a trajectory in state space. (2) The control law differs: funnel control uses a universal adaptive high-gain error feedback, whereas we synthesize a trajectory-specific feedback via convex optimization. (3) The control objective differs: funnel control focuses on output tracking with guaranteed transient behavior across a class of systems, while our goal is to certify robust trajectory tracking with formal guarantees under model uncertainty/nonlinearity and constraints. (4) The frameworks differ: funnel control is analytical and model-free, while ours is model-based and optimization-based.

This paper builds upon the foundational work on LMI-based robust control developed in the 2000s \cite{accikmecse2002robust,accikmecse2003robust,accikmecse2008stability}. This foundational line of work has recently been extended to LMI-based observer and controller design for incrementally quadratic nonlinear systems \cite{xu2020observer}. A preliminary version of the present work was introduced in \cite{10167750}. This paper introduces several substantial extensions: the class of nonlinear systems is expanded from Lipschitz to more general sector-bounded nonlinearities; a DLMI ensuring funnel invariance for such nonlinear systems is derived; a different formulation of the funnel dynamics is explored; and continuous-time constraint satisfaction (CTCS) is investigated.

\subsection{Statement of contributions}

\noindent \textbf{1)} We extend funnel synthesis to a broad class of nonlinear systems by incorporating time-varying $\delta$QC. In contrast to previous works limited to Lipschitz nonlinearities \cite{reynolds2021funnel}, we consider more general sector-bounded nonlinearities. In particular, we demonstrate how L-smooth nonlinearities can be explicitly handled within the $\delta$QC framework for reduced conservativeness.

\noindent \textbf{2)} We derive a DLMI that implies a dissipation inequality in the form of an input-to-state Lyapunov condition, ensuring the funnel's invariance. This generalizes existing LMI-based robust state-feedback synthesis methods \cite{accikmecse2011robust} by allowing the nonlinearity to depend explicitly on external bounded disturbances.

\noindent \textbf{3)} The paper presents a numerical optimal control-based solution method for solving the continuous-time funnel synthesis problem, in which the DLMI appears as one of the constraints. Unlike existing approaches that are restricted to represent funnels using zeroth-order hold or first-order hold interpolations, the proposed method supports higher-order representations, including second-order hold. It also allows the funnel profile to be aligned with the solution of a Lyapunov differential equation.

\noindent \textbf{4)} The proposed method systematically addresses CTCS using two convex approaches. While prior work often increases discretization nodes \cite{tobenkin2011invariant,majumdar2017funnel}, which leads to more decision variables, the first approach adds intermediate constraint-checking points without increasing variable count. The second reformulates the continuous-time LMI using an integral representation of constraint violation using an exterior penalty function, leading to nondifferentiable constraints involving the maximum eigenvalue function. To solve this, we introduce a subgradient-based SCvx method and establish theoretical support for the existence of subgradients.

\subsection{Notation}

We denote by $\mathbb{R}$ the set of real numbers, $\mathbb{R}_{+}$ the set of nonnegative real numbers, $\mathbb{R}_{++}$ the set of positive real numbers, $\mathbb{Z}$ the set of integers, and $\mathbb{R}^{n}$ the set of $n$-dimensional real column vectors. The set $\mathbb{Z}_{[a,b)}=\{z\in\mathbb{Z}:a\leq z<b\}$denotes the integers in the half-closed interval from a $a$ to b $b$. The set $\mathbb{S}^{n}$ denotes the set of $n\times n$ real symmetric matrices, $\mathbb{S}_{+}^{n}$ the set of $n\times n$ positive semidefinite (PSD) matrices, $\mathbb{S}_{++}^{n}$ the set of $n\times n$ PD matrices, and $\mathbb{S}^{n}$ the set of $n\times n$ real symmetric matrices. The set $\mtl_{2}[a,b]$ denotes the space of Lebesgue measurable functions $x(t)$ defined on an interval $[a,b]\subset\mathbb{R}$ such that $\left(\int_{a}^{b}x(t)^{\top}x(t)\mathrm{d}t\right)^{1/2}<\infty.$ A property is said to hold almost everywhere (also expressed as holding for almost every $x\in X$, or almost all $x\in X$) if the set of points in a measurable space $X$ where it fails has Lebesgue measure zero. For a measurable function $f:[t_{0},t_{f}]\to\mathbb{R}$, the essential supremum of $f$, denoted by $\esssup_{t\in[t_{0},t_{f}]}f(t)$, is the smallest real number $M$ such that $f(t)\leq M$ for almost all $t\in[t_{0},t_{f}]$. For notational brevity, we omit repeated symmetric terms and write expressions such as $(\star)^{\top}PA$ to indicate $A^{\top}PA$. In block matrices, we write $\left[\begin{array}{cc}
a & \star\\
b & c
\end{array}\right]$ where $\star$ represents the transpose of the corresponding off-diagonal block (e.g., $b^{\top}$). We use $\oplus$ to denote Minkowski sum and $\times$ to denote the Cartesian product of sets. The zero matrices having $m\times n$ size and identity matrix having $n\times n$ size are denoted by $0_{m\times n}$ and $I_{n}$, respectively. The subscript will be omitted when it is clear from the context. We use $\text{diag}(\cdot)$ to denote a diagonal matrix formed from its arguments, $\vecc{\cdot}$ to denote vectorization of a matrix by stacking its columns, and $\mat{\cdot}$ as the inverse operation of $\vecc{\cdot}$, reshaping a vector into a matrix. For time-varying signals$a(t)$ and $b(t)$, we write $(a,b)$ to denote the pair of functions $\{a(\cdot),b(\cdot)\}$. We omit the time variable $t$ when it is either clear from context or not essential to the discussion.



\section{Nonlinear Systems, Funnels, and Invariance Conditions\label{sec:2}}

\subsection{Nonlinear systems}

Consider continuous-time nonlinear dynamical systems of the form
\begin{align}
\dot{x}(t) & =f(t,x(t),u(t),w(t)),\quad t\in[t_{0},t_{f}],\label{eq:nonlinear_system}
\end{align}
where $f:[t_{0},t_{f}]\times\mathbb{R}^{n_{x}}\times\mathbb{R}^{n_{u}}\times\mathbb{R}^{n_{w}}\rightarrow\mathbb{R}^{n_{x}}$ is assumed to be continuous in $t$ and continuously differentiable in other arguments. The vectors $x(t)\in\mathbb{R}^{n_{x}}$ and $u(t)\in\mathbb{R}^{n_{u}}$ are the state and the control input, and $w(t)\in\mathbb{R}^{n_{w}}$ is the unknown, but bounded exogeneous disturbance. We assume that $u(\cdot)\in\mtl_{2}^{n_{u}}[t_{0},t_{f}]$ is piecewise continuous, and $w(\cdot)\in\mtl_{2}^{n_{w}}[t_{0},t_{f}]$ is continuous almost everywhere and essentially bounded, that is, 
\begin{equation}
\norm{w(\cdot)}_{\infty}\coloneqq\esssup_{t\in[t_{0},t_{f}]}\norm{w(t)}_{2} \leq w_{max},\label{eq:w_pointwise}
\end{equation}
for some $w_{max}\in\mathbb{R}_{+}$. The time instances $t_{0}$ and $t_{f}$ are the initial and the final time, respectively, with $t_{0}\leq t_{f}<\infty$.

Without loss of generality, we express the model \eqref{eq:nonlinear_system} in the linear fractional form:
\begin{subequations}
\label{eq:linear_fractional_form}
\begin{align}
\dot{x}(t) & =A_{o}(t)x+B_{o}(t)u+F_{o}(t)w+E\phi(t,q_{o}),\\
q_{o}(t) & =C_{o}x+D_{o}u+G_{o}w,  
\end{align}
\end{subequations}
where $A_{o}:[t_{0},t_{f}]\rightarrow\mathbb{R}^{n_{x}\times n_{x}},B_{o}:[t_{0},t_{f}]\rightarrow\mathbb{R}^{n_{x}\times n_{u}},F_{o}:[t_{0},t_{f}]\rightarrow\mathbb{R}^{n_{x}\times n_{w}}$ are continuous in time and uniformly bounded on $[t_0, t_f]$. A pair $(q_{o},\phi)$ represents the system's nonlinearity, where $\phi:[t_{0},t_{f}]\times\mathbb{R}^{n_{q_{o}}}\rightarrow\mathbb{R}^{n_{\phi}}$ is a known nonlinear function with its argument $q_{o}(t)\in\mathbb{R}^{n_{q_{o}}}$ that is a linear function of $x$, $u$, and $w$. The matrices $C_{o},D_{o},$ $G_{o}$, and $E$ are selector matrices used to structure $q_{o}$. It is worth noting that the linear fractional form \eqref{eq:linear_fractional_form} is general enough to illustrate \eqref{eq:nonlinear_system}, while isolating the nonlinearities for focused analysis. This is because simply selecting $E=I$, $\phi=f$, $q_{o}=[x^{\top},u^{\top},w^{\top}]^{\top}$ and setting $A_{o},B_{o},$ and $C_{o}$ to zero matrices of appropriate sizes recover the original system representation \eqref{eq:nonlinear_system}.

The nominal trajectory refers to a collection of state and open-loop control input trajectories $\{\bar{x}(t),\bar{u}(t)\}_{t=t_{0}}^{t_{f}}$, compactly denoted by $(\bar{x},\bar{u})$, satisfying the original system dynamics, that is, 
$
\dot{\bar{x}}(t)=f(t,\bar{x}(t),\bar{u}(t),0),
$
where the nominal disturbance $\bar{w}(t)$ is set as the zero vector for all $t$.

Now, we focus on the incremental behavior of the system, that is, how the state evolves relative to a nominal trajectory, also referred to as the difference dynamics. This is essential because the funnel is centered around the nominal trajectory, not the origin. To formalize this, define deviation variables as
\begin{align*}
\eta & \coloneqq x-\bar{x},\quad\xi\coloneqq u-\bar{u},\quad\delta q_{o}\coloneqq q_{o}-\bar{q}_{o},
\end{align*}
where $\bar{q}_{o}=C_{o}\bar{x}+D_{o}\bar{u}$.  Substituting into the linear fractional form yields:
\begin{subequations}
\label{eq:incremental_dyn_1}
\begin{align}
\dot{\eta}(t) & =A_{o}(t)\eta+B_{o}(t)\xi +F_{o}(t)w +E\delta\phi(t,\delta q_{o}),\\
\delta q_{o}(t) & =C_{o}\eta +D_{o}\xi +G_{o}w,\label{eq:def_q0}
\end{align}
\end{subequations}
where $\delta\phi(t,\delta q_{o})\coloneqq\phi(t,q_{o})-\phi(t,\bar{q}_{o})$. The expression above shows how the deviation evolves due to input deviations $\xi$, the disturbance $w$, and the incremental nonlinearity $\delta\phi$, which captures the difference in the nonlinear term relative to the nominal trajectory. We will further refine the expression \eqref{eq:incremental_dyn_1} using the Mean Value Theorem in the subsequent section.

\subsection{Funnel definition and invariance conditions}
We consider a scalar-valued quadratic storage function (Lyapunov function) $V:[t_{0},t_{f}]\times\mathbb{R}^{n_{x}}\rightarrow\mathbb{R}_{+}$ defined by
\begin{equation}
V(t,\eta(t))=\eta(t)^{\top}Q(t)^{-1}\eta(t),\label{eq:Lyapunov_function}
\end{equation}
where $Q:[t_{0},t_{f}]\rightarrow\mathbb{S}_{++}^{n_{x}}$ is piecewise continuous and continuously differentiable. The positive definiteness of $Q(t)$ ensures quadratic bounds, i.e., 
$\alpha_{1}\|\eta(t)\|^{2} \leq V(t,\eta(t)) \leq \alpha_{2}\|\eta(t)\|^{2}$ 
for some $\alpha_{1},\alpha_{2}>0$ and all $t\in[t_{0},t_{f}]$.

A \emph{state funnel} is defined as the sub-level set of $V$:
\begin{equation}
\mte_{\eta}(t)=\{\eta\mid\eta^{\top}Q(t)^{-1}\eta\leq c_{Q}\},\label{eq:state_funnel}
\end{equation}
where $c_{Q}\in\mathbb{R}_{+}$ is a support value that specifies the level constant of $V$ for $\mte_{\eta}$. We apply a linear time-varying feedback controller for the incremental system \eqref{eq:incremental_dyn_1} given by $\xi(t)=K(t)\eta(t)$ with the piecewise continuous feedback gain $K:[t_{0},t_{f}]\rightarrow\mathbb{R}^{n_{u}\times n_{x}}$. Then, by \cite[Lemma A.1.]{TaewanThesis}, $\eta\in\mte_{\eta}(t)$ implies $K(t)\eta\in\mte_{\xi}(t)$ where the \emph{input funnel} $\mte_{\xi}(t)\subseteq\mathbb{R}^{n_{u}}$ is defined as
\begin{equation}
\mte_{\xi}(t)=\{(c_{Q}K(t)Q(t)K(t)^{\top})^{\frac{1}{2}}y\mid\norm y_{2}\leq1,y\in\mathbb{R}^{n_{u}}\}.\label{eq:input_funnel}
\end{equation}
With the state and input funnels, a full \emph{funnel} is defined as
\begin{equation}
\mtf(t)=(\{\bar{x}(t)\}\oplus\mte_{\eta}(t))\times(\{\bar{u}(t)\}\oplus\mte_{\xi}(t)).\label{eq:funnel}
\end{equation}
That is, the funnel is a Cartesian product of state and input funnels, shifted to the nominal trajectory.
\begin{defn}
\label{def:funnel_invariance}The funnel defined in \eqref{eq:funnel} is invariant for the nonlinear system \eqref{eq:nonlinear_system} in a sense that if $x(\cdot)$ is a solution of \eqref{eq:nonlinear_system} with $x(t_{0})\in\{\bar{x}(t_{0})\}\oplus\mte_{\eta}(t_{0})$, the control law
\begin{equation}
u(t)=\bar{u}(t)+K(t)(x(t)-\bar{x}(t)),\label{eq:control_law}
\end{equation}
and $w(\cdot)\in\mtl_{2}^{n_{w}}[t_{0},t_{f}]$ satisfying \eqref{eq:w_pointwise}, then $(x(t),u(t))\in\mtf(t)$ for all $t\in[t_{0},t_{f}]$.
\end{defn}

\begin{defn}
\label{def:funnel_feasibility} Given a state constraint set $\mtx\subset\mathbb{R}^{n_{x}}$ and an input constraint set $\mtu\subset\mathbb{R}^{n_{u}}$, the funnel defined in \eqref{eq:funnel} is feasible if 
\begin{equation}
\mtf(t)\subseteq\mtx\times\mtu,\quad\forall t\in[t_{0},t_{f}].\label{eq:funnel_feasibility}
\end{equation}
\end{defn}

\begin{defn}
Given a fixed nominal trajectory $(\bar{x},\bar{u})$, funnel synthesis refers to a procedure for computing the functions $Q:[t_{0},t_{f}]\rightarrow\mathbb{S}_{++}^{n_{x}}$ and $K:[t_{0},t_{f}]\rightarrow\mathbb{R}^{n_{u}\times n_{x}}$, ensuring the invariance and feasibility of funnel as defined in Definitions \ref{def:funnel_invariance} and \ref{def:funnel_feasibility}, respectively.
\end{defn}

\section{Continuous-time Funnel Synthesis Problem\label{sec:3}}

\subsection{Structured nonlinearity and mean value theorem}

Recall that the nonlinear function $\phi$ in \eqref{eq:linear_fractional_form} is given by
\[
\phi(t,q_{o})=\left[\begin{array}{c}
\phi_{1}(t,q_{1})\\
\vdots\\
\phi_{n_{\phi}}(t,q_{n_{\phi}})
\end{array}\right],
\]
where each $\phi_{i}:[t_{0},t_{f}]\times\mathbb{R}^{n_{q_{i}}}\rightarrow\mathbb{R}$ corresponds to the $i$-th row of $\phi$, and $q_{i}\in\mathbb{R}^{n_{q_{i}}}$ is a subset of the entries of $q_{o}$. That is, $q_{i}$ is not necessarily equal to $q_{o}$, but rather represents the specific argument passed to $\phi_{i}$ drawn from $q_{o}$, since $\phi_{i}$ may depend only on part of $q_{o}$.

For a more compact representation, we introduce an index $[i]$ for all $i\in\mathbb{Z}_{[1,n_{c}]}$ to represent the $i$-th nonlinearity channel, which is a group of component indices treated together in the reformulation. Let $\phi_{[i]}$ denote the stacked function composed of all $\phi_{j}$ with $j\in[i]$. The associated argument $q_{[i]}\in\mathbb{R}^{n_{q_{[i]}}}$ is defined as the minimal vector containing all distinct variables that appear as arguments of the functions $\phi_{j}$ in channel $[i]$. If two or more $\phi_{j}$ in the same channel share the same argument, that argument appears only once in $q_{[i]}$. For example, if $\phi_{1}$ and $\phi_{2}$ are in channel $[1]$ and $q_{1}=q_{2}$, then $q_{[1]}=q_{1}$ not $[q_{1},q_{2}]^{\top}$. We define selector matrices $C_{i}$, $D_{i}$, $G_{i}$ so that $q_{[i]}=C_{i}x+D_{i}u+F_{i}w$. This reformulation is called structured nonlinearity, whose details can be found in \cite[Chapter 6.2.1]{reynolds2020computational}.

Let $\bar{q}_{[i]}\coloneqq C_{[i]}\bar{x}+D_{[i]}\bar{u}$ be the nominal value of $q_{[i]}$, and define the increment $\delta q_{[i]}=q_{[i]}-\bar{q}_{[i]}$. Applying the Mean Value Theorem \cite{Rudin1976} to each channel gives
\[
\phi_{[i]}(t,q_{[i]})=\phi_{[i]}(t,\bar{q}_{[i]})+\phi_{[i]}'(t,\tilde{q}_{[i]})\delta q_{[i]},
\]
where $\tilde{q}_{[i]}(t)\in\mathbb{R}^{n_{q_{[i]}}}$ lies on the line segment between $q_{[i]}(t)$ and $\bar{q}_{[i]}(t)$, i.e., $\tilde{q}_{[i]}\in\text{Co}(\{q_{[i]},\bar{q}_{[i]}\})$, where $\text{Co(\ensuremath{\cdot)}}$ denotes the convex hull. Consider the Jacobian $\phi_{[i]}':[t_{0},t_{f}]\times\mathbb{R}^{n_{q_{[i]}}}\rightarrow\mathbb{R}^{1\times n_{q_{[i]}}}$ with respect to the second argument. By adding and subtracting the term $\phi_{[i]}'(t,\bar{q}_{[i]})\delta q_{[i]}$ to the right-hand side, we obtain the followings:
\begin{subequations}
\label{eq:mean_value_theorem}
\begin{align}
\phi_{[i]}(t,q_{[i]}) & =\phi_{[i]}(t,\bar{q}_{[i]})+\phi{}_{[i]}'(t,\bar{q}_{[i]})\delta q_{[i]}+\delta p_{[i]},\label{eq:phi_expansion}\\
\delta p_{[i]}(t) & =\left(\phi_{[i]}'(t,\tilde{q}_{[i]})-\phi{}_{[i]}'(t,\bar{q}_{[i]})\right)\delta q_{[i]},\label{eq:delta_po}
\end{align}
\end{subequations}
for each $i$ in $\mathbb{Z}_{[1,n_{c}]}$ where the function $\delta p_{[i]}(t)\in\mathbb{R}^{n_{p_{[i]}}}$ and $\sum_{i=1}^{n_{c}}n_{p_{[i]}}=n_{\phi}$. 

The incremental form of dynamics \eqref{eq:incremental_dyn_1} can then be rewritten as
\begin{subequations}
\label{eq:incremental_dyn_2}
\begin{align}
\dot{\eta}(t) & =A(t)\eta+B(t)\xi+F(t)w+E\delta p(t,\delta q,\\
\delta q(t) & =C\eta+D\xi+Gw,
\end{align}
\end{subequations}
where the vectors $\delta q\in\mathbb{R}^{n_{q}}$ and $\delta p\in\mathbb{R}^{n_{p}}$ are the stacked vector of all channel-wise $\delta q_{[i]}$ and $\delta p_{[i]}$, respectively, with total dimension $n_{q}=\sum_{i=1}^{n_{c}}n_{q_{[i]}}$ and $n_{p}=\sum_{i=1}^{n_{c}}n_{p_{[i]}}$. The constant matrices $C$, $D$, and $G$ are defined by the vertical stacking of individual selector matrices:
\[
C=\left[\begin{array}{c}
C_{1}\\
\vdots\\
C_{n_{c}}
\end{array}\right],\quad D=\left[\begin{array}{c}
D_{1}\\
\vdots\\
D_{n_{c}}
\end{array}\right],\quad G=\left[\begin{array}{c}
G_{1}\\
\vdots\\
G_{n_{c}}
\end{array}\right].
\]

The time-varying matrices $A(t),B(t),$ and $F(t)$ are given by
\begin{subequations}
\label{eq:A_B_F}
\begin{align}
A(t) & =A_{o}(t)+E\phi'(t)C_{o},\\
B(t) & =B_{o}(t)+E\phi'(t)D_{o},\\
F(t) & =F_{o}(t)+E\phi'(t)G_{o},
\end{align}
\end{subequations}
where $\phi'(t) \coloneqq \bigl[\phi'_{[i]}(t,\bar{q}_{[i]}),\;\ldots,\;\phi'_{[n_c]}(t,\bar{q}_{[n_c]})\bigr]^{\top}$.


\begin{example}
\label{exa:unicycle}Consider the dynamics of a unicycle model given by
\[
\dot{x}=\left[\begin{array}{c}
\dot{x}_{1}\\
\dot{x}_{2}\\
\dot{x}_{3}
\end{array}\right]=\left[\begin{array}{c}
u_{1}\cos x_{3}\\
u_{1}\sin x_{3}\\
u_{2}
\end{array}\right],
\]
where $x_{1}$ and $x_{2}$ are $x$- and $y$-positions, $x_{3}$ is the yaw angle, and $u_{1}$ and $u_{2}$ are the linear and the angular velocities, respectively. This model can be written in the form of \eqref{eq:incremental_dyn_1} using the following matrices:
\begin{align*}
A_{o} & =0_{3\times3},\quad B_{o}=\left[\begin{array}{c}
0_{2\times2}\\
\begin{array}{cc}
0 & 1\end{array}
\end{array}\right],\quad E=\left[\begin{array}{c}
I_{2}\\
0_{1\times2}
\end{array}\right],\\
C_{o} & =\left[\begin{array}{cc}
0_{2\times2} & \begin{array}{c}
1\\
0
\end{array}\end{array}\right],\quad D_{o}=\left[\begin{array}{cc}
0 & 0\\
1 & 0
\end{array}\right],
\end{align*}
and $q_{o}=[x_{3},u_{1}]^{\top}$. The nonlinear function is given by $\phi=[u_{1}\cos x_{3},u_{1}\sin x_{3}]^{\top}$. Choosing channels $\phi_{[1]}=\phi_{1}$ and $\phi_{[2]}=\phi_{2}$, we have $q_{o}=q_{[1]}=q_{[2]}$, $C_{1}=C_{2}=C_{o}$, and $D_{1}=D_{2}=D_{o}$. The stacked Jacobian evaluated on the nominal trajectory is 
\begin{equation}
\phi'(t)=\left[\begin{array}{c}
\phi{}_{[1]}'(t,\bar{q}_{[1]})\\
\phi{}_{[2]}'(t,\bar{q}_{[2]})
\end{array}\right]=\left[\begin{array}{cc}
-\bar{u}_{1}\sin\bar{x}_{3} & \cos\bar{x}_{3}\\
\bar{u}_{1}\cos\bar{x}_{3} & \sin\bar{x}_{3}
\end{array}\right]\label{eq:phi_hat_unicycle}
\end{equation}
Thus, the time-varying matrices $A(t)$ and $B(t)$ become
\[
A(t)=\left[\begin{array}{cc}
0_{3\times2} & \begin{array}{c}
-\bar{u}_{1}\sin\bar{x}_{3}\\
\bar{u}_{1}\cos\bar{x}_{3}\\
0
\end{array}\end{array}\right],\quad B(t)=\left[\begin{array}{cc}
\cos\bar{x}_{3} & 0\\
\sin\bar{x}_{3} & 0\\
0 & 1
\end{array}\right].
\]

\end{example}

\subsection{Characterization of nonlinearity using $\delta$QC}

The incremental form of dynamics system expressed in \eqref{eq:incremental_dyn_2} can be viewed as an uncertain LTV system. The pair $(\delta q,\delta p)$, originating from the nonlinearity $\phi$ in the original system \eqref{eq:nonlinear_system}, now appears as an uncertainty in the incremental system. To characterize this uncertainty, we employ the framework of $\delta$QCs \cite{accikmecse2008stability,d2012incremental}.
\begin{defn}
Let $\delta p:[t_{0},t_{f}]\times\mathbb{R}^{n_{q}}\rightarrow\mathbb{R}^{n_{p}}$ be an uncertain nonlinear function for the incremental system \eqref{eq:incremental_dyn_2}. A symmetric matrix-valued function $M:[t_{0},t_{f}]\rightarrow\mathbb{S}^{(n_{p}+n_{q})}$ is called a time-varying incremental multiplier matrix if it satisfies the following $\delta$QC:
\begin{equation}
\left[\begin{array}{c}
\delta q(t)\\
\delta p(t,\delta q(t))
\end{array}\right]^{\top}M(t)\left[\begin{array}{c}
\delta q(t)\\
\delta p(t,\delta q(t))
\end{array}\right]\geq0,\label{eq:def_iQC}
\end{equation}
for all $\delta q\in\mtl_{2}^{n_{q_{o}}}[t_{0},t_{f}]$ and $t\in[t_{0},t_{f}]$.
We denote by $\mtm(t)\subset\mathbb{S}^{n_{p}+n_{q}}$ the set of all incremental multiplier matrices that satisfy this condition for the given nonlinearity $\delta p$ at a fixed time $t\in[t_{0},t_{f}]$.
\end{defn}

\begin{example}
\label{exa:The-incrementally-sector} Let $K_{1},K_{2}:[t_{0},t_{f}]\rightarrow\mathbb{R}^{n_{p}\times n_{q}}$ be two given matrix-valued functions. The nonlinearity $\delta p_{o}$ is said to satisfy a time-varying incremental sector bound if, for some invertible and symmetric weight $S=S^{\top}\in\mathbb{S}^{n_{p}}$, the following inequality holds:
\[
(\delta p(t)-K_{1}(t)\delta q(t))^{\top}S^{-1}(K_{2}(t)\delta q(t)-\delta p(t))\geq0,
\]
for each $t\in[t_{0},t_{f}]$.
This condition admits a $\delta$QC representation as in \eqref{eq:def_iQC} with the time-varying multiplier set 
\begin{equation}
\mathcal{M}(t) = \left\{ \lambda 
\left[\begin{smallmatrix}
- K_{1}^{\top} S^{-1} K_{2} - K_{2}^{\top} S^{-1} K_{1} & \star \\
S^{-1}(K_{1}+K_{2}) & -2S^{-1}
\end{smallmatrix}\right]
\;\middle|\; \lambda \in \mathbb{R}_{+} \right\},
\label{eq:incremental_MM_sector_bounded}
\end{equation}
where the explicit time dependence of $K_{1}, K_{2}$ is omitted for brevity. As a special case, if $\delta p$ satisfied a norm bound $\norm{\delta p(t)}\leq\gamma(t)\norm{\delta q(t)}$ with $\gamma(t)>0$, then \eqref{eq:incremental_MM_sector_bounded} holds with $K_{1}(t)=-\gamma(t)I$, $K_{2}=\gamma(t)I$, and $S=I$.
\end{example}

To enable the convexification of the funnel synthesis problem, we adopt a block-diagonal parameterization of the incremental multiplier matrices \cite{accikmecse2011robust,xu2020observer}. We assume that for every $t\in[t_{0},t_{f}]$, there exist matrices $N_{1}(t)\in\mathbb{S}^{n_{q}}$, $N_{2}(t)\in\mathbb{S}^{n_{p}}$, and a nonsingular matrix $T(t)\in\mathbb{R}^{(n_{q}+n_{p})\times(n_{q}+n_{p})}$ such that
\begin{equation}
T(t)^{\top}\left[\begin{array}{cc}
N_{1}(t)^{-1} & 0\\
0 & -N_{2}(t)^{-1}
\end{array}\right]T(t)\in\mtm(t),\label{eq:block_diagonalization}
\end{equation}
The transformation matrix $T(t)$ is partitioned as
\[
T(t)=\left[\begin{array}{cc}
T_{11}(t) & T_{12}(t)\\
T_{21}(t) & T_{22}(t)
\end{array}\right],
\]
and we assume that $T_{22}(t)$ is nonsingular for all $t\in[t_{0},t_{f}]$. We define the feasible set-valued function $\mtn:[t_{0},t_{f}]\rightarrow\mathbb{S}_{n_{q}}\times\mathbb{S}_{n_{p}}$ as the set of all admissible pairs $(N_{1}(t),N_{2}(t))$ that satisfy the block-diagonal parameterization condition \eqref{eq:block_diagonalization}.
\begin{example}
The incremental multiplier matrices defined in Example \ref{exa:The-incrementally-sector} can be expressed in block-diagonal form using
\begin{equation}
T(t)=\left[\begin{array}{cc}
K_{2}-K_{1} & 0\\
K_{2}+K_{1} & -2I
\end{array}\right],\quad\mtn(t)=\{(\lambda S,\lambda S)\mid\lambda>0\}.\label{eq:multiplier_sectorbounded}
\end{equation}
\end{example}

\subsection{Lipschitz and L-smooth nonlinearities}

Here we provide the characterization of the uncertain term $\delta p$ using the $\delta$QC \eqref{eq:def_iQC}, based on local Lipschitz and L-smooth properties that are two common assumptions for general nonlinear systems. 

Suppose that each nonlinear component $\phi_{[i]}$ has bounded Jacobian $\phi_{[i]}'$ over any compact set $\Omega\subset\mathbb{R}^{n_{q_{[i]}}}$ for each $t\in [t_0,t_f]$ so that $\phi_{[i]}$ is locally pointwise Lipschitz. For each nonlinear channel $i\in\mathbb{Z}_{[1,n_{c}]}$, there exists a time-varying constant $\gamma_{i}(t)\in\mathbb{R}_{+}$ such that
\begin{equation}
\norm{\phi_{[i]}'(t,q_{[i]})-\phi{}_{[i]}'(t,\bar{q}_{[i]})}_{2}\leq\gamma_{i}(t),\ \forall q_{[i]},\bar{q}_{[i]}\in\Omega,\label{eq:Lipshitz}
\end{equation}
for each $t\in[t_{0},t_{f}]$. If $\phi_{[i]}$ is Lipschitz continuous with constant $\gamma_{i}/2$, then its Jacobian is bounded by $\gamma_{i}/2$, and hence the above bound holds with $\gamma_{i}$.

On the other hand, if $\phi$ is locally L-smooth, then its Jacobian is Lipschitz continuous. That is, for any compact set $\Omega\subset\mathbb{R}^{n_{q_{[i]}}}$, there exists a time-varying constant $\beta_{i}(t)\in\mathbb{R}_{+}$ such that for all $i\in\mathbb{Z}_{[1,n_{c}]}$,
\begin{equation}
\norm{\phi_{[i]}'(t,q_{[i]})-\phi{}_{[i]}'(t,\bar{q}_{[i]})}_{2}\leq\beta_{i}(t)\norm{\delta q_{[i]}(t)}_{2}, \ \forall q_{[i]},\bar{q}_{[i]}\in\Omega,\label{eq:Lsmooth}
\end{equation}
for each $t\in[t_{0},t_{f}]$.

\begin{rem}
The representation with \eqref{eq:Lsmooth} presents a direct way to control the level of nonlinearity, whereas the representation with \eqref{eq:Lipshitz} does not exploit this flexibility. It is obvious that the representation with \eqref{eq:Lipshitz} has a constant upper bound of the term, $\norm{\phi_{[i]}'(t,q_{[i]})-\phi{}_{[i]}'(t,\bar{q}_{[i]})}_{2}$, for each $t$ whereas the upper bound in \eqref{eq:Lsmooth} goes to zero as the region around the nominal trajectory shrinks to zero ($\bar{q}_{[i]}\leftarrow q_{[i]}$). This implies that we can reduce the region, in which we seek for a funnel, arbitrarily to reduce the level of nonlinearity to zero and hence increase the chance of finding a funnel. This is not the case with \eqref{eq:Lipshitz} because the level of nonlinearity may not be reduced no matter how small the region is. 
\end{rem}

With the transformation matrix chosen as the identity $T(t)=I$, the valid set $\mtn(t)$ corresponding to the Lipchitz nonlinearity in \eqref{eq:Lipshitz} is given by 
\begin{align}
\mtn(t) & =\{(N_{1},N_{2}):N_{1}=\blkdiag(\lambda_{1}I_{n_{q_{[1]}}},\ldots,\lambda_{n_{c}}I_{n_{q_{[n_{c}]}}}),\nonumber \\
  N_{2} & =\diag(\lambda_{1}\gamma_{1}^{2}(t),\ldots,\lambda_{n_{c}}\gamma_{n_{c}}^{2}(t)),\;\lambda_{i}\in\mathbb{R}_{+},\;i\in\mathbb{Z}_{[1,n_{c}]}\}.\label{eq:N_Lipschitz}
\end{align}
Similarly, for the L-smooth nonlinearity in \eqref{eq:Lsmooth}, the valid set $\mtn(t)$ is given by
\begin{align}
\mtn(t) & =\{(N_{1},N_{2}):N_{1}=\blkdiag(\lambda_{1}I_{n_{q_{[1]}}},\ldots,\lambda_{n_{c}}I_{n_{q_{[n_{c}]}}}),\nonumber \\
 & N_{2}=\diag(\lambda_{1}s_{1}(t)\beta_{1}^{2}(t),\ldots,\lambda_{n_{c}}s_{n_{c}}(t)\beta_{n_{c}}^{2}(t)),\nonumber \\
 & \norm{\delta q_{[i]}(t)}_{2}^{2}\leq s_{i}(t),\;\lambda_{i}\in\mathbb{R}_{+},\;i\in\mathbb{Z}_{[1,n_{c}]}\},\label{eq:N_Lsmooth}
\end{align}
where $s:[t_{0},t_{f}]\rightarrow\mathbb{R}_{+}^{n_{c}}$ serves as one of the decision variables. To ensure the validity of the multiplier matrix for the L-smooth nonlinearity, each inequality $\norm{\delta q_{[i]}}_{2}^{2}\leq s_{i}$ must be satisfied. Assuming $G_{i}=0$, it follow that $\norm{\delta q_{[i]}}_{2}=\norm{C_{i}\eta(t)+D_{i}\xi(t)}_{2}$. Therefore, the condition $\norm{\delta q_{[i]}}_{2}^{2}\leq s_{i}$ can be equivalently enforced by the following LMI: 
\begin{equation}
\left[\begin{array}{cc}
s_{i}(t)I_{n_{x}} & *\\
(C_{i}Q(t)+D_{i}Y(t))^{\top} & Q(t)
\end{array}\right]\succeq0.\label{eq:sv_LMI}
\end{equation}

\begin{example}
\label{exa:unicycle_2}Recall the Jacobian $\phi'$ of the nonlinearity $\phi$ of the unicycle model (see Example \ref{exa:unicycle}) given in \eqref{eq:phi_hat_unicycle}. Suppose the control input for the velocity is bounded such that $\abs{u_{1}}\leq v_{max}$ for some $v_{max}\in\mathbb{R}_{+}$. Then, the row-wise norms of $\phi'$ can be deduced by
\begin{align*}
\norm{\phi_{[1]}'}_{2} & =\sqrt{u_{1}^{2}\sin^{2}x_{3}^{2}+\cos^{2}x_{3}^{2}}\leq\max\{1,v_{max}\},\\
\norm{\phi_{[2]}'}_{2} & =\sqrt{u_{1}^{2}\cos^{2}x_{3}^{2}+\sin^{2}x_{3}^{2}}\leq\max\{1,v_{max}\}.
\end{align*}
Since the Jacobian is bounded, each $\phi_{[i]}$ is Lipschitz continuous with constant at most $\max\{1,v_{max}\}$, and thus the corresponding constant $\gamma_{i}$ in \eqref{eq:Lipshitz} is $\gamma_{i}=2\max\{1,v_{max}\}$ for $i=1,2$. $\frac{\partial}{\partial q_{[1]}}\phi_{[1]}'(t,q_{[1]})$

In addition, since $\phi_{i}'$ is differentiable, we can compute:
\begin{align*}
\frac{\partial\phi_{[1]}'(t,q_{[1]})}{\partial q_{[1]}}&=\left[\begin{array}{cc}
-u_{1}\cos x_{3} & -\sin x_{3}\\
-\sin x_{3} & 0
\end{array}\right],\\
\frac{\partial\phi_{[2]}'(t,q_{[2]})}{\partial q_{[2]}}&=\left[\begin{array}{cc}
-u_{1}\sin x_{3} & \cos x_{3}\\
\cos x_{3} & 0
\end{array}\right].
\end{align*}
Considering norms, we find
\[
\bigg|\bigg|\frac{\partial\phi_{[1]}'}{\partial q_{[1]}}\bigg|\bigg|_{2}\leq\max\{\sqrt{2},v_{max}\},\;\bigg|\bigg|\frac{\partial\phi_{[2]}'}{\partial q_{[2]}}\bigg|\bigg|_{2}\leq\max\{\sqrt{2},v_{max}\}.
\]
Hence, the valid L-smooth constants $\beta_{i}$ in \eqref{eq:Lsmooth} are $\beta_{i}=\max\{\sqrt{2},v_{max}\}$ for $i=1,2$.
\end{example}

\subsection{Invariance condition by DLMI}

We establish funnel invariance using a dissipation inequality based on a Lyapunov function in the following Lemma.
\begin{lem}
\label{lem:state_funnel_invariance}Let $Q:[t_{0},t_{f}]\rightarrow\mathbb{S}_{++}^{n_{x}}$, $K:[t_{0},t_{f}]\rightarrow\mathbb{R}^{n_{u}\times n_{x}}$, $\lambda_{w}:[t_{0},t_{f}]\rightarrow\mathbb{R}_{+}$ be piecewise continuous functions. Define the Lyapunov function $V:[t_{0},t_{f}]\times\mathbb{R}^{n_{x}}\rightarrow\mathbb{R}_{+}$ as in \eqref{eq:Lyapunov_function}. Suppose that all trajectories of the incremental system \eqref{eq:incremental_dyn_2} with the feedback control $\xi(t)=K(t)\eta(t)$ satisfy, for almost all $t\in[t_{0},t_{f}]$, the following inequalities :
\begin{subequations}
\label{eq:dissipative_inequality}
\begin{align}
\dot{V}(t,\eta)+\alpha V(t,\eta)-\lambda_{w}(t)w(t)^{\top}w(t) & \leq0,\label{eq:Lyapunov_condition}\\
0\leq\lambda_{w}(t) & \leq\alpha,\label{eq:lambda_d}
\end{align}
\end{subequations}
for some decay rate $\alpha\in\mathbb{R}_{++}$. Then, the state funnel $\mte_{\eta}(t)$, defined in \eqref{eq:state_funnel}, with the support value $c_{Q}=w_{max}^{2}$ is invariant. That is, if $\eta(\cdot)$ is a solution of \eqref{eq:incremental_dyn_2} with $\eta(t_{0})\in\mte_{\eta}(t_{0})$, then $\eta(t)\in\mte_{\eta}(t)$ for all $t\in[t_{0},t_{f}]$.
\end{lem}

\begin{proof}
Start with the dissipation inequality \eqref{eq:Lyapunov_condition} multiply both sides by $e^{\alpha(t-t_{0})}$:
\[
e^{\alpha(t-t_{0})}\dot{V}(t,\eta)+\alpha e^{\alpha(t-t_{0})}V(t,\eta)\leq\lambda_{w}(t)e^{\alpha(t-t_{0})}w(t)^{\top}w(t).
\]
Noting that the left-hand side is the derivative of $e^{\alpha(t-t_{0})}V(t,\eta)$, we write
\[
\frac{d}{dt}\left(e^{\alpha(t-t_{0})}V(t,\eta)\right)\leq\lambda_{w}(t)e^{\alpha(t-t_{0})}w(t)^{\top}w(t).
\]
Integrate both sides from $t_{0}$ to $t\in[t_{0},t_{f}]$:
\begin{equation*}
e^{\alpha(t-t_{0})}V(t,\eta)-V(t_{0},\eta) \leq
\int_{t_{0}}^{t}\lambda_{w}(\tau)\,e^{\alpha(\tau-t_{0})}\norm{w(\tau)}_2^2\mathrm{d}\tau.
\label{eq:int_bound}
\end{equation*}
Using the pointwise bound $\norm w_{2}^{2}\leq w_{max}^{2}$ from \eqref{eq:w_pointwise} and $\lambda_{w}(t)\leq\alpha$ from \eqref{eq:lambda_d}, we further bound the integral:
\begin{align*}
\int_{t_{0}}^{t}\lambda_{w}(\tau)e^{\alpha(\tau-t_{0})}\norm{w(\tau)}_2^2\mathrm{d}\tau &\leq \alpha w_{max}^{2}\int_{t_{0}}^{t}e^{\alpha(\tau-t_{0})}\mathrm{d}\tau, \\ & = w_{max}^{2}(e^{\alpha(t-t_{0})}-1).
\end{align*}
Therefore,
\[
e^{\alpha(t-t_{0})}V(t,\eta)\leq V(t_{0},\eta)+w_{max}^{2}(e^{\alpha(t-t_{0})}-1).
\]
Multiply both sides by $e^{-\alpha(t-t_{0})}$:
\begin{align*}
V(t,\eta) & \leq e^{-\alpha(t-t_{0})}V(t_{0},\eta)+(1-e^{-\alpha(t-t_{0})})w_{max}^{2}.
\end{align*}
If $\eta(t_{0})\in\mte_{\eta}(t_{0})$, so $V(t_{0},\eta)\leq w_{max}^{2}$, and we have $V(t,\eta) \leq w_{max}^{2}$. Therefore, we can conclude that $\eta(t)\in\mte_{\eta}(t)$ for all $t\in[t_{0},t_{f}]$.
\end{proof}
The following lemma shows the invariance of the state funnel $\mte_{\eta}$ implies the invariance of the full funnel $\mtf$.
\begin{lem}
\label{lem:funnel_invariance}If the state funnel $\mte_{\eta}(t)$ is invariant for the incremental system \eqref{eq:incremental_dyn_2} under the linear feedback gain $K(t)$, then the funnel $\mtf(t)$, defined in \eqref{eq:funnel}, is invariant for the uncertain nonlinear system \eqref{eq:nonlinear_system} in the sense described in Definition \ref{def:funnel_invariance}. 
\end{lem}

\begin{proof}
Since $x(t_{0})\in\{\bar{x}(t_{0})\}\oplus\mte_{\eta}(t_{0})$, it follows that $\ensuremath{\eta(t_{0})=x(t_{0})-\bar{x}(t_{0})\in\mte_{\eta}(t_{0})}$. By the assumed invariance of the state $\mte_{\eta}(t)$, we have $\eta(t)\in\mte_{\eta}(t)$ for all $t\in[t_{0},t_{f}]$. Under the feedback law $u(t)=\bar{u}(t)+K(t)\eta(t)$, it follows that the input deviation is $\xi(t)=K(t)\eta(t)$, and hence $\xi(t)\in\mte_{\xi}(t)$. Therefore, $(x(t),u(t))\in\mtf(t)$ for all $t\in[t_{0},t_{f}]$.
\end{proof}
Now, we derive a sufficient condition for funnel invariance by formulating the dissipation inequality \eqref{eq:dissipative_inequality} as a DLMI. This condition enables the computation of the funnel-shaping matrices $Q(t)$ and feedback gains $K(t)$ using convex optimization techniques.
\begin{thm}
Consider the incremental system given in \eqref{eq:incremental_dyn_2}. Suppose there exist functions $Q:[t_{0},t_{f}]\rightarrow\mathbb{S}_{++}^{n_{x}}$ , $Y:[t_{0},t_{f}]\rightarrow\mathbb{R}^{n_{u}\times n_{x}}$, $\lambda_{w}:[t_{0},t_{f}]\rightarrow\mathbb{R}_{+}$, $N_{1}:[t_{0},t_{f}]\rightarrow\mathbb{S}^{n_{q}}$, and $N_{2}:[t_{0},t_{f}]\rightarrow\mathbb{S}^{n_{p}}$, such that $(N_{1}(t),N_{2}(t))\in\mtn(t)$, $0\leq\lambda_{w}(t)\leq\alpha$, and the following DLMI holds for almost all $t\in[t_{0},t_{f}]$:
\begin{equation}
\left[\begin{array}{cc}
H_{11}(t)+H_{11}(t)^{\top}-\dot{Q}(t) & \star\\
H_{12}(t) & H_{22}(t)
\end{array}\right]\preceq0,\label{eq:DLMI_main}
\end{equation}
where the matrix blocks are defined as:
\begin{subequations}
\label{eq:DLMI_matrix_blocks}
\begin{align}
H_{11}(t) & =A(t)Q(t)+B(t)Y(t) \nonumber \\ 
-E&T_{22}(t)^{-1}T_{21}(t)(CQ(t)+DY(t))+\frac{\alpha}{2}Q(t),\label{eq:def_H11}\\
H_{12}(t) & =\left[\begin{smallmatrix}
N_{2}(t)T_{22}(t)^{-\top}E^{\top}\\
F(t)^{\top}-D^{\top}T_{21}(t)^{\top}T_{22}(t)^{-T}E^{\top}\\
\Sigma(t)(CQ(t)+DY(t))
\end{smallmatrix}\right]\\
H_{22}(t) & =\left[\begin{smallmatrix}
-N_{2}(t) & \star & \star\\
0 & -\lambda_{w}(t)I & \star\\
\Lambda N_{2} & \Sigma(t)D & -N_{1}(t)
\end{smallmatrix}\right]
\end{align}
\end{subequations}
and the auxiliary matrices are: 
\[
\Lambda(t)\coloneqq T_{12}(t)T_{22}(t)^{-1}, \Sigma(t)\coloneqq T_{11}(t)-T_{12}(t)T(t)_{22}^{-1}T_{21}(t).
\]
Then, with the feedback gain defined as $K(t)=Y(t)Q(t)^{-1}$, then the dissipation inequality \eqref{eq:dissipative_inequality} holds almost all $t\in[t_{0},t_{f}]$ and the funnel $\mathcal{F}(t)$ defined in \eqref{eq:funnel} is invariant in the sense described in Definition \ref{def:funnel_invariance}.
\end{thm}

\begin{proof}
Throughout, the time argument $t$ is often omitted whenever clear. Define the closed-loop matrices $A_{cl}=A+BK$ and $C_{cl}=C+DK$. Then the DLMI \eqref{eq:DLMI_main} becomes:
\begin{align*}
\left[\begin{smallmatrix}
\left[\begin{smallmatrix}
    Q(A_{cl}-ET_{22}^{-1}T_{21}C_{cl})^{\top} \\ + (A_{cl}-ET_{22}^{-1}T_{21}C_{cl})Q-\dot{Q}+\alpha Q
\end{smallmatrix}\right]
 & \star & \star & \star\\
N_{2}T_{22}^{-\top}E^{\top} & -N_{2} & \star & \star\\
F^{\top}-D^{\top}T_{21}^{\top}T_{22}^{-T}E^{\top} & 0 & -\lambda_{w}I & \star\\
\Sigma C_{cl}Q & \Lambda N_{2} & \Sigma D & -N_{1}
\end{smallmatrix}\right] & \preceq0.
\end{align*}
Multiply on the left and right by $\text{diag}(Q^{-1},N_{2}^{-1},I,I)$ and apply Schur complement:
\begin{align*}
\left[\begin{smallmatrix}
\left[\begin{smallmatrix}
    Q(A_{cl}-ET_{22}^{-1}T_{21}C_{cl})^{\top} \\ + (A_{cl}-ET_{22}^{-1}T_{21}C_{cl})Q-\dot{Q}+\alpha Q
\end{smallmatrix}\right]
& \star & \star\\
T_{22}^{-\top}E^{\top} & -N_{2} & \star\\
F^{\top}Q^{-1}-D^{\top}T_{21}^{\top}T_{22}^{-T}E^{\top}Q^{-1} & 0 & -\lambda_{w}I
\end{smallmatrix}\right]+\\
\left[\begin{smallmatrix}
\Sigma C_{cl} & \Lambda & \Sigma D\end{smallmatrix}\right]^{\top}\begin{smallmatrix}N_{1}^{-1}\end{smallmatrix}\left[\begin{smallmatrix}
\Sigma C_{cl} & \Lambda & \Sigma D\end{smallmatrix}\right] & \preceq0.
\end{align*}
More compactly, we can write:
\[
\left[\begin{smallmatrix}
A_{cl}^{\top}Q^{-1}+Q^{-1}A_{cl}-Q^{-1}\dot{Q}Q^{-1}+\alpha Q^{-1}-\Theta & \star & \star\\
T_{22}^{-\top}E^{\top}Q^{-1} & 0 & \star\\
F^{\top}Q^{-1}-D^{\top}T_{21}^{\top}T_{22}^{-T}E^{\top}Q^{-1} & 0 & -\lambda_{w}I
\end{smallmatrix}\right]+\Psi\preceq0,
\]
where
\begin{align*}
\Theta & \coloneqq(ET_{22}^{-1}T_{21}C_{cl})^{\top}Q^{-1}+Q^{-1}(ET_{22}^{-1}T_{21}C_{cl}),\\
\Psi & \coloneqq\left[\begin{smallmatrix}
C_{cl} & 0 & D\\
0 & I & 0
\end{smallmatrix}\right]^{\top}\left[\begin{smallmatrix}
\Sigma & \Lambda\\
0 & I
\end{smallmatrix}\right]^{\top}\left[\begin{smallmatrix}
N_{1}^{-1} & 0\\
0 & -N_{2}^{-1}
\end{smallmatrix}\right]\left[\begin{smallmatrix}
\Sigma & \Lambda\\
0 & I
\end{smallmatrix}\right]\left[\begin{smallmatrix}
C_{cl} & 0 & D\\
0 & I & 0
\end{smallmatrix}\right].
\end{align*}
To eliminate the cross-coupling terms and simplify the structure, define:
\[
T_{post}\coloneqq\left[\begin{smallmatrix}
I & 0 & 0\\
T_{21}C_{cl} & T_{22} & T_{21}D\\
0 & 0 & I
\end{smallmatrix}\right],\quad T_{pre}\coloneqq T_{post}^{\top}.
\]
Post- and pre-multiplying by the above by $T_{pre}$ and $T_{post}$, respectively, yields:
\begin{align*}
\left[\begin{smallmatrix}
A_{cl}^{\top}Q^{-1}+Q^{-1}A_{cl}-Q^{-1}\dot{Q}Q^{-1}+\alpha Q^{-1} & \star & \star\\
E^{\top}Q^{-1} & 0 & \star\\
F^{\top}Q^{-1} & 0 & -\lambda_{w}I
\end{smallmatrix}\right] & +\\
\left[\begin{smallmatrix}
C_{cl} & 0 & D\\
0 & I & 0
\end{smallmatrix}\right]^{\top}\begin{smallmatrix}T^{\top}\end{smallmatrix}\left[\begin{smallmatrix}
N_{1}^{-1} & 0\\
0 & -N_{2}^{-1}
\end{smallmatrix}\right]
\begin{smallmatrix}T\end{smallmatrix}
\left[\begin{smallmatrix}
C_{cl} & 0 & D\\
0 & I & 0
\end{smallmatrix}\right] & \preceq0.
\end{align*}
This matrix inequality implies that for all $\eta\in\mathbb{R}^{n_{x}}$, $\delta p\in\mathbb{R}^{n_{p}}$, and $w\in\mathbb{R}^{n_{w}}$, we have
\[
\dot{V}+\alpha V+\left[\begin{smallmatrix}
\delta q\\
\delta p
\end{smallmatrix}\right]^{\top}T^{\top}\left[\begin{smallmatrix}
N_{1}^{-1} & 0\\
0 & -N_{2}^{-1}
\end{smallmatrix}\right]T\left[\begin{smallmatrix}
\delta q\\
\delta p
\end{smallmatrix}\right]-\lambda_{w}\norm{w}_2^2 \leq0.
\]
It follows from the definition of the $\delta$QC \eqref{eq:def_iQC} that the middle term is nonpositive. Thus, the dissipation inequality \eqref{eq:dissipative_inequality}, $\dot{V}+\alpha V-\lambda_{w}w^{\top}w\leq0$, holds almost all $t\in[t_{0},t_{f}]$. By Lemma \ref{lem:state_funnel_invariance}, this implies the invariance of the state funnel $\mte_{\eta}(t)$, and consequently, the invariance of the full funnel $\mtf(t)$ follows from Lemma \ref{lem:funnel_invariance}.
\end{proof}

\subsection{Funnel constraints}

We develop time-varying LMI conditions implying the constraint satisfaction as defined in Definition \eqref{def:funnel_feasibility}, requiring that all states and control inputs within the funnel respect the state and input constraints. 

Let the feasible sets be defined as:
\begin{subequations}
\label{eq:constraint_sets}
\begin{align}
\mtx & =\bigcap_{i=1}^{m_{x}}\mtx^{i},\quad\mtu=\bigcap_{j=1}^{m_{u}}\mtu^{j}\\
\mtx^{i} & =\{x\mid h_{i}(x)\leq0\},\quad i=\mathbb{Z}_{[1,m_{x}]},\\
\mtu^{j} & =\{u\mid g_{j}(u)\leq0\},\quad j=\mathbb{Z}_{[1,m_{u}]},
\end{align}
\end{subequations}
where each $h_{i}:\mathbb{R}^{n_{x}}\rightarrow\mathbb{R}$ and $g_{j}:\mathbb{R}^{n_{u}}\rightarrow\mathbb{R}$ are assumed to be continuously differentiable. Our goal is to ensure that the state and input funnels remain inside the feasible sets: 
\begin{subequations}
\label{eq:what_funnel_constraint_aim}
\begin{align}
(\bar{x}(t)\oplus\mte_{\eta}(t)) & \subseteq\mtx,\\
(\bar{u}(t)\oplus\mte_{\xi}(t)) & \subseteq\mtu,\quad\forall t\in[t_{0},t_{f}],
\end{align}
\end{subequations}
which can be compactly written as $\mtf(t)\subseteq\mtx\times\mtu$. However, when it comes to general nonlinear functions $h_{i}$ and $g_{j}$, it is not tractable to reformulate the conditions \eqref{eq:what_funnel_constraint_aim} into equivalent LMIs. We therefore linearize $h_{i}$ and $g_{j}$ around the nominal trajectory $(\bar{x},\bar{u})$, to obtain the following time-varying convex polytope approximations:
\begin{subequations}
\label{eq:polyhedral_constraint_set}
\begin{align}
\mtp_{x} & =\bigcap_{i=1}^{m_{x}}\mtp_{x}^{i},\quad\mtp_{u}=\bigcap_{j=1}^{m_{u}}\mtp_{u}^{i}\\
\mtp_{x}^{i}(t) & =\{x\mid(a_{i}^{h}(t))^{\top}x\leq b_{i}^{h}(t)\},\quad i=\mathbb{Z}_{[1,m_{x}]},\\
\mtp_{u}^{u}(t) & =\{u\mid(a_{j}^{g}(t))^{\top}u\leq b_{j}^{g}(t)\},\quad j=\mathbb{Z}_{[1,m_{u}]},
\end{align}
\end{subequations}
where
\begin{align*}
a_{i}^{h}(t) & =\frac{\partial h_{i}}{\partial x}\bigg|_{x=\bar{x}(t)},\quad b_{i}^{h}(t)=a_{i}^{h}(t)^{\top}\bar{x}(t)-h_{i}(\bar{x}),\\
a_{j}^{g}(t) & =\frac{\partial g_{j}}{\partial u}\bigg|_{u=\bar{u}(t)},\quad b_{j}^{g}(t)=a_{i}^{g}(t)^{\top}\bar{u}(t)-g_{j}(\bar{u}).
\end{align*}

We now provide LMI conditions that ensures the funnel remains within the linearized state and input constraint sets \eqref{eq:polyhedral_constraint_set}.
\begin{lem}
Suppose $\bar{x}(t)\in\mtp_{x}(t)$ and $\bar{u}(t)\in\mtp_{u}(t)$ for all $t\in[t_{0},t_{f}]$. Then, the inclusion $\mtf(t)\subseteq\mtp_{x}(t)\times\mtp_{u}(t)$ is equivalent to satisfying the following LMIs for all $t\in[t_{0},t_{f}]$, $i=\mathbb{Z}_{[1,m_{x}]}$, and $j=\mathbb{Z}_{[1,m_{u}]}$:
\begin{subequations}
\label{eq:LMIs_for_constraints}
\begin{align}
0 & \preceq\left[\begin{smallmatrix}
(b_{i}^{h}(t)-a_{i}^{h}(t){}^{\top}\bar{x}(t))^{2} & c_{Q}a_{i}^{h}(t)^{\top}Q(t)\\
c_{Q}Q(t)a_{i}^{h}(t) & c_{Q}Q(t)
\end{smallmatrix}\right],\label{eq:LMI_state_constraint}\\
0 & \preceq\left[\begin{smallmatrix}
(b_{j}^{g}(t)-a_{j}^{g}(t){}^{\top}\bar{u}(t))^{2} & c_{Q}a_{i}^{g}(t)^{\top}Y(t)\\
c_{Q}Y(t)^{\top}a_{i}^{h}(t) & c_{Q}Q(t)
\end{smallmatrix}\right].\label{eq:LMI_input_constraint}
\end{align}
\end{subequations}
\end{lem}

\begin{proof}
We derive the LMI \eqref{eq:LMI_input_constraint} for the input constraint; the state constraint follows analogously. The inclusion $\bar{u}\oplus\mte_{\xi}\subset\mtp_{u}$ means that for all $\xi\in\mte_{\epsilon}$, we must have $(a_{j}^{g})^{\top}(\bar{u}+\xi)\leq b_{j}^{g}$. This is equivalent to
\[
\max_{\xi\in\mte_{\xi}}(a_{j}^{g})^{\top}(\bar{u}+\xi)\leq b_{j}^{g}.
\]
Using the definition of the input funnel $\mte_{\xi}$ in \eqref{eq:input_funnel}, we can write
\[
\max_{\norm y_{2}\leq1}(a_{j}^{g})^{\top}(c_{Q}KQK^{\top})^{\frac{1}{2}}y\leq b_{j}^{g}-(a_{j}^{g})^{\top}\bar{u},
\]
where the term related to $\bar{u}$ is moved to the right-hand side since $\bar{u}$ is fixed. The solution of the maximization problem is equal to the norm $\norm{(c_{Q}KQK^{\top})^{\frac{1}{2}}a_{j}^{g}}_{2}$. Hence, we obtain: 
\[
\norm{(c_{Q}KQK^{\top})^{\frac{1}{2}}a_{j}^{g}}_{2}\leq b_{j}^{g}-(a_{j}^{g})^{\top}\bar{u}.
\]
Since the left-hand side is always nonnegative, squaring both sides preserves the inequality:
\[
c_{Q}(a_{j}^{g})^{\top}KQK^{\top}a_{j}^{g}\leq(b_{j}^{g}-(a_{j}^{g})^{\top}\bar{u})^{2}.
\]
Applying Schur complement generates:
\[
\left[\begin{array}{cc}
(b-(a_{j}^{g})^{\top}\bar{u})^{2} & \sqrt{c_{Q}}a^{\top}K\\
\sqrt{c_{Q}}K^{\top}a & Q^{-1}
\end{array}\right]\succeq0
\]

Multiplying both sides by $\diag(1,\sqrt{c_{Q}}Q)$ and using $Y=KQ$ completes the proof.
\end{proof}
The inclusion $\mtp_{x}^{i}(t)\subseteq\mtx^{i}$ (or $\mtp_{u}^{j}(t)\subset\mtu^{j}$) is not guaranteed for general nonlinear functions $h_{i}$ (or $g_{j}$). However, this inclusion holds under certain structural assumptions. Specifically, when the nonlinear functions are concave, the linearized sets are inner assumptions.
\begin{lem}
Let $h_i:\mathbb{R}^{n_x}\!\to\!\mathbb{R}$ be twice differentiable and let $\bar x\in\mathbb{R}^{n_x}$.
Assume there exists a symmetric matrix $P\succeq 0$ such that
\begin{equation}
\label{eq:P_dominates_all_H}
P-\tfrac{1}{2}\nabla^2 h_i(\tilde x)\succeq 0 \quad \text{for all } \tilde x\in\mathbb{R}^{n_x}.
\end{equation}
Define $a_i^h \triangleq \nabla h_i(\bar x)$. Then, for all $x\in\mathbb{R}^{n_x}$,
\begin{equation}
\label{eq:global_upper_bound}
h_i(x)\;\le\; h_i(\bar x) + (a_i^h)^\top(x-\bar x) \;+\; (x-\bar x)^\top P (x-\bar x).
\end{equation}
\end{lem}

\begin{proof}
By Taylor’s theorem at $\bar x$, there exists $\tilde x$ on the segment between $\bar x$ and $x$ such that
\[
h_i(x)=h_i(\bar x)+(a_i^h)^\top(x-\bar x)+\tfrac{1}{2}(x-\bar x)^\top\nabla^2 h_i(\tilde x)(x-\bar x).
\]
By \eqref{eq:P_dominates_all_H}, the last term satisfies 
$\tfrac{1}{2}(x-\bar x)^\top\nabla^2 h_i(\tilde x)(x-\bar x)\le (x-\bar x)^\top P(x-\bar x)$,
which yields \eqref{eq:global_upper_bound}.
\end{proof}

\begin{corr}
\label{cor:polytope_inclusion}
Suppose $\bar x(t)\in\mathcal X^{i}$ where $\mathcal X^{i}$ is defined in \eqref{eq:constraint_sets}.
If, in addition, $h_i$ is concave so that $\nabla^2 h_i(\tilde x)\preceq 0$ for all $\tilde x$, then taking $P=0$ in \eqref{eq:P_dominates_all_H} yields
\[
h_i(x)\;\le\; h_i(\bar x(t)) + \nabla h_i(\bar x(t))^\top\!\big(x-\bar x(t)\big).
\]
Therefore, any $x$ satisfying the polyhedral constraint \eqref{eq:polyhedral_constraint_set},
\[
\big(\nabla h_i(\bar x(t))\big)^\top x \;\le\; b_i^h \;=\; \big(\nabla h_i(\bar x(t))\big)^\top \bar x(t) - h_i(\bar x(t)),
\]
also satisfies $h_i(x)\le 0$, i.e., $\mathcal P_x^{i}(t)\subseteq \mathcal X^{i}$.
\end{corr}


\begin{rem}
If $h_{i}$ is convex instead of concave, the linearized set $\mtp_{x}^{i}$ forms an outer approximation of $\mtx^{i}$. In this case, instead of using $\mtp_{x}^{i}$ derived in \eqref{eq:polyhedral_constraint_set}, techniques from polyhedral inner approximation, such as those discussed in \cite{guessab2014approximations}, can be used to construct a conservative polytope such that $\mtp\subset\mtx^{i}$. Then, one can enforce the constraint using LMIs similar to those in \eqref{eq:LMIs_for_constraints}.
\end{rem}

\subsection{Funnel cost functions}

A typical objective in funnel design is to maximize the size of the funnel entry to enlarge the set of initial states from which the system can be safely controlled. The proposed method does not restrict the choice of the cost function, as long as it remains convex with respect to the decision variables. In this subsection, we outline several convex cost functions that can be used to promote larger funnel entries while maintaining computational tractability.

We consider cost functions $J:\mathbb{S}_{++}^{n_{x}}\rightarrow\mathbb{R}$ that depend on the ellipsoidal shape of the funnel entry $Q(t_{0})$. To reflect directional preference or scale differences in the state space, we define the weighted matrix $Q_{w0}\coloneqq W_{0}Q(t_{0})W_{0}$, where $W_{0}\in\mathbb{S}_{++}^{n_{x}}$ is a diagonal weighting matrix. The examples of convex cost functions include:
\begin{equation}
J(Q(t_{0}))\in\{-\trace{Q_{w0}},-\log\det(Q_{w0}),\trace{Q_{w0}^{-1}}\}.\label{eq:funnel_cost}
\end{equation}
The first cost function promotes larger funnel entries by maximizing the sum of lengths of weighted principal axes. The second one is proportional to maximizing the volume of the weighted funnel entry \cite{boyd2004convex}. The third one is to minimize the sum of the eigenvalues of $Q_{w0}^{-1}$, which encourages increasing the eigenvalues of $Q_{w0}$, thereby enlarging the funnel entry.

\subsection{Funnel synthesis problem}

We are now ready to formulate the funnel synthesis problem as a continuous-time optimization problem.
\begin{problem}
Continuous-time funnel synthesis.
\begin{subequations}
\label{eq:cont_funl}
\begin{align}
    \underset{\resizebox{0.3\hsize}{!}{$Q(\cdot),Y(\cdot),N_{1}(\cdot),N_{2}(\cdot),\lambda_{d}(\cdot)$}}{\operatorname{minimize}}~~& J(Q(t_{0})) \\[-0.1cm]
    \operatorname{subject~to}~~~~~&  \forall\, t\in[t_0,t_f], \\
 & Q(t)\succ0,\label{eq:cont_PD}\\
 & \eqref{eq:DLMI_main}, \eqref{eq:LMIs_for_constraints}, \\
 & (N_{1}(t),N_{2}(t))\in\mtn(t),\label{eq:cont_funl_multiplier_matrices}\\
 & Q(t_{f})\preceq Q_{f}\label{eq:cont_funl_boundary}. 
\end{align}
\end{subequations}

The constraint \eqref{eq:cont_PD} enforces positive definiteness of the matrix $Q(t)$. The dissipation inequality implied by the DLMI \eqref{eq:DLMI_main} imposes the funnel invariance condition. The pointwise LMIs \eqref{eq:LMIs_for_constraints} ensure the funnel feasibility. The inclusion \eqref{eq:cont_funl_multiplier_matrices} enforces valid multiplier matrices for the nonlinearity. Finally, the terminal condition \eqref{eq:cont_funl_boundary} with $Q_{f}\in\mathbb{S}_{++}^{n_{x}}$ limits the size of the funnel at the final time.
\end{problem}

\begin{rem}
\label{rem:alpha}The resulting optimization problem is convex in continuous time, except for the decay rate $\alpha$. Since $\alpha$ is a scalar, it can be efficiently selected via line search. In the special case where the disturbance $w$ is absent, setting $\alpha=0$ is sufficient for the funnel's invariance and leads to a fully convex formulation.
\end{rem}

In the rest of this subsection, we specify how the inclusion \eqref{eq:cont_funl_multiplier_matrices} can be converted into pointwise (i.e., at each time $t$) linear constraints or LMIs. Consider the sector bounded nonlinearity illustrated in \eqref{eq:multiplier_sectorbounded}. Then, the inclusion $(N_{1},N_{2})\in\mtn$ is equivalent to the following linear equality constraints:
\begin{equation}
N_{1}(t)=\lambda_{sec}(t)S,\quad N_{2}(t)=\lambda_{sec}(t)S,\label{eq:inclusion_sector_bound}
\end{equation}
for any positive scalar-valued function $\lambda_{sec}:[t_{0},t_{f}]\rightarrow\mathbb{R}_{+}$. In this case, the function $\lambda_{sec}(t)$ can be treated as an our decision variable with the constraint on $\lambda_{sec}(t)\geq0$, while preserving the convexity of our funnel synthesis problem \eqref{eq:cont_funl}.

For the Lipschitz nonlinearity \eqref{eq:Lipshitz} and its associated multiplier matrix \eqref{eq:N_Lipschitz}, the inclusion \eqref{eq:cont_funl_multiplier_matrices} is equivalent to the following linear equality constraints:
\begin{subequations}
\label{eq:inclusion_Lipschitz}
\begin{align}
N_{1}(t) & =\blkdiag(\lambda_{1}^{\gamma}(t)I_{n_{q_{[1]}}},\ldots,\lambda_{n_{c}}^{\gamma}(t)I_{n_{q_{[n_{c}]}}}),\\
N_{2}(t) & =\diag(\lambda_{1}^{\gamma}(t)\gamma_{1}^{2}(t),\ldots,\lambda_{n_{c}}^{\gamma}(t)\gamma_{n_{c}}^{2}(t)),
\end{align}
\end{subequations}
where $\lambda_{i}^{\gamma}:[t_{0},t_{f}]\rightarrow\mathbb{R}_{+}$ are positive-valued functions that can be optimized as decision variables with the constraint $\lambda_{i}^{\gamma}(t)\geq0$.

On the other hand, for the L-smooth nonlinearity \eqref{eq:Lsmooth} and its multiplier matrix \eqref{eq:N_Lsmooth}, the inclusion \eqref{eq:cont_funl_multiplier_matrices} can be converted into the following two linear equality and one LMI constraint:
\begin{subequations}
\label{eq:inclusion_Lsmooth}
\begin{align}
N_{1}(t) & =\blkdiag(\lambda_{1}^{\beta}I_{n_{q_{[1]}}},\ldots,\lambda_{n_{c}}^{\beta}I_{n_{q_{[n_{c}]}}}),\\
N_{2}(t) & =\diag(\lambda_{1}^{\beta}s_{1}(t)\gamma_{1}^{2}(t),\ldots,\lambda_{n_{c}}^{\beta}s_{n_{c}}(t)\gamma_{n_{c}}^{2}(t)),\\
 & \left[\begin{smallmatrix}
s_{i}(t)I_{n_{x}} & *\\
(C_{i}Q(t)+D_{i}Y(t))^{\top} & Q(t)
\end{smallmatrix}\right]\succeq0,\quad\forall i\in\mathbb{Z}_{[1,n_{c}]}.\label{eq:inclusion_Lsmooth_LMI}
\end{align}
\end{subequations}
where $\lambda_{i}^{\beta}\in\mathbb{R}_{+}$ are positive constants and $s_{i}:[t_{0},t_{f}]\rightarrow\mathbb{R}_{+}$ are scalar-valued functions that can be jointly optimized as decision variables.

\section{Solution methods \label{sec:4}}

What distinguishes the funnel synthesis problem \eqref{eq:cont_funl} from conventional LMI-based robust controller synthesis \cite{accikmecse2011robust} is that it is formulated as a continuous-time optimization problem, thereby involving DLMI and pointwise LMIs rather than time-invariant LMIs. As a result, it cannot be directly solved using off-the-shelf SDP solvers. To enable tractable computation, a dedicated solution method is required to reformulate the problem into a finite dimensional SDP, and this reformulation must approximate the original problem as closely as possible, ideally without introducing conservatism.

Our solution methods are based on numerical optimal control techniques. To reformulate the funnel synthesis problem \eqref{eq:cont_funl} as an equivalent optimal control problem, we first introduce the notion of funnel dynamics that plays a role analogous to the system dynamics in standard optimal control problem.

\subsection{Funnel Dynamics}

To derive funnel dynamics, we rewrite the DLMI condition \eqref{eq:DLMI_matrix_blocks} as an equivalent differential matrix equation (DME) that serves as the funnel dynamics along with a pointwise LMI constraint. We consider two types of funnel dynamics: Lyapunov-type and $\fdyn$-type, allowing us the flexibility to choose between them. 

The DLMI condition in \eqref{eq:DLMI_matrix_blocks} can be equivalently written as
\begin{subequations}\label{eq:Lyapunov_type_funldyn_H}
\begin{align}
\dot{Q}(t)=H_{11}(t)+H_{11}(t)^{\top} & +Z_{1}(t),\label{eq:Lyapunov_type_funldyn}\\
\left[\begin{array}{cc}
-Z_{1}(t) & \star\\
H_{12}(t) & H_{22}(t)
\end{array}\right] & \preceq0,\label{eq:Lyapunov_type_H}
\end{align}
\end{subequations}
where a slack symmetric matrix-valued function $Z_{1}:[t_{0},t_{f}]\rightarrow\mathbb{S}^{n_{x}}$ is introduced for this equivalent conversion. We refer to \eqref{eq:Lyapunov_type_funldyn} as the Lyapunov-type funnel dynamics as it structurally resembles Lyapunov differential equations \cite{behr2019solution}. 

Alternatively, the DLMI condition \eqref{eq:DLMI_matrix_blocks} can be rewritten as
\begin{subequations}
\label{eq:Basic_type_funldyn_H}
\begin{align}
\dot{Q}(t)=Z_{2}(t),\label{eq:Basic_type_funldyn}\\
\left[\begin{array}{cc}
H_{11}(t)+H_{11}(t)^{\top}-Z_{2}(t) & \star\\
H_{12}(t) & H_{22}(t)
\end{array}\right] & \preceq0,\label{eq:Basic_type_H}
\end{align}
\end{subequations}
where $Z_{2}:[t_{0},t_{f}]\rightarrow\mathbb{S}^{n_{x}}$ is a symmetric slack matrix-valued function. We refer to this formulation \eqref{eq:Basic_type_funldyn_H} as the $\fdyn$-type funnel dynamics as the slack variable $Z_{2}$ directly defines the evolution of $Q$.
\begin{rem}
It is important to note that both formulations are mathematically equivalent. However, once control parameterization is applied after discretization, the resulting solution space differs depending on the chosen funnel dynamics. For instance, piecewise constant parameterization of $Z_{2}(t)$ yields a piecewise linear $Q(t)$, whereas the same applied to $Z_{1}(t)$ does not, as $Q(t)$ is implicitly defined as a solution of differential equation \eqref{eq:Lyapunov_type_funldyn}. 
\end{rem}

Both formulations, \eqref{eq:Lyapunov_type_funldyn} and \eqref{eq:Basic_type_funldyn}, define LTV systems where $Q$ serves as the state and all other decision variables including $Y$ and $Z$ (which corresponds to $Z_{1}$in \eqref{eq:Lyapunov_type_funldyn_H} or $Z_{2}$ in \eqref{eq:Basic_type_funldyn_H}) acts as control inputs. Consequently, the block-form LMIs \eqref{eq:Lyapunov_type_H} and \eqref{eq:Basic_type_H} naturally represent state-input constraints.

To simplify subsequent analysis, we introduce vectorized variables as follows:
\[
q_{v}\coloneqq\vecc Q,\quad y_{v}=\vecc Y,\quad z_{v}=\vecc Z.
\]
Using these variables, both forms of the funnel dynamics, \eqref{eq:Lyapunov_type_funldyn} and \eqref{eq:Basic_type_funldyn}, can be expressed in the vectorized form:
\begin{equation}
\dot{q}_{v}(t)=\underbrace{A_{q}(t)q_{v}(t)+B_{qy}(t)y_{v}(t)+B_{qz}(t)z_{v}(t)}_{\coloneqq F(t,q_{v},y_{v},z_{v})},\label{eq:vectorized_funl_dynamics}
\end{equation}
where $A_{q}(t)\in\mathbb{R}^{n_{q}\times n_{q}}$, $B_{qy}(t)\in\mathbb{R}^{n_{q}\times n_{y}}$, and $B_{qz}(t)\in\mathbb{R}^{n_{q}\times n_{q}}$ are time-varying matrices that can be constructed using Kronecker products \cite{brewer2003kronecker,10167750} with dimensions $n_{q}=n_{x}^{2}$, $n_{y}=n_{x}n_{u}$, and $n_{z}=n_{q}$. We also define the vectorized form of the multiplier matrices $N_{1}$ and $N_{2}$ as:
\[
m_{1,v}\coloneqq\vecc{N_{1}},\quad m_{2,v}=\vecc{N_{2}}.
\]

\subsection{Discretization}

We consider a uniform time grid defined as
\begin{subequations}
\label{eq:uniform_time_mesh}
\begin{align}
t_{k} & =t_{0}+\frac{k}{N}(t_{f}-t_{0}),\quad k\in\mathbb{Z}_{[0,N]},\\
\triangle_{k} & =\triangle(t_{k}),
\end{align}
\end{subequations}
where $N\in\mathbb{Z}_{++}$. The symbol $\triangle$ serves as a placeholder for any time-varying variable, and we denote its value at $t=t_{k}$ by $\triangle_{k}\coloneqq\triangle(t_{k})$. Each $\triangle_{k}$ is referred to as a node point. With the uniform time grid \eqref{eq:uniform_time_mesh}, we use the first-order hold (FOH) interpolation for the control inputs in the funnel dynamics, defined by
\begin{subequations}
\label{eq:FOH_interpolation}
\begin{align}
\square(t) & =\lambda_{k}^{m}(t)\square_{k}^{m}+\lambda_{k}^{p}(t)\square_{k+1}^{p},\quad\forall t\in[t_{k},t_{k+1}),\\
\lambda_{k}^{m}(t) & =\frac{t_{k+1}-t}{t_{k+1}-t_{k}},\quad\lambda_{k}^{p}(t)=\frac{t-t_{k}}{t_{k+1}-t_{k}},
\end{align}
\end{subequations}
for all $k\in\mathbb{Z}_{[0,N)}$. The symbol $\square$ is a placeholder for any funnel input variables (e.g. $y_{v}$, $z_{v}$, $m_{v}$) for the funnel dynamics. To employ the zeroth-order hold (ZOH), it suffices to set $\square_{k}^{m}=\square_{k+1}^{p}$ for all $k\in\mathbb{Z}_{[0,N-1]}$. Also, continuity at each node point $t=t_{k}$ can be enforced by $\square_{k}^{p}=\square_{k}^{m}$ for all $k\in\mathbb{Z}_{[1,N-1]}$.

We adopt a multiple-shooting scheme \cite{bock1984multiple} to enforce the funnel dynamics. Over each subinterval $[t_{k,}t_{k+1})$, the funnel state evolves according to
\[
q_{v}(t)=q_{v}(t_{k})+\int_{t_{k}}^{t}F(\tau,q_{v}(\tau),y_{v}(\tau),z_{v}(\tau))\mathrm{d}\tau,
\]
where $F$ is defined in \eqref{eq:vectorized_funl_dynamics}. Enforcing continuity across subintervals yields the funnel dynamics constraint:
\begin{equation}
q_{v,k+1}=q_{v}(t_{k})+\int_{t_{k}}^{t_{k+1}}F(\tau,q_{v}(\tau),y_{v}(\tau),z_{v}(\tau))\mathrm{d}\tau.\label{eq:const_on_funl_dynamics_integration}
\end{equation}
Since $F$ is linear in its arguments except $t$, and FOH interpolation is used, the right-hand side of \eqref{eq:const_on_funl_dynamics_integration} is affine in the decision variables $q_{v,k},y_{v,k}^{m},y_{v,k+1}^{p},z_{v,k}^{m}$, and $z_{v,k+1}^{p}$. Thus, we can exactly rewrite the constraint \eqref{eq:const_on_funl_dynamics_integration} as
\begin{equation}
q_{v,k+1}=A_{k}^{q}q_{v,k}+B_{yk}^{m}y_{v,k}^{m}+B_{yk}^{p}y_{v,k+1}^{p}+B_{zk}^{m}z_{v,k}^{m}+B_{zk}^{p}z_{v,k+1}^{p},\label{eq:discrete_const_on_funl_dynamics}
\end{equation}
where $A_{k}^{q}$,$B_{yk}^{m/p}$, and $B_{zk}^{m/p}$ are appropriately sized matrices that can be obtained via variational discretization \cite{lin2014control}.

\subsection{Positive Definiteness of $Q$ \label{subsec:Positive-Definiteness-of-Q}}

To ensure the validity of the Lyapunov function \eqref{eq:Lyapunov_function} and the associated funnel \eqref{eq:state_funnel}, it is essential to maintain $Q(t)\succ0$ for all $t\in[t_{0},t_{f}]$. Since $q_{v,k}$ is our decision variable, we can enforce $Q(t)\succ0$ for each $t=t_{k}$ for all $k\in\mathbb{Z}_{[0,N]}$. However, the positive definiteness of $Q(t)$ within the subinterval $(t_{k},t_{k+1})$ depends on the funnel dynamics and the corresponding control inputs, which may lead to inter-node constraint violations. Here, we present cases where the funnel state $Q$ preserves its positive definiteness for all $t$ by only imposting $Q(t)\succ0$ at each node point.
\begin{lem}
\label{lem:PD_of_Q_with_direct}Consider the $\fdyn$-type funnel dynamics \eqref{eq:Basic_type_funldyn}, and let the variable $Z(t)$ (or its vectorized form $z_{v}(t)$) be parameterized by the FOH scheme \eqref{eq:FOH_interpolation} such that 
\begin{equation}
Z_{k}^{m}\succeq Z_{k+1}^{p}\label{eq:Z_Z}
\end{equation}
 for all $k\in\mathbb{Z}_{[0,N-1]}$. Then, imposing $Q(t)\succ0$ for each $t=t_{k}$ for all $k\in\mathbb{Z}_{[0,N]}$ ensures the solution of \eqref{eq:Basic_type_funldyn} preserves the positive definiteness, that is, $Q(t)\succ0$ for all $t\in[t_{0},t_{f}]$.
\end{lem}

\begin{proof}
Let the step size be $d_{k}=t_{k+1}-t_{k}>0$ and introduce the local time $\theta_{t}=t-t_{k}$. With this variable, the FOH expression becomes
\[
Z(\theta_{t})=Z_{k}^{m}+Z_{k}^{r}\theta_{t},
\]
where
\[
Z_{k}^{r}=\frac{Z_{k+1}^{p}-Z_{k}^{m}}{dt_{k}}.
\]
Note that by assumption $Z_{k}^{m}\succeq Z_{k+1}^{p}$, so $Z_{k}^{r}\preceq0$. Integrating the $\fdyn$-funnel dynamics $\dot{Q}(t)=Z(t)$ from $t_{k}$ to $t$ give
\begin{align}
Q(t_{k}+\theta_{t}) & =Q(t_{k})+\theta_{t}Z_{k}^{m}+\frac{1}{2}\theta_{t}^{2}Z_{k}^{r}.\label{eq:Q_t}
\end{align}
Setting $t=t_{k+1}$ (so, $\theta_{t}=d_{k}$) reproduces the node value
\begin{equation}
Q(t_{k+1})=Q(t_{k})+d_{k}Z_{k}^{m}+\frac{1}{2}d_{k}^{2}Z_{k}^{2}\succ0.\label{eq:Q_kp1}
\end{equation}
Define a convex combination $\tilde{Q}(t)$ of $Q(t_{k})$ and $Q(t_{k+1})$ for $t\in[t_{k},t_{k+1}]$:

\begin{equation}
\tilde{Q}(t)\coloneqq(1-\frac{\theta_{t}}{d_{k}})Q(t_{k})+\frac{\theta_{t}}{d_{k}}Q(t_{k+1}).\label{eq:Q_tilde}
\end{equation}
Since $\tilde{Q}(t)$ is the convex combination of two PD matrices $Q(t_{k})$ and $Q(t_{k+1})$, $\tilde{Q}(t)\succ0$ for all $t\in[t_{k},t_{k+1}]$. Substituting \eqref{eq:Q_kp1} into $\tilde{Q}$ generates
\[
\tilde{Q}(t)=Q(t_{k})+\theta_{t}Z_{k}^{m}+\frac{1}{2}\theta_{t}d_{k}Z_{k}^{r}.
\]
Comparing this with \eqref{eq:Q_t} and using $0\leq\theta_{t}\leq d_{k}$ together with $Z_{k}^{2}\preceq0$ gives $\tilde{Q}(t)\preceq Q(t)$. Hence, $Q(t)\succ0$ for all $t\in[t_{k},t_{k+1}]$. We can apply this argument for all $k\in\mathbb{Z}_{[0,N-1]}$, so $Q(t)\succ0$ for all $t\in[t_{0},t_{f}]$.
\end{proof}
\begin{rem}
\label{rem:riccati} It is known that the solution of a differential Riccati equation preserves the positive definiteness of $Q$ \cite{dieci1994positive}. However, the Lyapunov-type funnel dynamics \eqref{eq:Lyapunov_type_funldyn} are not in exact Riccati form, and therefore do not guarantee $Q(t)\succ0$ over the entire interval. One possible way to leverage the positivity-preserving property of Riccati-type dynamics is to treat the feedback gain $K$ as a decision variable instead of $Y$. This, however, introduces bilinear terms between $K$ and $Q$, leading to a nonconvex formulation. To maintain convexity, we retain the use of Lyapunov-type dynamics with $Y$ as the decision variable, and instead address potential inter-sample constraint violations directly in the next section.
\end{rem}

\subsection{Nodal constraint satisfaction}

In addition to the funnel dynamics, we must satisfy the invariance condition \eqref{eq:Lyapunov_type_H} or \eqref{eq:Basic_type_H}, along with the state and input constraints \eqref{eq:LMIs_for_constraints} for all $t\in[t_{0},t_{f}]$. For notational simplicity, we collectively express these conditions in the unified form:
\begin{equation}
L_{l}(t,q_{v}(t),\zeta_{v}(t))\preceq0,\quad i\in\mathbb{Z}_{[1,m_{L}]},\label{eq:continuous-time-LMI}
\end{equation}
where $\zeta_{v}\coloneqq\{y_{v},z_{v},n_{v}\}$ collects decision variables other than $q$, and the indices $l$ are partitioned as follows:
\begin{itemize}
\item $l\in\mathbb{Z}_{[1,m_{inv}]}$: invariance condition, \eqref{eq:Lyapunov_type_H}or \eqref{eq:Basic_type_H} where $m_{inv}=1$,
\item $l\in\mathbb{Z}_{[m_{inv}+1,m_{l}]}$: state constraints \eqref{eq:LMI_state_constraint} and input constraints \eqref{eq:LMI_input_constraint} where $m_{l}=m_{inv}+m_{x}+m_{u}$.
\end{itemize}
Enforcing \eqref{eq:continuous-time-LMI} continuously over the entire horizon $[t_{0},t_{f}]$ is nontrivial, as decision variables are only defined at discrete node points $t_{k}$. A practical approximation is to enforce constraints only at the nodes. Specifically, the nodal constraint satisfaction requires:
\begin{subequations}
\label{eq:nodewise_LMI}
\begin{align}
L_{l}(t_{k},q_{v,k},\zeta_{v,k}^{m}) & \preceq0,\label{eq:nodewise_left}\\
L_{l}(t_{k+1},q_{v,k+1},\zeta_{v,k+1}^{p}) & \preceq0,\label{eq:nodewise_right}\\
\forall k\in\mathbb{Z}_{[0,N-1]},\quad\forall l= & \mathbb{Z}_{[1,m_{l}]}.\nonumber 
\end{align}
\end{subequations}

With the first-order interpolation \eqref{eq:FOH_interpolation}, we allow discontinuities of $\zeta_{v}(t)$ at node points (except at the boundaries), so the constraints are enforced at both left and right evaluations of each node, given by \eqref{eq:nodewise_left} and \eqref{eq:nodewise_right}, respectively.

\subsection{Discrete funnel synthesis problem}

We now derive a finite-dimensional optimization problem based on the discretization of the continuous-time funnel synthesis problem \eqref{eq:cont_funl} over the uniform time grid \eqref{eq:uniform_time_mesh}. This discretized formulation is referred to as the discrete funnel synthesis problem.
\begin{problem}
Discrete funnel synthesis.
\begin{subequations}
\label{eq:funnel_synthesis_discrete}
\begin{align}
    \underset{\resizebox{0.2\hsize}{!}{$q_{v,k},\zeta_{v,k}^{m},\zeta_{v,k}^{p}$}}{\operatorname{minimize}}~~& J(Q_{0}) \\[-0.1cm]
    \operatorname{subject~to}~~~& \eqref{eq:discrete_const_on_funl_dynamics}, \ \eqref{eq:nodewise_LMI}, \ \forall k\in\mathbb{Z}_{[0,N-1]},\\
 & (\mat{(m_{1,v})_{k}},\mat{(m_{2,v})_{k}})\in\mtn(t_{k}),\label{eq:discrete_funl_multiplier_matrices} \nonumber \\
 & \qquad \qquad \qquad \qquad \forall k\in\mathbb{Z}_{[0,N]}, \\
 & \mat{q_{N}}\preceq Q_{f}.\label{eq:discrete_terminal}
\end{align}
\end{subequations}
\end{problem}

The continuous-time set multiplier set inclusion \eqref{eq:cont_funl_multiplier_matrices} is imposed at each node point $t_{k}$, leading to the pointwise inclusion condition in \eqref{eq:discrete_funl_multiplier_matrices}. This nodewise enforcement is sufficient under the following conditions. For sector bounded nonlinearity \eqref{eq:multiplier_sectorbounded}, if the inclusion constraint \eqref{eq:inclusion_sector_bound} hold at each node, it holds for all $t\in[t_{0},t_{f}]$ by convexity and FOH interpolation. For Lipschitz nonlinearity \eqref{eq:N_Lipschitz}, if we define the conservative Lipschitz constant on each subinterval as
\[
\gamma_{k}\coloneqq\max_{t\in[t_{k},t_{k+1}]}\gamma(t),
\]
then enforcing the corresponding multiplier condition \eqref{eq:inclusion_Lipschitz} only at the node points ensures constraint satisfaction throughout the entire interval. For L-smooth nonlinearity \eqref{eq:N_Lsmooth}, a similar bounding can be applied by defining conservative value
\[
\beta_{k}\coloneqq\max_{t\in[t_{k},t_{k+1}]}\beta(t).
\]
However, unlike the sector-bounded and Lipschitz cases, the time-varying LMI \eqref{eq:inclusion_Lsmooth_LMI} associated with L-smooth nonlinearity must be enforced over the entire interval, not just at the nodes. In our formulation, we address this by including an additional LMI condition in $L_{l}$ as described in \eqref{eq:continuous-time-LMI}, ensuring that the L-smooth constraints are enforced properly.

The discrete funnel synthesis problem \eqref{eq:funnel_synthesis_discrete} is a finite-dimensional SDP, which can be efficiently solved using standard solvers such as MOSEK \cite{andersen2000mosek} or Clarabel \cite{goulart2024clarabel}.

\section{Continuous-time constraint satisfaction \label{sec:5}}

In this section, we present two techniques to enhance the continuous-time satisfaction of time-varying LMI constraints \eqref{eq:continuous-time-LMI} beyond enforcement only at discrete node points.

\subsection{Constraint reformulation}

In numerical optimal control, continuous-time inequality constraints (e.g., $g(x(t))\leq0$) are often difficult to enforce directly at all times. One remedy is constraint reformulation \cite{lin2014control,ELANGO2025112464} , where violations over a time interval are quantified and penalized or constrained via an integrated form. 

We adopt a similar idea to our continuous-time LMI constraints. Specifically, noting that the matrix inequality $L_{l}(t,q_{v}(t),\zeta_{v}(t))\preceq0$ is equivalent to the scalar inequality $\lambda_{max}(L_{l}(t,q_{v}(t),\zeta_{v}(t))\leq0$, the satisfaction of \eqref{eq:continuous-time-LMI} over the entire subinterval $t\in[t_{k},t_{k+1}]$ is equivalent to the integral condition:
\begin{align}
\int_{t_{k}}^{t_{k+1}}\max(0,\lambda_{max}(L_{l}(t,q_{v}(t),\zeta_{v}(t)))^{p_{l}}dt\leq0,\label{eq:CTCS_by_max_function}
\end{align}
where the exponent satisfies $p_{l}\geq1$. This equivalence holds because $\lambda_{max}(\cdot)$ is continuous \cite[Chapter 2, Section 5.1]{kato2013perturbation}, and $q(\cdot)$, and $\zeta(\cdot)$ are continuous on each subinterval $[t_{k},t_{k+1}]$; the detailed proof of this equivalence is given in \cite[Lemma 2]{ELANGO2025112464}. Since $q(t)$ and $\zeta(t)$ are linear functions of the decision variables $q_{v,k}$, $\zeta_{v,k}^{m}$, and $\zeta_{v,k+1}^{p}$ within the subinterval, we define the left-hand side of \eqref{eq:CTCS_by_max_function} as:

\begin{equation}
h_{k}^{l}(\chi_{k})\leq0,\quad k\in\mathbb{Z}_{[0,N-1]},\label{eq:def_of_h}
\end{equation}
where $\chi_{k}\coloneqq(q_{v,k},\zeta_{v,k}^{m},\zeta_{v,k+1}^{p})$ collects the relevant decision variables on $[t_{k},t_{k+1}]$, and $h_{k}^{l}$ is an integral functional.
\begin{lem}
The function $h_{k}^{l}(\cdot)$ is convex for each subinterval $[t_{k},t_{k+1}]$.
\end{lem}

\begin{proof}
For each $t\in(t_{k},t_{k+1})$, the matrix-valued function $L_{l}$ is affine in $\chi_{k}$. Since $\lambda_{max}(\cdot)$ is convex, and $\max(0,\cdot)$ is convex and nondecreasing, their composition $\max(0,\lambda_{max}(L_{l}(t,q(t),\zeta(t))$ is convex in $\chi_{k}$ \cite[Chapter 3.2.4]{boyd2004convex}. Finally, integration over $t\in[t_{k},t_{k+1}]$ preserves the convexity, so $h_{k}^{l}(\chi_{k})$ is convex.
\end{proof}
\begin{rem}
Although $h_{k}^{l}(\cdot)$ is convex, it does not take the form of a standard conic constraint, such as second-order cone or a PSD cone. Hence, it cannot be directly handled by conventional conic optimization solvers.
\end{rem}

\subsection{Intermediate constraint-checking points\label{subsec:Intermediate-constraint-checking}}

One way to make $h_{k}^{l}$ tractable is to approximate the integral by a finite sum using standard integration techniques, such as Simpson's method with the midpoint rule \cite{burden2011numerical}. This yields the approximation:
\begin{align*}
& \int_{t_{k}}^{t_{k+1}}\max(0,\lambda_{max}(L_{l}(t,q_{v}(t),\zeta_{v}(t)))^{p_{l}}dt \\ 
& \approx\sum_{s=0}^{N_{s}+1}c_{s}\max(0,\lambda_{max}(L_{l}(t_{k,s},q_{v}(t_{k,s}),\zeta_{v}(t_{k,s})))^{p_{l}},
\end{align*}
where $N_{s}$ is the number of intermediate points, $c_{s}\in\mathbb{R}_{++}$ are integration weights, $t_{k,s}$ for $s\in\mathbb{Z}_{[1,N_{s}]}$ are the evaluation points with $t_{k,0}=t_{k}$ and $t_{k,N_{s}+1}=t_{k+1}$. Note that the following equivalence holds:
\begin{align*}
\sum_{s=0}^{N_{s}+1}c_{s} & \max(0,\lambda_{max}(L_{l}(t_{k,s},q_{v}(t_{k,s}),\zeta_{v}(t_{k,s})))^{p_{l}}  \leq0\\
& \Leftrightarrow L_{l}(t_{k,s},q_{v}(t_{k,s}),\zeta_{v}(t_{k,s}))  \preceq0,\quad\forall s\in\mathbb{Z}_{[0,N_{s}+1]}.
\end{align*}
This equivalence holds because each $c_{s}$ is nonnegative, and the integrand $\max(0,\lambda_{\max}(\cdot))^{p}$ is nonnegative as well. Therefore, the entire sum is nonpositive if and only if each term in the sum is zero, which implies that all matrix inequalities are individually satisfied. 

As a result, approximating the integral by a finite sum leads to the introduction of intermediate constraint-checking points $t_{k,s}\in(t_{k,}t_{k+1})$, at which the LMI constraints \eqref{eq:continuous-time-LMI} are enforced. Under the FOH interpolation \eqref{eq:FOH_interpolation} and the linear funnel dynamics \eqref{eq:vectorized_funl_dynamics}, $q(t_{k,s})$ and $\zeta(t_{k,s})$ at any $t_{k,s}\in[t_{k},t_{k+1}]$ become linear functions of the nodal variables $(q_{k},\zeta_{k,}\zeta_{k+1})$. Therefore, the conditions 
\begin{equation}
L_{l}(t_{k,s},q(t_{k,s}),\zeta(t_{k,s}))\preceq0,\label{eq:additional_LMI_checking_points}
\end{equation}
remain LMIs in the decision variables and can be directly incorporated into the discrete problem formulation \eqref{eq:funnel_synthesis_discrete}. 

The intermediate points $t_{k,s}$ can be uniformly distributed over each subinterval using
\[
t_{k,s}=t_{k}+\frac{s}{N_{s}+1}(t_{k+1}-t_{k}),\quad s\in\mathbb{Z}_{[1,N_{s}]},
\]
where $N_{s}\in\mathbb{R}_{++}$ is the number of intermediate points per subintervals.
\begin{rem}
Adding the intermediate constraint-checking points \eqref{eq:additional_LMI_checking_points} is more efficient compared to merely increasing the number $k$ of node points. While both strategies improve constraint satisfaction over $[t_{0},t_{f}]$, increasing node points leads to more decision variables and thus higher computational cost. In contrast, adding checking points only increases the number of LMI constraints without expanding the optimization variables, offering a more computationally efficient approach. Another benefit of this approach is that it preserves the number of control updates, thereby avoiding unnecessary increases in control frequency, which is often a desirable property in practical applications.
\end{rem}

\subsection{Derivation of subgradient}

Instead of approximating the integral, the second approach is to obtain the (sub) gradient of $h_{k}^{l}(\cdot)$ and apply the successive convexification method illustrated in \cite{ELANGO2025112464}. In this and next subsections, we establish some related theoretical results.

We define the function $g_{k}^{l}$ as
\begin{align}
g_{k}^{l}(\chi_{k})\coloneqq\int_{t_{k}}^{t_{k+1}}\tilde{g}_{k}^{l}(t,\chi_{k})\mathrm{d}t,\label{eq:def_of_subgradient_of_h}
\end{align}
where $\tilde{g}_{k}^{l}\in\partial_{\chi_{k}}\{\max(0,\lambda_{max}(L_{l}(t,q_{v},\zeta_{v})))^{p_{l}}\}$ denotes a subgradient of the integrand with respect to $\chi_{k}$ at time $t$. 
\begin{lem}
If $\tilde{g}_{k}^{l}(t,\chi_{k})$ is Lebesgue integrable over $[t_{k},t_{k+1}]$, then $g_{k}^{l}(\chi_{k})$ is a subgradient of $h_{k}^{l}$ at $\chi_{k}$.
\end{lem}

\begin{proof}
For notational brevity, we abuse the notation $L_{l}$ as a matrix-valued function defined by:
\[
L_{l}(t,\chi_{k})\coloneqq L_{l}(t,q_{v}(t),\zeta_{v}(t)),\quad t\in[t_{k},t_{k+1}].
\]
Since $q_{v}$ and $\zeta_{v}$ are affine in $\chi_{k}$, $L_{l}$ is also affine in $\chi_{k}$ for each fixed $t$. Let $\bar{\chi}_{k}$ be any point in the domain. For each fixed $t\in[t_{k},t_{k+1}]$, by the definition of the subgradient, we have
\begin{align*}
\max(0,\lambda_{max}(L_{l}(t,\bar{\chi}_{k})))^{p_{l}} \geq & \max(0,\lambda_{max}(L_{l}(t,\chi_{k})))^{p_{l}} \\ 
& +\tilde{g}_{k}^{l}(t,\chi_{k})^{\top}(\bar{\chi}^{k}-\chi_{k}).
\end{align*}
This inequality holds pointwise in $t$, and both sides are Lebesgue integrable since $\tilde{g}_{k}^{l}(t,\chi_{k})$ is assumed to be Lebesgue integrable. We can therefore integrate both sides over the interval $[t_{k},t_{k+1}]$, yielding:
\begin{align*}
h_{k}^{l}(\bar{\chi}_{k}) & \geq h_{k}^{l}(\chi_{k})+\int_{t_{k}}^{t_{k+1}}\tilde{g}_{k}^{l}(t,\chi_{k})^{\top}(\bar{\chi}^{k}-\chi_{k})\mathrm{d}t,\\
 & =h_{k}^{l}(\chi_{k})+g_{k}^{l}(\chi_{k})^{\top}(\bar{\chi}^{k}-\chi_{k}).
\end{align*}
Therefore, $g_{k}^{l}(\chi_{k})$ is a subgradient of $h_{k}^{l}$ at $\chi_{k}$.
\end{proof}

About the integrability of $\tilde{g}_{k}^{l}$, we can establish the following result.
\begin{thm}
\label{thm: integrable}
For a subinterval $[t_{k},t_{k+1}]$ and a power $p_{l}\geq1$, there exists a Lebesgue measurable selection $\tilde{g}_{k}^{l}(t,\chi_{k})\in\partial_{\chi_{k}}\{\max(0,\lambda_{max}(L_{l}(t,\chi_{k})))^{p_{l}}\}$ for almost every $t\in[t_{k},t_{k+1}]$ and this selection is bounded on the subinterval. Hence, $\tilde{g}_{k}^{l}(t,\chi_{k})$ is Lebesgue integrable.
\end{thm}

\begin{proof}
We apply results from \cite{rockafellar1998variational}. The function $\lambda_{max}(L_{l}(t,\chi_{k}))$ is measurable in $t$ for any fixed $\chi_{k}$, and continuous in $\chi_{k}$ for each $t\in(t_{k},t_{k+1})$. Therefore, it qualifies as a normal integrand \cite[Example 14.29]{rockafellar1998variational}. According to Proposition 14.44, the composition $\max(0,\lambda_{max}(L_{l}(t,\chi_{k})))^{p_{l}}$ is also a normal integrand. Then, by \cite[Theorem 14.56]{rockafellar1998variational}, the subgradient mappings $\partial_{\chi_{k}}\{\max(0,\lambda_{max}(L_{l}(t,\chi_{k})))^{p_{l}}\}$ are closed-valued and measurable. Finally, \cite[Corollary 14.6]{rockafellar1998variational} ensures that this measurable set-valued map admits a measurable selection. Because $L_{l}(t,\chi_{k})$ is affine in $\chi_{k}$ (hence Lipschitz), $\lambda_{max}$ is 1-Lipschitz on symmetric matrices, and the scalar maps $s\rightarrow\max(0,s)$ and $s\rightarrow s^{p_{l}}$ are Lipschitz on every bounded interval, their composition is Lipschitz on the compact set $[t_{k},t_{k+1}]\times\operatorname{dom}(\chi)$ where $\operatorname{dom}(\chi)\subset\mathbb{R}^{n_{q_{v}}+n_{y_{v}}+n_{z_{v}}}$ denotes the bounded domain of $\chi_{k}$. Consequently, its subgradient is bounded. Hence, the measurable selection $\tilde{g}_{k}^{l}(t,\chi_{k})$ is integrable.
\end{proof}

While Theorem \ref{thm: integrable} guarantees that a integrable selection $\tilde{g}_{k}^{l}(t,\chi_{k})$ exists almost everywhere on $[t_{k},t_{k+1}]$, the lemma is purely existential; it does not specify which subgradient is selected or how to compute it. For implementation, we choose a specific subgradient of the integrand $\max(0,\lambda_{max}(L_{l}(t,q_{v},\zeta_{v})))^{p_{l}}$ at $\bar{\chi}_{k}$, obtained from the chain rule. Let $v_{i}$ be the eigenvectors associated with $\lambda_{max}(L_{l}(t,\chi_{k}))$ and $M(t,\chi_{k})$ the algebraic multiplicity of $\lambda_{max}(L_{l})$. Define the matrix
\[
W(t,\chi_{k})=\frac{1}{M(t,\chi_{k})}\sum_{i:\lambda_{i}=\lambda_{max}(L_{l})}v_{i}v_{i}^{\top},
\]
which averages the outer products of the eigenvectors associated with the largest eigenvalue. The subgradient we use at $\bar{\chi}_{k}$ is
\begin{equation}
\tilde{g}_{k}^{l}(t,\chi_{k})=\left[\begin{smallmatrix}
\frac{\partial\max(0,\lambda_{max}(L_{l}(t)))^{p_{l}}}{\partial\lambda_{max}(L_{l}(t))}\cdot\frac{\partial\vecc{L_{l}(t)}}{\partial\chi_{k}}\cdot\vecc{W(t,\chi_{k})}\end{smallmatrix}\right]_{\chi_{k}=\bar{\chi}_{k}},\label{eq:ghat_computation}
\end{equation}
where $L_{l}(t,\chi_{k})$ is compactly written as $L_{l}(t)$. When $p_{l}>1$ and the largest eigenvalue of $L_{l}$ is simple ($M(t,\chi_{k})=1$), the integrand $\max(0,\lambda_{max}(L_{l}(t,q_{v},\zeta_{v})))^{p_{l}}$ is differentiable. In that case the expression above coincides with the exact gradient of the integrand.

\subsection{Successive convexification with subgradients \label{subsec:Successive-convexification-with}}

We now apply the successive convexification (SCvx) framework proposed in \cite{ELANGO2025112464} to solve the convex constraint $h_{k}^{l}(\chi_{k})\leq0$, using the subgradient $g_{k}^{l}(\chi_{k})$ derived in the previous subsection. The SCvx method can be interpreted as a combination of exact penalization and prox-linear method \cite{drusvyatskiy2019efficiency} that is sequential convex programming for convex-composite minimization. The goal is to solve the discrete funnel synthesis problem with the $\epsilon$-relaxed CTCS constraint:
\begin{subequations}
\label{eq:scvx_mainproblem}
\begin{align}
    \underset{\resizebox{0.2\hsize}{!}{$q_{v,k},\zeta_{v,k}^{m},\zeta_{v,k}^{p}$}}{\operatorname{minimize}}~~& J(Q_{0}) \\[-0.1cm]
    \operatorname{subject~to}~~~& \eqref{eq:discrete_const_on_funl_dynamics},\eqref{eq:nodewise_LMI},\eqref{eq:discrete_funl_multiplier_matrices},\eqref{eq:discrete_terminal}\label{eq:scvx_same_constraints-1}\\
 & h_{k}^{l}(\chi_{k})\leq\epsilon,\quad\forall k\in\mathbb{Z}_{[0,N-1]},
\end{align}
\end{subequations}
where $\epsilon\in\mathbb{R}_{++}$ is a small constant introduced to ensure linear independence constraint qualification (LICQ), as discussed in \cite[Lemma 10]{ELANGO2025112464}. In practice, $\epsilon$ is typically chosen to be sufficiently small, for example, less than $10^{-6}$, so that it does not significantly affect the quality of continuous-time constraint satisfaction. Note that the constraints in \eqref{eq:scvx_same_constraints-1} include those already present in the discrete funnel synthesis problem \eqref{eq:funnel_synthesis_discrete}. 

At each iteration, the original constraint for the CTCS is replaced by its first-order approximation using a subgradient. Specifically, at the $i$-th iteration, the following convex SDP subproblem is solved:
\begin{problem}
(SDP subproblem)
\begin{subequations}
\label{eq:scvx_subproblem}
\begin{align}
    \underset{\resizebox{0.2\hsize}{!}{$q_{v,k},\zeta_{v,k}^{m},\zeta_{v,k}^{p}$}}{\operatorname{minimize}}~~&
    \begin{smallmatrix}
        J(Q_{0})+w_{h}\sum_{k=0}^{N-1}\max(0,v_{k}^{l}) \\ 
        +w_{tr}\left(\sum_{k=0}^{N-1}\norm{\chi_{k}-\bar{\chi}_{k}}_{2}^{2}+\norm{q_{v,N}-\bar{q}_{v,N}}_{2}^{2}\right)
    \end{smallmatrix} \label{eq:scvx_cost_function} \\[-0.1cm]
 \operatorname{subject~to}~~~& \eqref{eq:discrete_const_on_funl_dynamics},\eqref{eq:nodewise_LMI},\eqref{eq:discrete_funl_multiplier_matrices},\eqref{eq:discrete_terminal}\label{eq:scvx_same_constraints}\\
 & h_{k}^{l}(\bar{\chi}_{k})+g_{k}^{l}(\bar{\chi}_{k})^{\top}(\chi_{k}-\bar{\chi}_{k})-\epsilon=v_{k}^{l},\label{eq:scvx_linearized_constraint} \\
 & \forall k\in\mathbb{Z}_{[0,N-1]}, \nonumber
\end{align}
\end{subequations}
where the reference solutions $\bar{\chi}_{k}$ for all $k\in\mathbb{Z}_{[0,N-1]}$ and $\bar{q}_{v,N}$ come from an initial guess or the previous iteration's solution. The penalization terms in the objective reflect two key components of the SCvx method. The first term weighted by $w_{h}\in\mathbb{R}_{++}$ enforces the linearized constraint through exact penalization. The second term weighted by $w_{tr}$ promotes proximity to the linearization point $\bar{\chi}_{k}$ by penalizing the trust-region. This term also aligns the subproblem structure with the prox-linear method framework \cite{drusvyatskiy2019efficiency}.
\end{problem}



\begin{corr}
\label{cor:scvx_convergence}
Suppose that, for almost every $t \in [t_0, t_f]$, the maximum eigenvalue of $L_l(t,\chi_k)$ is simple, so that the subgradient $g_k^l(\bar{\chi}_k)$ coincides with the true gradient of $L_l(t,\chi_k)$. Then, under the assumptions of~\cite[Theorem~30]{ELANGO2025112464}, the sequence generated by the SCvx method converges to a solution of~\eqref{eq:scvx_mainproblem}.
\end{corr}

\begin{proof}
The result follows directly from the convergence statement of~\cite[Theorem~30]{ELANGO2025112464}, applied to the case where $g_k^l(\bar{\chi}_k)$ coincides with the true gradient of $L_l(t,\chi_k)$.
\end{proof}

\subsection{Summary of algorithm \label{subsec:Summary-of-algorithm}}

We summarize the proposed method for solving the continuous-time funnel synthesis problem with CTCS. The method proceeds in two stages. First, we solve the discrete funnel synthesis problem with additional constraint enforcement at a finite number of $N_{s}$ intermediate constraint-checking points, as described in Section \ref{subsec:Intermediate-constraint-checking}. The resulting solution is then used as an initial guess for the Scvx method with subgradients, presented in Section \ref{subsec:Successive-convexification-with}. The overall algorithm is as follows:
\begin{algorithm}[H]
\caption{Funnel Synthesis with SCvx}
\begin{algorithmic}[1]
\State Choose $N$, $N_{s}$, and initialize problem data
\State Solve the discrete funnel synthesis problem \eqref{eq:funnel_synthesis_discrete} with additional constraint enforcement \eqref{eq:additional_LMI_checking_points} at $N_{s}$ intermediate points per subinterval
\State Set an initial guess $\bar{\chi}_{k}$ for all $k\in\mathbb{Z}_{[0,N-1]}$, and $\bar{q}_{v,N}$
\For{$i = 1, \ldots, N_{iter}$}
    \State Compute the subgradient $g_{k}^{l}(\bar{\chi}_{k})$ for each constraint $h_{k}^{l}$ using \eqref{eq:ghat_computation}
    \State Solve the SDP subproblem \eqref{eq:scvx_subproblem}.
    \State Update the solutions $\bar{\chi}_{k}\leftarrow\chi_{k}$ for all $k\in\mathbb{Z}_{[0,N-1]}$ and $\bar{q}_{v,N}\leftarrow q_{v,N}$.
\EndFor
\State Return solution variables
\end{algorithmic}
\end{algorithm}

\section{Numerical examples \label{sec:6}}

In this section, we validate the proposed method through two numerical examples involving obstacle avoidance: control of a unicycle and a 6-DoF quadrotor. In both cases, a Mosek solver is used to solve the SDP \eqref{eq:funnel_synthesis_discrete} and \eqref{eq:scvx_subproblem}. All simulations are written in Julia and executed on a MacBook with an Apple M1 Pro processor.

\subsection{Unicycle}

We consider a unicycle model illustrated in Example \ref{exa:unicycle}. The simulation environment for all unicycle experiments is configured as follows. The time horizon is set to $t_{0}=0$ and $t_{f}=5$ seconds, uniformly divided into $N=9$ subintervals. Two circular obstacles, each with a radius of 0.5\,m, define state constraints, and input constraints are given by $0\leq u_{1}\leq2$ (m/s) and $\abs{u_{2}}\leq2$ (rad/s), as specified in \eqref{eq:constraint_sets}. The cost function is chosen as $\trace{Q_{w0}^{-1}}$, with the weight matrix $W_{0}=I_{3}$ as defined in \eqref{eq:funnel_cost}. The nominal trajectory starts at $\bar{x}(t_{0})=(0,0,0)$ and ends at $\bar{x}(t_{f})=(5,5,0)$. The final funnel matrix $Q_{f}$, used as the terminal constraint in \eqref{eq:discrete_terminal}, is set to diag(0.08, 0.08, 0.06).

\begin{figure} \begin{center}  \includegraphics[width=0.7\linewidth]{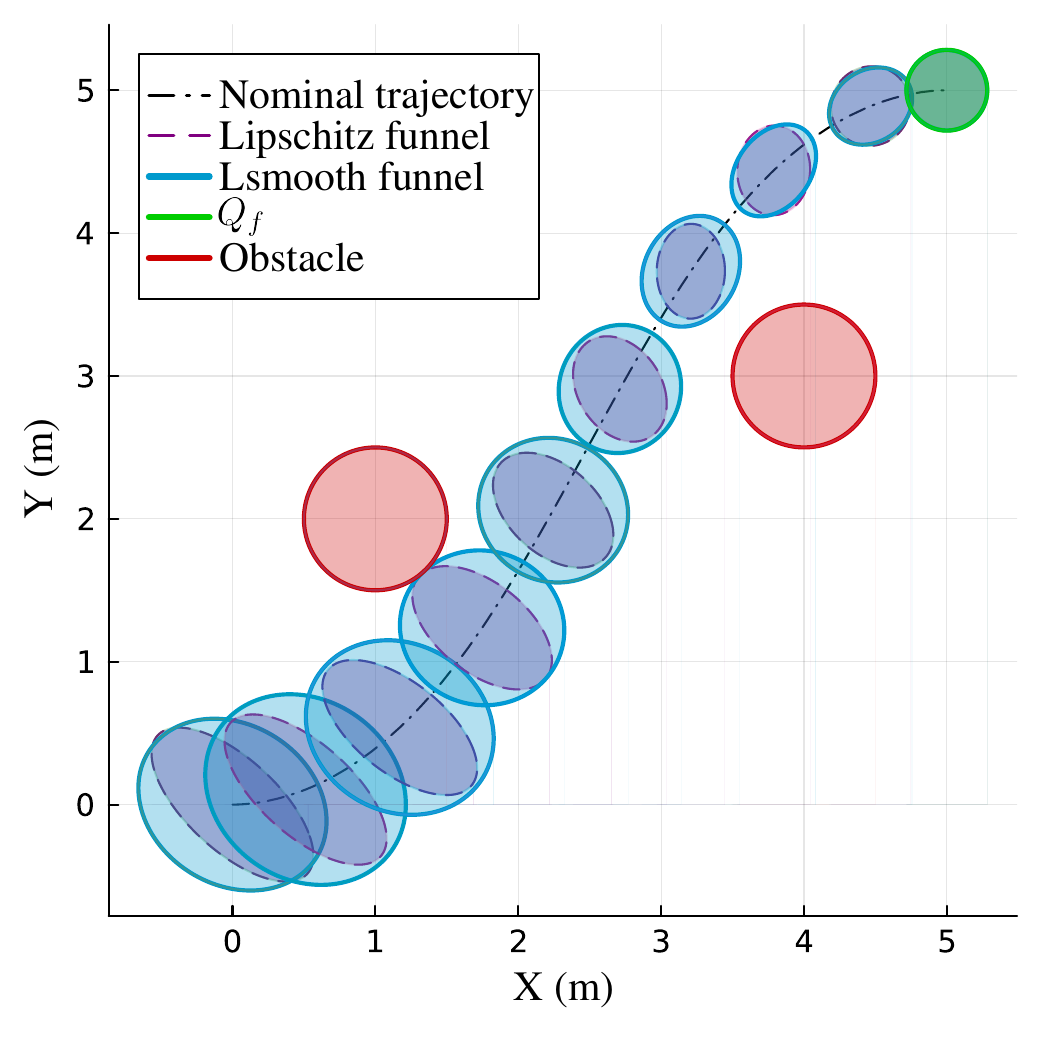} \caption{The synthesized state funnel projected on $x$ ($x_1$) and $y$ ($x_2$) position coordinates.} \label{fig:01_synthesized_funnel_Lip_vs_Lsm} \end{center} \end{figure}

\subsubsection{Comparison of funnels computed under Lipschitz and L-smooth conditions}

We first present results comparing two modeling approaches for the nonlinearity: the Lipschitz condition \eqref{eq:Lipshitz} and L-smooth condition \eqref{eq:Lsmooth}. Specifically, we consider the use of global constants, where $\gamma_{i}$ and $\beta_{i}$ are valid over the entire domain $\Omega=\mathbb{R}^{n_{q_{i}}}$ in \eqref{eq:Lipshitz}and \eqref{eq:Lsmooth} for each $i=1,2$ and $t\in[t_{0},t_{f}]$. As discussed in Example \ref{exa:unicycle_2}, the input constraint $\abs{u_{1}}\leq2$ implies that the global constants are $\gamma_{i}=4$ and $\beta_{i}=2$ for all $i=1,2$. For the L-smooth case, $\lambda_{i}^{\beta}$ is set to 0.3 for all $i=1,2$, as discussed in \eqref{eq:inclusion_Lsmooth}. In this comparison, the bounded disturbance is set to zero $w(t)=0$ for all $t\in[t_{0,}t_{f}]$, and $\fdyn$-type funnel dynamics \eqref{eq:Basic_type_funldyn} is used. The decay rate $\alpha$ is set to zero, as invariance can be guaranteed in the absence of bounded disturbances, as discussed in Remark \ref{rem:alpha}. With zero disturbance, the level constant $c_{Q}$ in \eqref{eq:state_funnel} can be any positive value, and it is chosen to be 1 in this case.

The synthesized funnel using $\beta_{i}(t)=2$ for all $t\in[t_{0},t_{f}]$ is shown in Fig.~\ref{fig:01_synthesized_funnel_Lip_vs_Lsm}. It is found that the SDP problem \eqref{eq:funnel_synthesis_discrete} becomes infeasible when using the global Lipschitz constant, indicating that this formulation is overly conservative compared to the L-smooth case. To continue the comparison, we reduce the Lipschitz constant until the SDP becomes feasible, which occurs at $\gamma_{i}=0.7$, which a value significantly smaller than the true global value. The resulting funnel is also shown in Fig.~\ref{fig:01_synthesized_funnel_Lip_vs_Lsm}. Notably, even though the L-smooth formulation uses the exact global value, it yields a larger funnel with a cost of 8.17. In contrast, the funnel computed with the Lipschitz condition has a significantly higher cost of 14.76, further demonstrating that the L-smooth condition is much less conservative.

\begin{figure} \begin{center}  \includegraphics[width=0.7\linewidth]{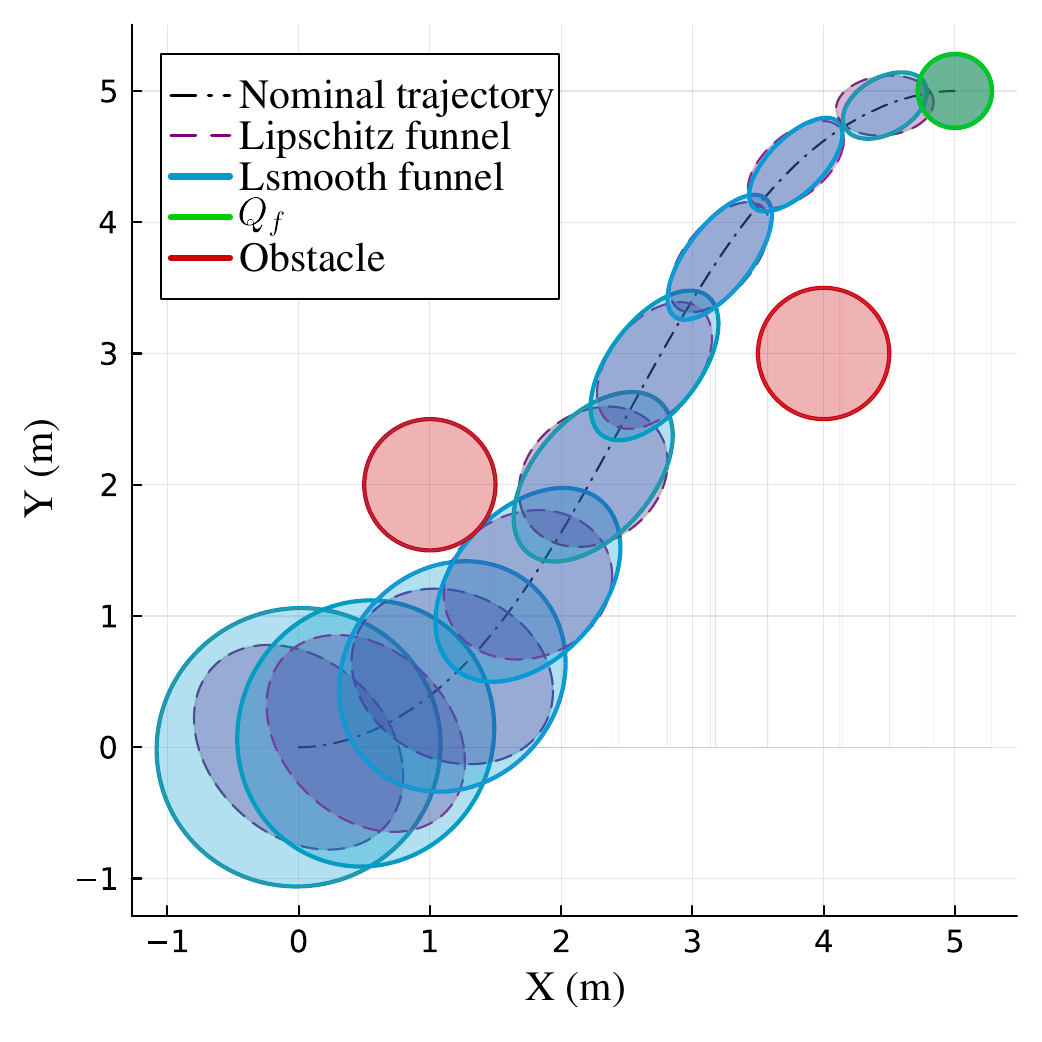} \caption{The synthesized state funnel with local Lipschiptz and L-smooth constants.} \label{fig:02_funnel} \end{center} \end{figure}

 To further compare the Lipschitz and L-smooth conditions, we consider the use of local constants where the domain $\Omega$ is strictly subset of Euclidean space. To this end, we first compute a funnel without considering the nonlinearity, which is equivalent to setting $\gamma_{i}=0$ for all $i=1,2$. Since the nonlinearity is ignored, the funnel is larger, but the invariance no longer holds exactly for the original nonlinear system. This funnel is only used to collect state and input samples, from which the Lipschitz and L-smooth constants are estimated for each subinterval; the corresponding values $\gamma_{k}$ and $\beta_{k}$ are used in the funnel computation \eqref{eq:funnel_synthesis_discrete}. As the estimation of these constants is not the focus of this paper, we refer the reader to \cite{reynolds2020temporally,reynolds2021funnel} for further details. For the L-smooth case, $\lambda_{i}^{\beta}$ is set to 1.0 for all $i=1,2$. The resulting funnel are shown in Fig.~\ref{fig:02_funnel}. It is shown that the L-smooth case achieves a larger funnel with a lower cost 2.96 compared to the Lipschitz case 5.11, again illustrating the reduced conservatism of the L-smooth condition.

\subsubsection{Comparison of funnel dynamics}

\begin{figure}
\centering
\begin{subfigure}[t]{0.49\linewidth}
\centering
\includegraphics[width=\linewidth]{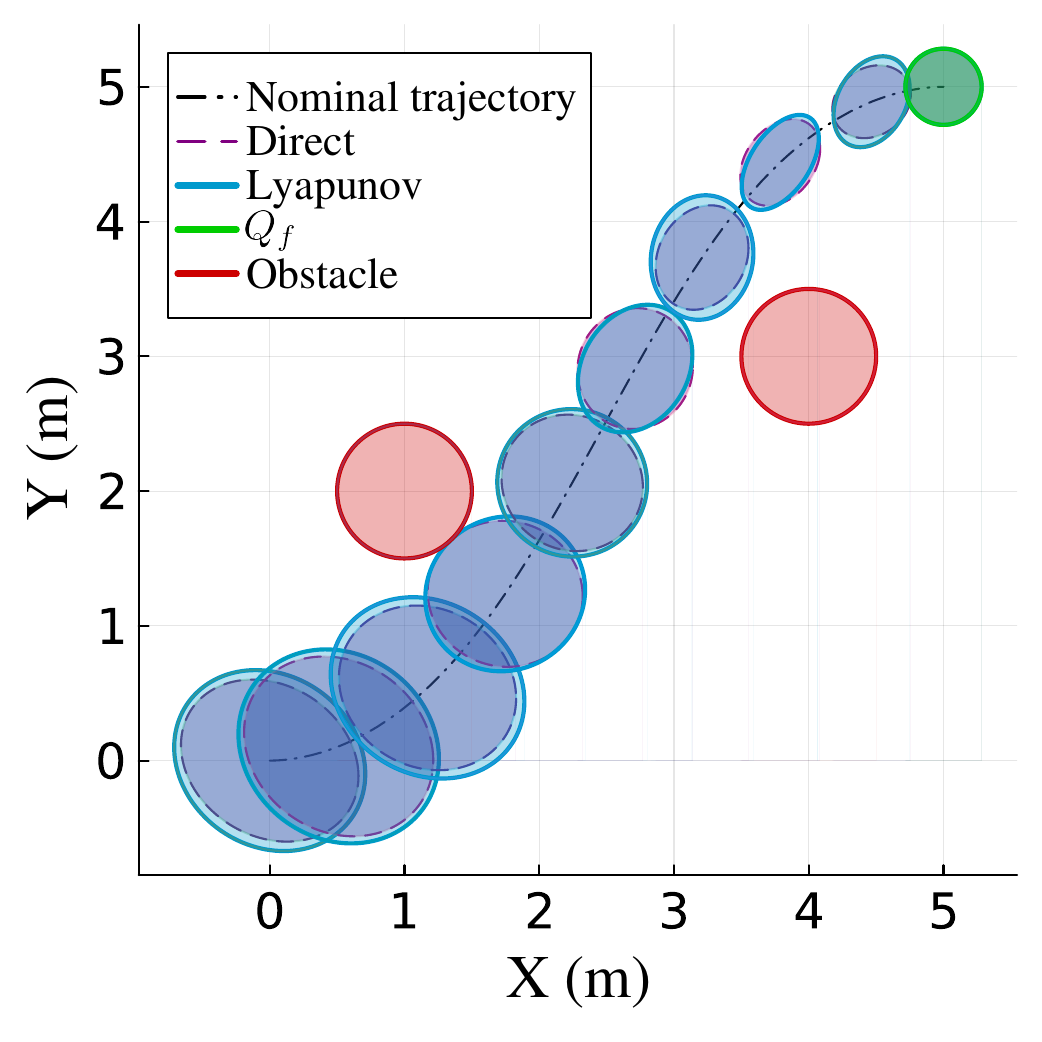}
\end{subfigure}
\begin{subfigure}[t]{0.49\linewidth}
\centering
\includegraphics[width=\linewidth]{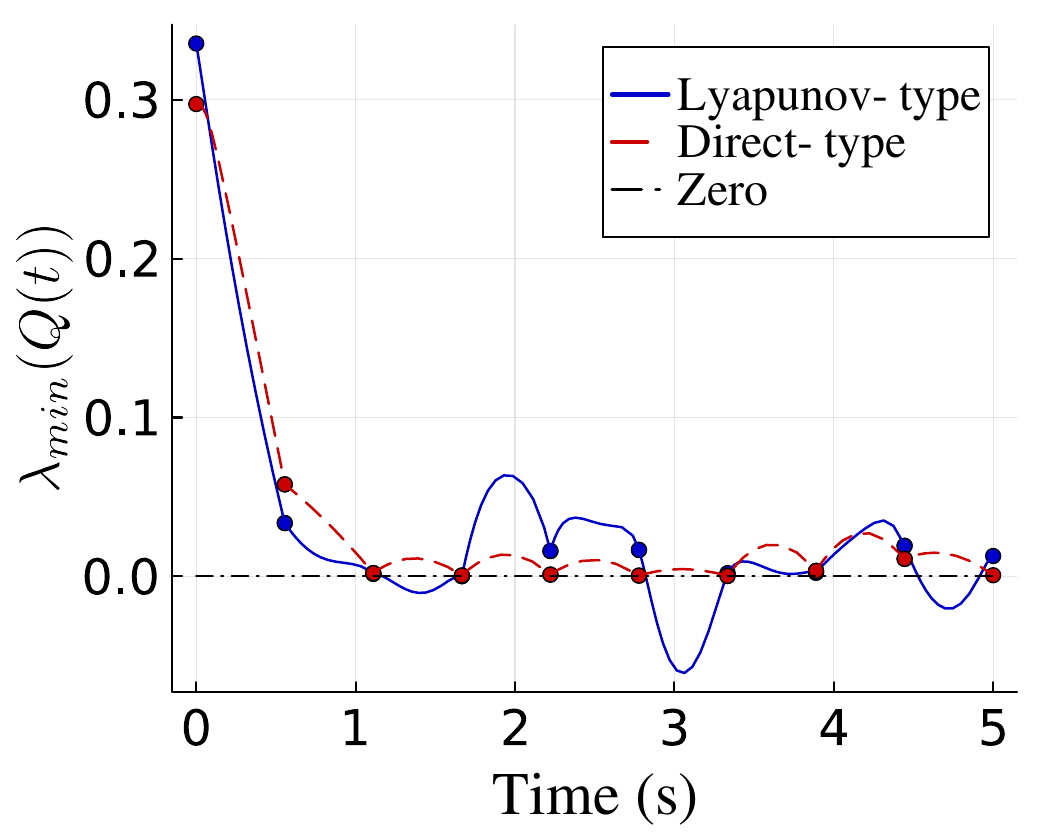}
\end{subfigure}
\caption{(Top): The synthesized state funnel with Lyapunov-type and $\fdyn$-type dynamics. (Bottom): The time history of minimum eigenvalues of $Q(t)$.}
\label{fig:03_comparison_lyp_direct}
\end{figure}

In the next simulations, we aim to compare the choice of funnel dynamics: Lyapunov-type \eqref{eq:Lyapunov_type_funldyn} and $\fdyn$-type \eqref{eq:Basic_type_funldyn}. For this comparison, we choose the use of global L-smooth constants $\beta_{i}=2$ for all $i=1,2$. The bounded disturbance is ignored, and $\lambda_{i}^{\beta}$ is set to 0.3 for all $i=1,2$. The computed funnels are illustrated in the upper panel of Fig. \ref{fig:03_comparison_lyp_direct}. It is shown that the funnel computed with the Lyapunov-type dynamics achieves a larger entry size with a lower cost 7.17 compared to the funnel with the Direct-type dynamics with a cost 8.17. We conjecture that this is because the Lyapunov-type formulation effectively encodes the DLMI \eqref{eq:DLMI_main} through a Lyapunov differential equation \eqref{eq:Lyapunov_type_funldyn}, allowing the funnel\textquoteright s evolution over time to be captured more naturally and smoothly. In contrast, the direct-type formulation constrains the funnel evolution to follow a second-order hold structure, which may limit its expressiveness. However, as discussed in Section \ref{subsec:Positive-Definiteness-of-Q}, the solution of the Lyapunov-type dynamics does not preserve the positive definiteness of $Q$. The lower panel of Fig.~\ref{fig:03_comparison_lyp_direct} shows the time evolution of the minimum eigenvalues of $Q(t)$ for both funnel dynamics formulations. The results reveal that, for the Lyapunov-type dynamics, the minimum eigenvalue occasionally falls below zero, indicating a loss of positive definiteness and hence invalidating the funnel. In contrast, as guaranteed by Lemma \ref{lem:PD_of_Q_with_direct}, the minimum eigenvalue remains strictly positive under the $\fdyn$-type dynamics, ensuring the funnel\textquoteright s validity throughout the time horizon.

\subsubsection{CTCS by introducing intermediate constraint-checking points}

\begin{figure}   
\centering   
\begin{subfigure}[t]{0.49\linewidth}
\includegraphics[width=\linewidth]{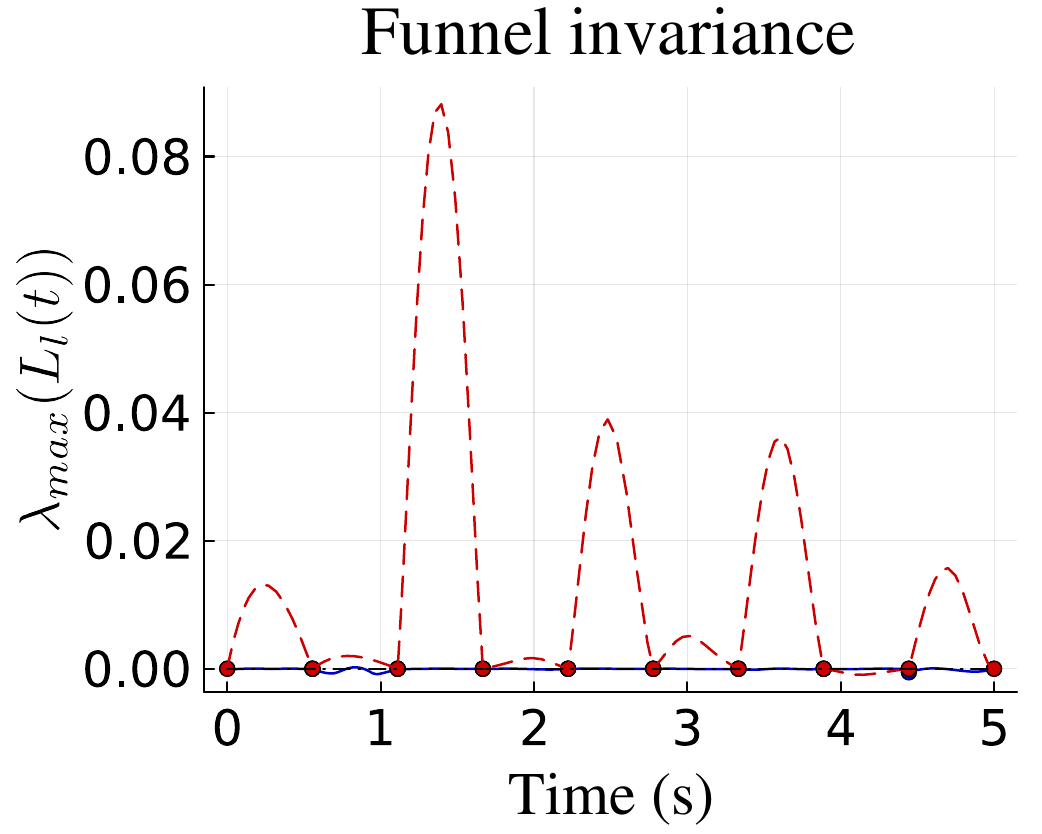}   \end{subfigure}
\begin{subfigure}[t]{0.49\linewidth}     
\includegraphics[width=\linewidth]{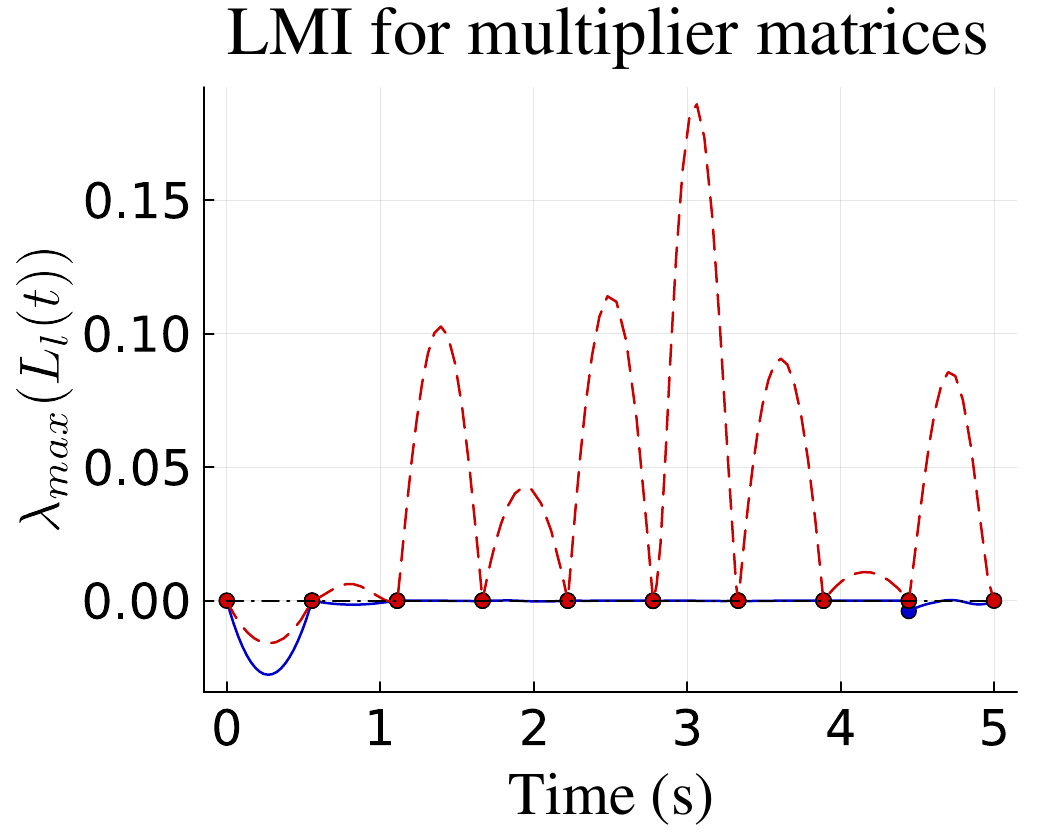}
\end{subfigure}
\begin{subfigure}[t]{0.49\linewidth}     
\includegraphics[width=\linewidth]{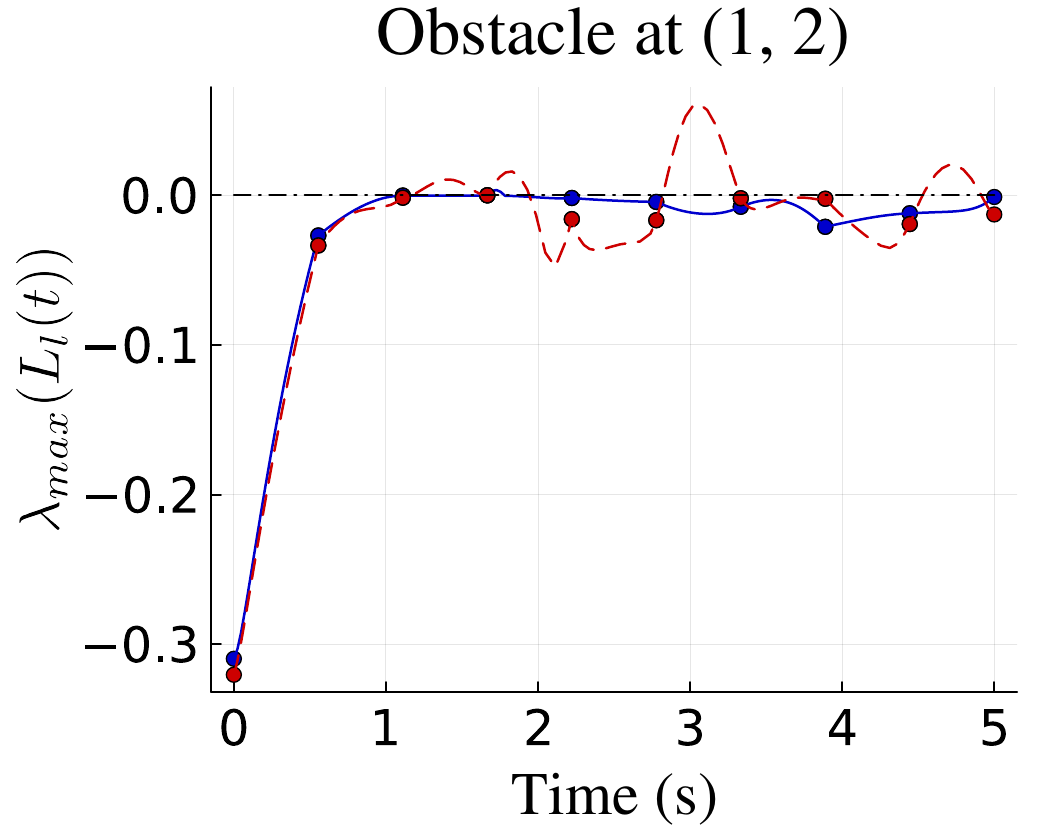}
\end{subfigure} 
\begin{subfigure}[t]{0.49\linewidth}     
\includegraphics[width=\linewidth]{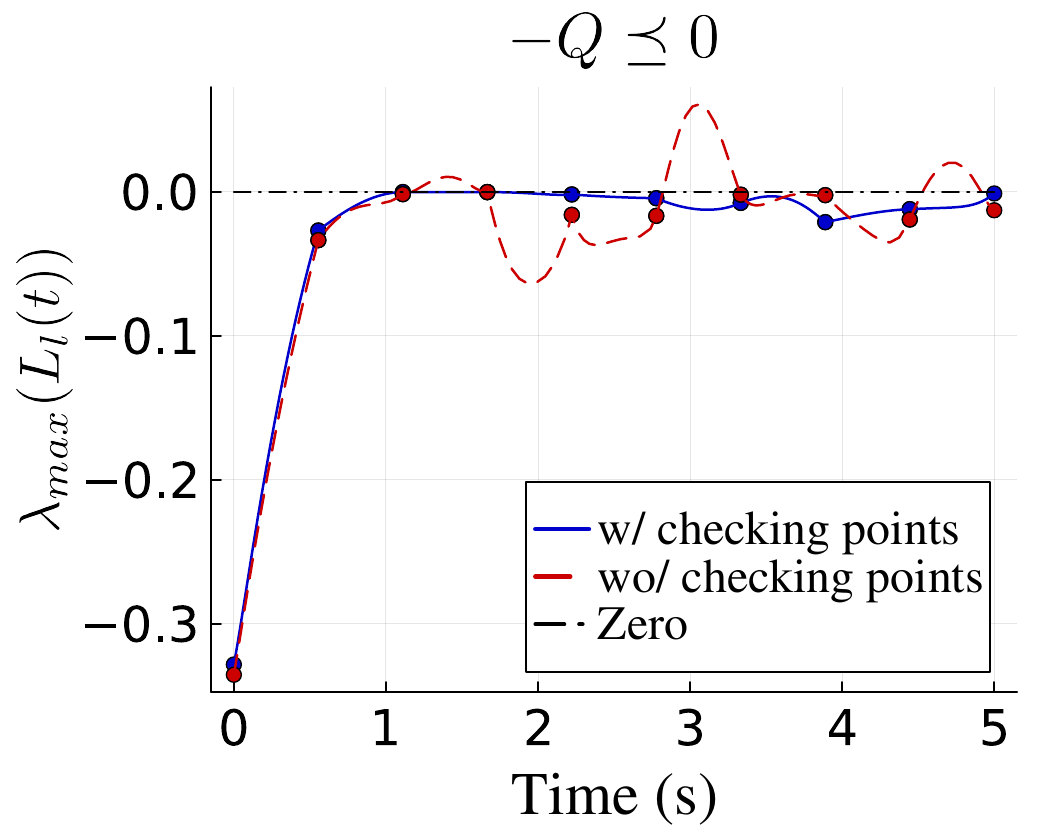}   
\end{subfigure}
\caption{Time evolution of the maximum eigenvalue of each pointwise-in-time LMI constraint $L_l(t)$ associated with the corresponding constraint labeled in each subplot title. “Obstacle 1” denotes the obstacle at (1,2) in the xy-coordinate plane.}   
\label{fig:05_constraint_summary}
\end{figure}

Next, we demonstrate the effectiveness of introducing intermediate constraint-checking points as described in \eqref{eq:additional_LMI_checking_points}. In this simulation, we use the Lyapunov-type funnel dynamics, and all other settings are kept identical to the previous experiment. We add $N_{s}=4$ checking points per subinterval. The time evolution of the maximum eigenvalues of representative pointwise-in-time LMI constraints from \eqref{eq:continuous-time-LMI} is shown in Fig. \ref{fig:05_constraint_summary}. As illustrated, for each constraint $L_{l}$, the maximum eigenvalue remains below zero across the entire time horizon, indicating that CTCS is achieved.

\begin{figure}
\begin{center}
\includegraphics[width=0.7\linewidth,trim={0cm 0.0cm 0cm 0cm},clip]{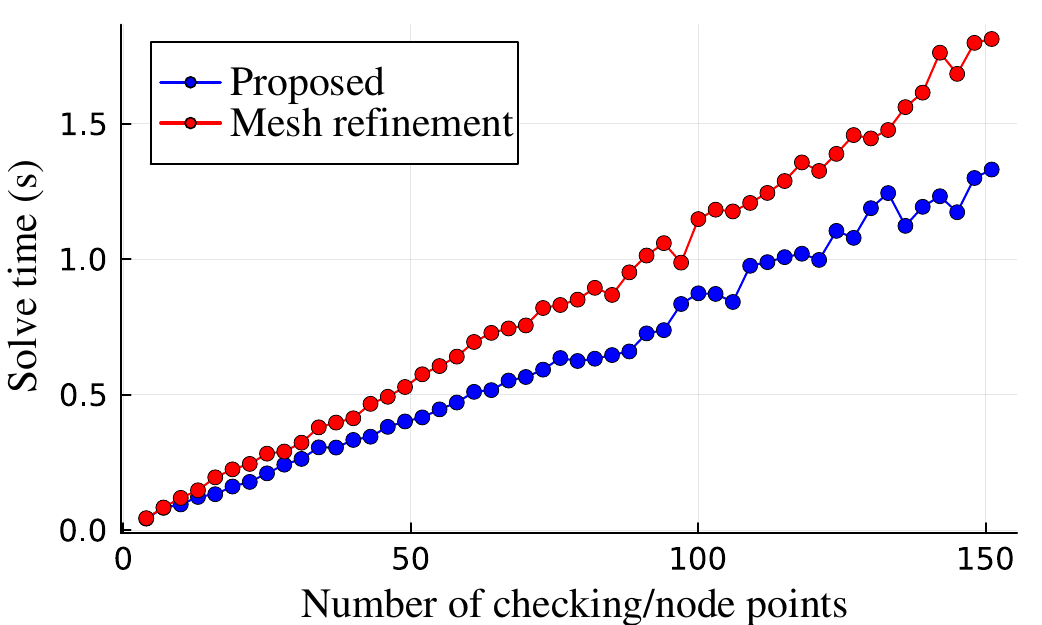} \caption{Mosek solve time versus number of constraints for two approaches: increasing node points and adding intermediate constraint-checking points (proposed).} \label{fig:06_comparison_ICP_vs_mesh} \end{center} \end{figure}

Figure \ref{fig:06_comparison_ICP_vs_mesh} compares the computational cost, measured by Mosek\textquoteright s solve time, of two strategies for improving CTCS: increasing the number of node points and introducing intermediate constraint-checking points. While both approaches enhance constraint coverage, increasing node points enlarges the decision variable space, which leads to higher computational complexity. In contrast, adding intermediate checking points preserves the number of decision variables and only increases the number of LMI constraints. Although the solve-time difference may appear modest, this is partly due to Mosek\textquoteright s internal handling, where it introduces auxiliary variables for the added constraints. Still, the intermediate-checking-point method consistently achieves lower solve times and avoids increasing the control update rate, which is often a desirable property in practice.

\subsubsection{CTCS by SCvx with subgradients}

\begin{figure} \begin{center}  \includegraphics[width=6.0cm,trim={0cm 0.0cm 0cm 0cm},clip]{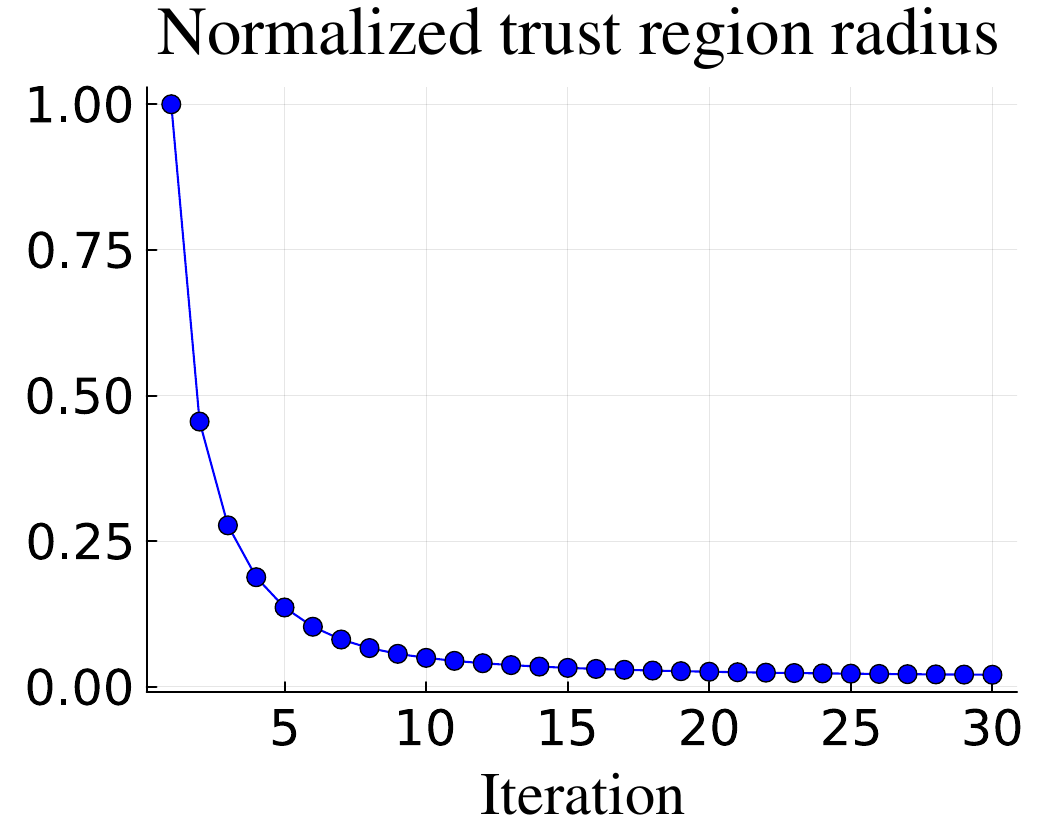} \caption{Normalized trust region penalization term from the SCvx cost function \eqref{eq:scvx_cost_function} over 30 iterations. The values are normalized to 1 at the first iteration to clearly illustrate the convergence behavior.} \label{fig:07_trust_region} \end{center} \end{figure}

In the following simulation, we validate the proposed SCvx-based approach with subgradients for CTCS, as described in Section \ref{subsec:Successive-convexification-with}. We employ the $\fdyn$-type funnel dynamics, while keeping all other settings identical to the previous experiment. As noted in Section \ref{subsec:Summary-of-algorithm}, the solution to the problem \eqref{eq:funnel_synthesis_discrete} with one intermediate checking point ($N_{s}=1$) for each subinterval is used as the initial guess for the SCvx iteration. The weights $w_{vc}$ and $w_{tr}$ in \eqref{eq:scvx_cost_function} are set to 20 and $10^{5}$, respectively. A total of 30 iterations are performed, with an average Mosek solve time of 0.2 seconds per iteration. The small epsilon in \eqref{eq:scvx_linearized_constraint} is set to $\text{\ensuremath{\epsilon=}}10^{-6}$. The time evolution of the maximum eigenvalues of the pointwise-in-time LMI constraints associated with the invariance condition \eqref{eq:Basic_type_H} and the validity of the multiplier matrices \eqref{eq:inclusion_Lsmooth_LMI} is shown in \ref{fig:07_scvx_ctcs_eigenvalues}. The plots indicate that as the iterations proceed, the constraint violation between node points diminishes. The trust region penalization term in \eqref{eq:scvx_cost_function} is plotted in Fig.~\ref{fig:07_trust_region}. The values are normalized such that the initial iteration starts at 1, for the purpose of clearly illustrating the convergence behavior.

\begin{figure}   \centering   \begin{subfigure}[t]{0.49\linewidth}     \includegraphics[width=\linewidth]{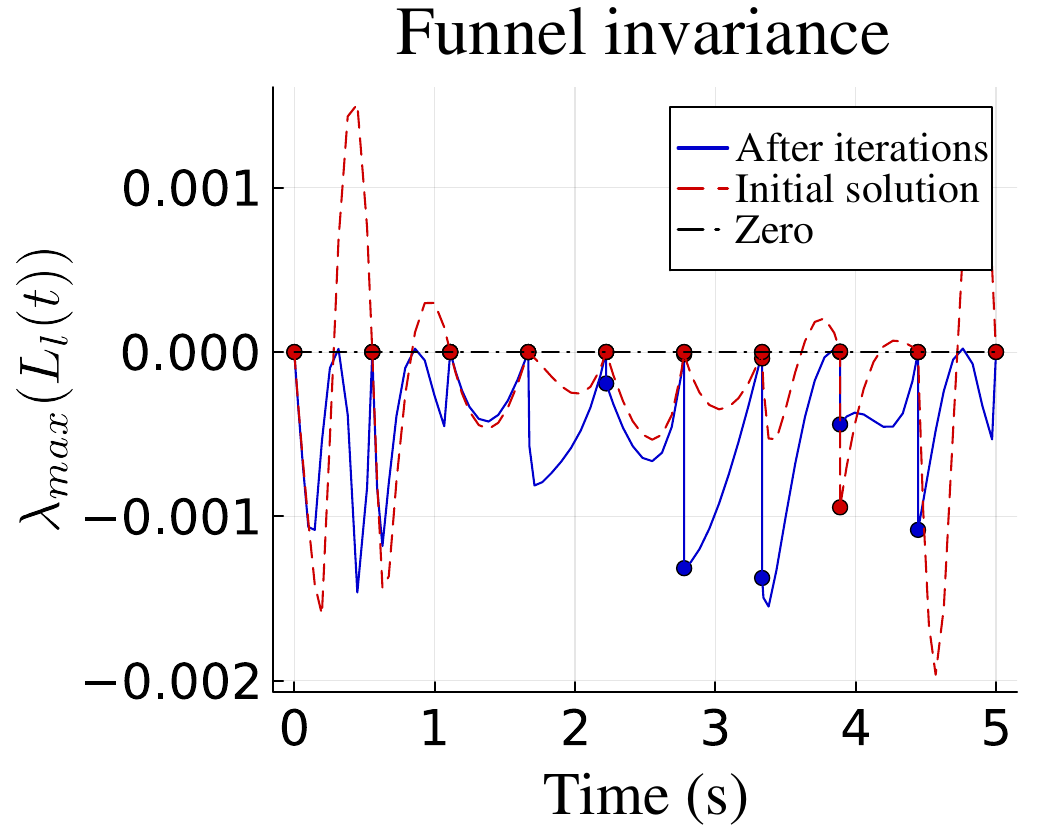}   \end{subfigure}   \begin{subfigure}[t]{0.49\linewidth}     \includegraphics[width=\linewidth]{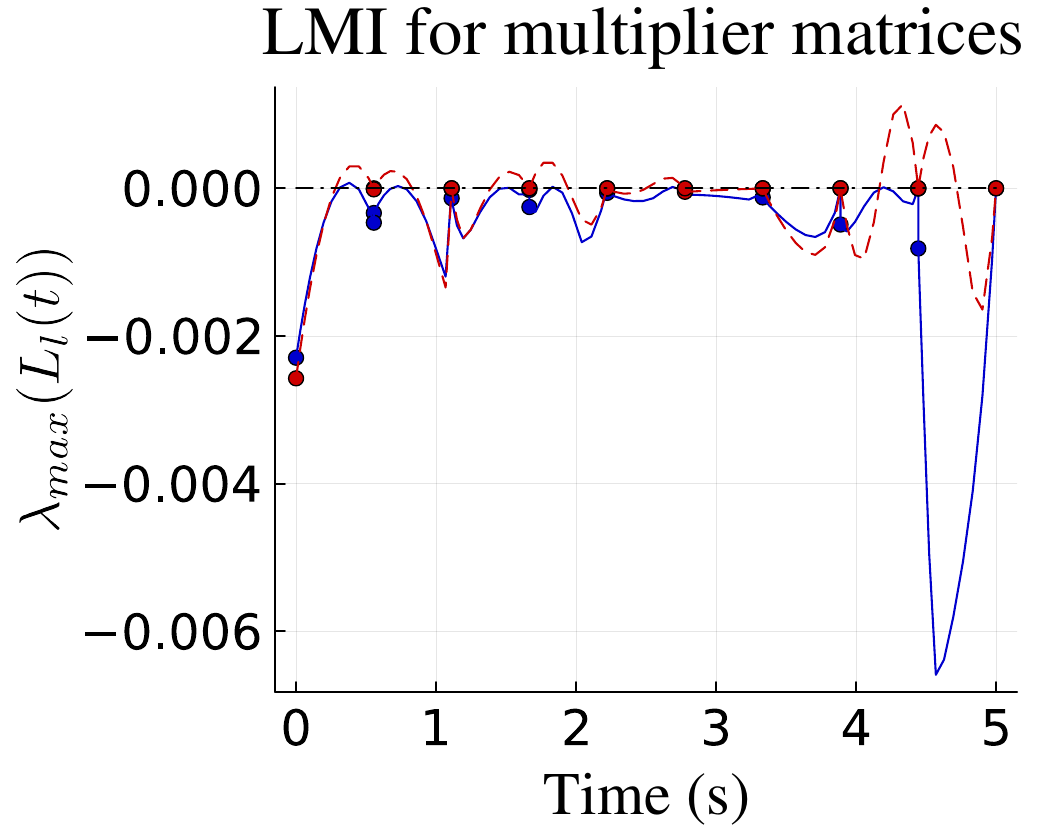}   \end{subfigure}   \caption{Time evolution of maximum eigenvalues for pointwise-in-time LMI constraints in the SCvx-based CTCS approach. Left: Funnel invariance constraint \eqref{eq:Basic_type_H}, Right: valid multiplier constraint \eqref{eq:inclusion_Lsmooth_LMI}.}   \label{fig:07_scvx_ctcs_eigenvalues} \end{figure}

\subsubsection{Funnel under bounded disturbance}

Next, we consider the case where the bounded disturbance exerts to the unicycle model:

\[
\left[\begin{array}{c}
\dot{x}_{1}\\
\dot{x}_{2}\\
\dot{x}_{3}
\end{array}\right]=\left[\begin{array}{c}
u_{1}\cos x_{3}\\
u_{1}\sin x_{3}\\
u_{2}
\end{array}\right]+\left[\begin{array}{c}
0.02w_{1}\\
0.02w_{2}\\
0
\end{array}\right],
\]
where the state and input variables are the same as those in Example \eqref{exa:unicycle}, and $w=[w_{1},w_{2}]^{\top}$ represents a bounded disturbance with $w_{max}=1$ such that $\norm w_{2}\leq1$. Following the discussion in Lemma 1\ref{lem:state_funnel_invariance}, the level constant $c_{Q}$ in the state funnel definition \eqref{eq:state_funnel} is chosen to $c_{Q}=w_{\max}^{2}$. The decay rate is set to$\alpha=0.1$. The $\fdyn$-type funnel dynamics \eqref{eq:Basic_type_funldyn} is employed and two intermediate checking points ($N_{s}=2$) are used for each subinterval. All other settings including the cost function, constraint formulations, and simulation parameters are kept identical to those used in the previous undisturbed simulations.

\begin{figure}
\begin{center}

\begin{subfigure}[t]{0.62\linewidth}
\includegraphics[width=\linewidth]{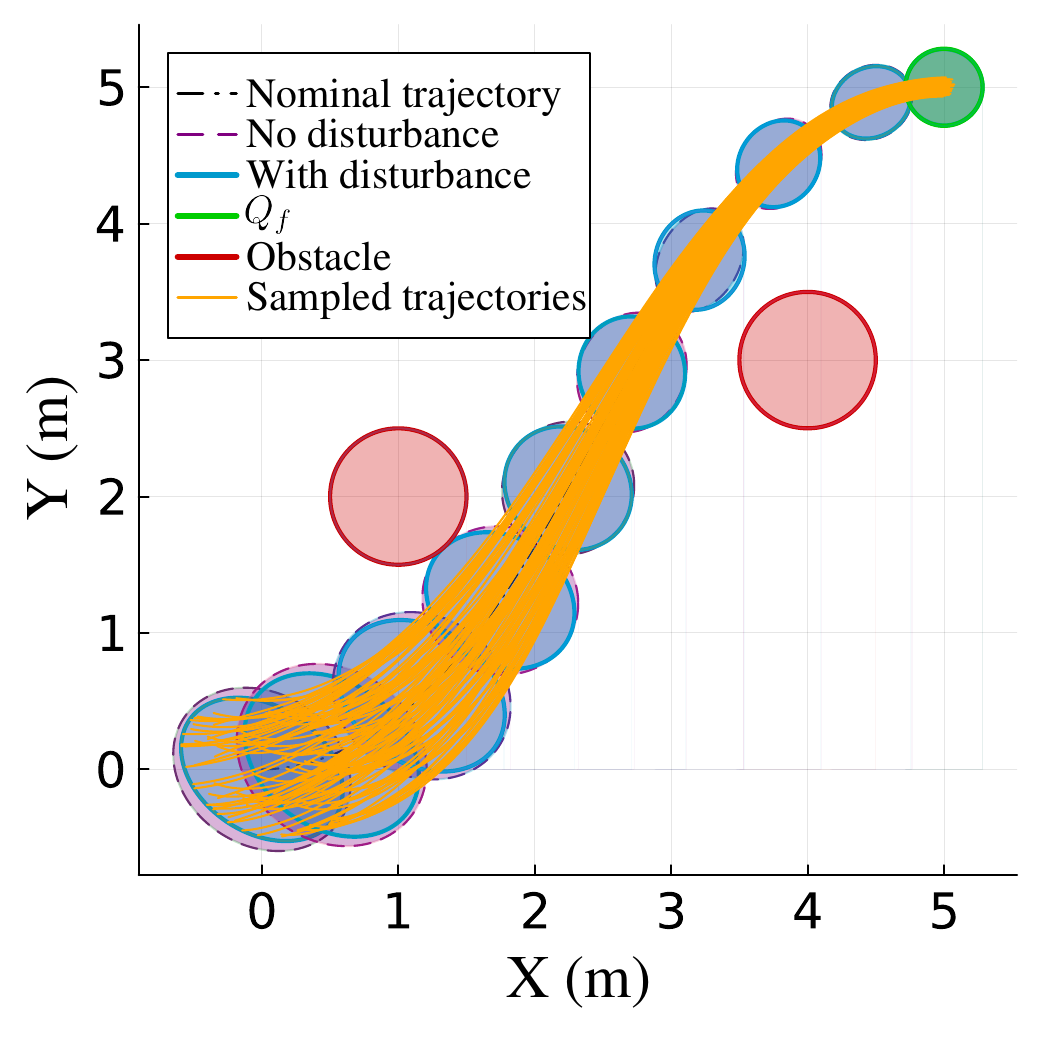}
\end{subfigure}
\begin{subfigure}[t]{0.36\linewidth}
\includegraphics[width=\linewidth]{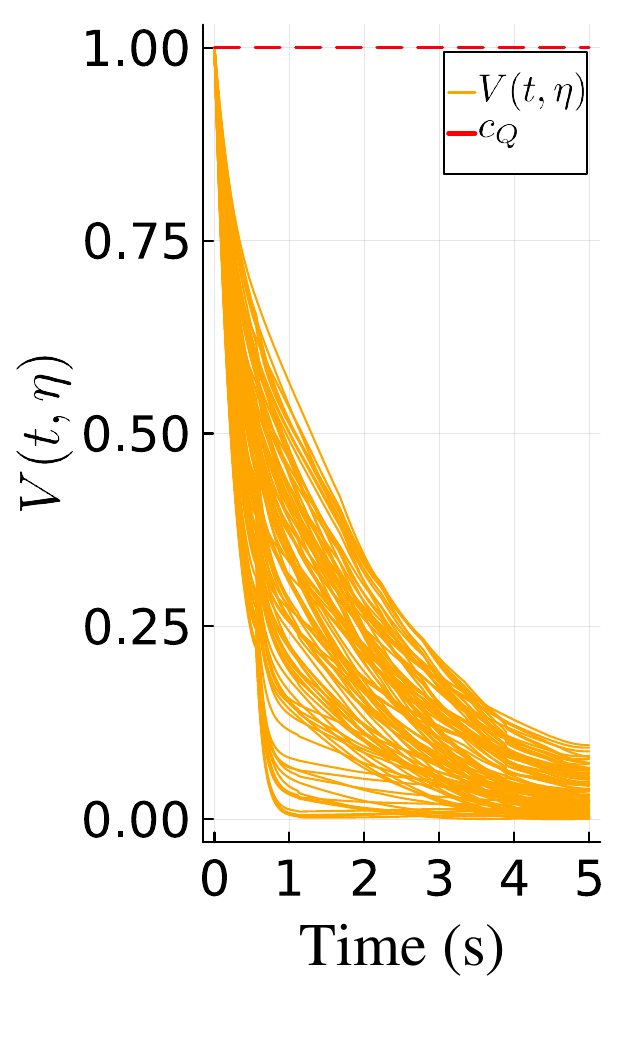}
\end{subfigure}
\caption{(Left): Comparison of synthesized state funnels with and without bounded disturbance; sampled trajectories are generated using the funnel synthesized under disturbance. (Right): The Lyapunov function values of sampled trajectories.} \label{fig:08_funnel} \end{center}
\end{figure}

\begin{figure} \begin{center}  \includegraphics[width=0.95\linewidth,trim={0cm 0cm 0.5cm 0cm},clip]{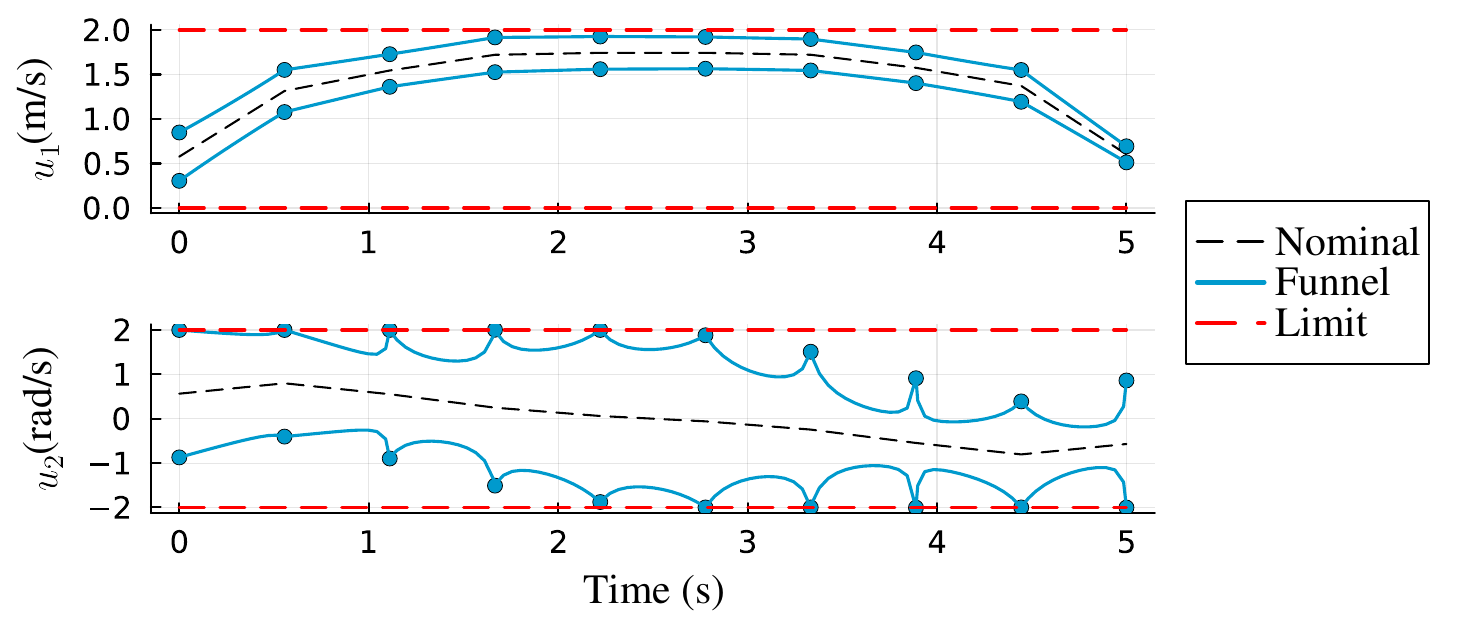} \caption{The synthesized input funnel projected onto each input dimension.} \label{fig:08_input_funnel} \end{center} \end{figure}

The computed the state funnel centered around the nominal state is illustrated in Figure \ref{fig:08_funnel}. It can be seen that the computed funnel without the disturbance is larger than the funnel computed under disturbance. This is expected, as the presence of external disturbances necessitates a more conservative funnel to ensure invariance is maintained despite the disturbance. The computed input funnel centered around the nominal input is illustrated in Figure \ref{fig:08_input_funnel}. The result shows that the input funnel projected onto each input dimension satisfies the input limit constraints over entire horizon. 

\begin{figure}
\centering 
\begin{subfigure}[t]{0.49\linewidth}
\includegraphics[width=\linewidth]{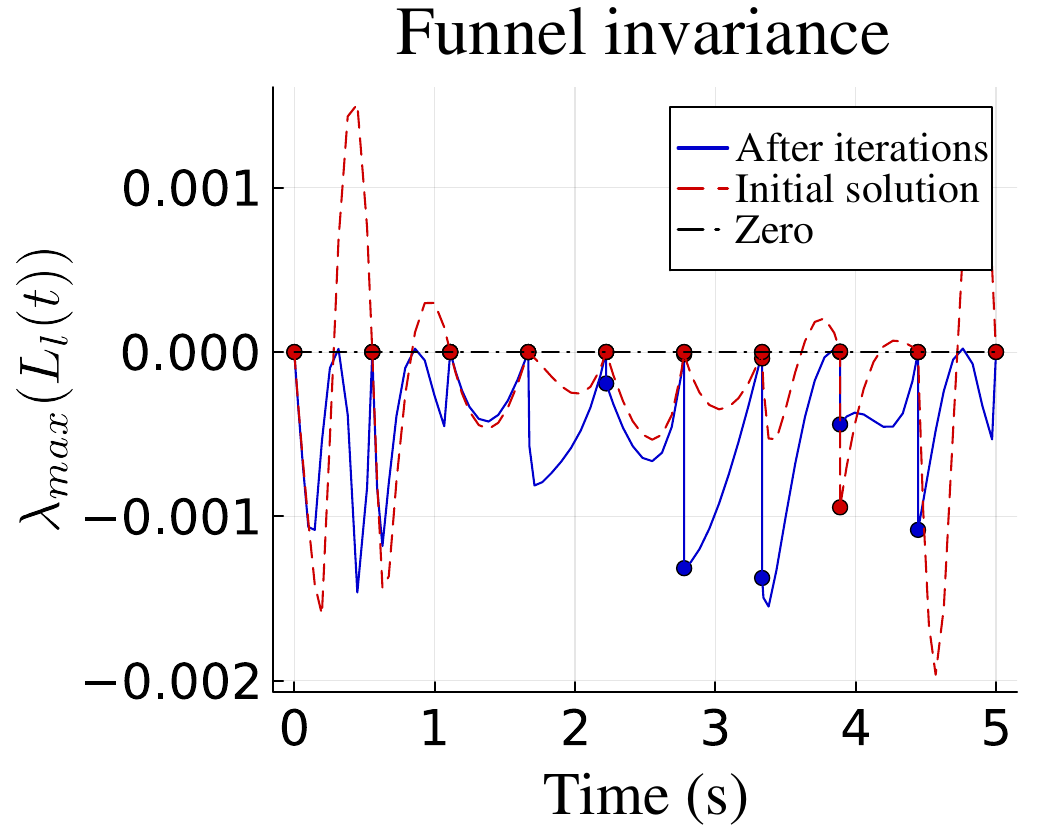}
\end{subfigure}
\begin{subfigure}[t]{0.49\linewidth}
\includegraphics[width=\linewidth]{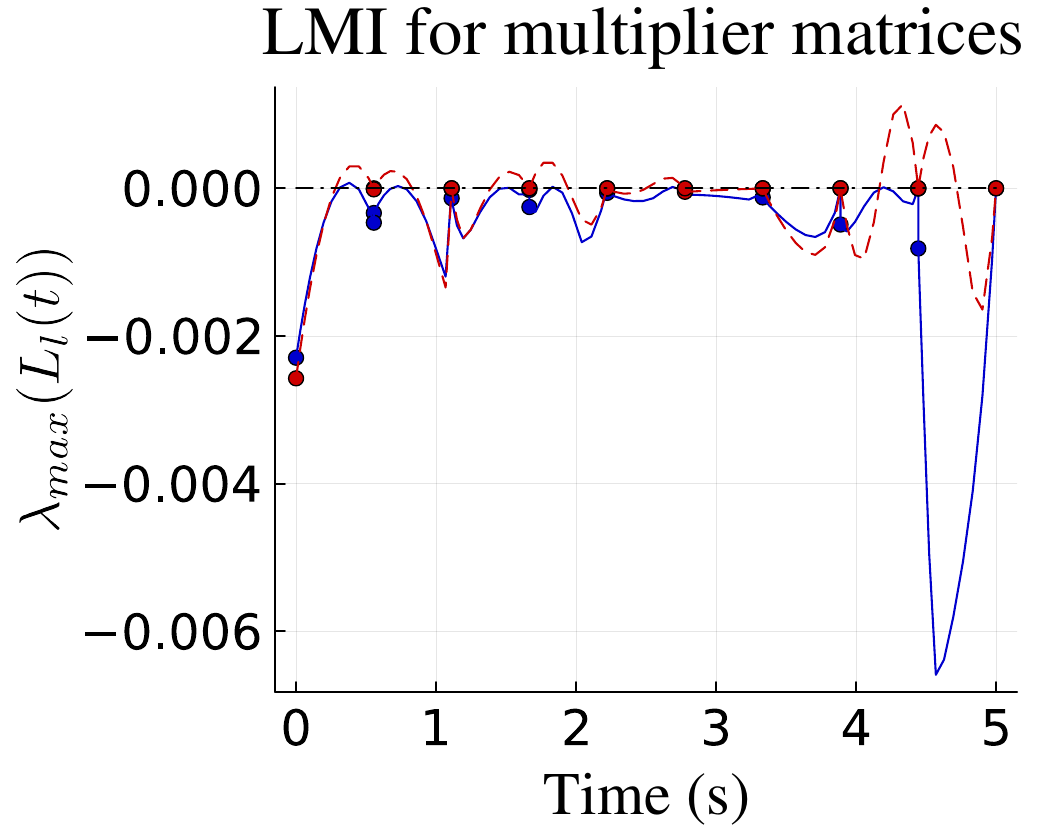}
\end{subfigure}
\caption{Time evolution of maximum eigenvalues for pointwise-in-time LMI constraints in the SCvx-based CTCS approach. Left: Funnel invariance constraint \eqref{eq:Basic_type_H}, Right: valid multiplier constraint \eqref{eq:inclusion_Lsmooth_LMI}.}
\label{fig:07_scvx_ctcs_eigenvalues}
\end{figure}

Finally, we check the invariance of the funnel by generating 100 sample trajectories. The initial states are randomly sampled from the boundary of the computed funnel at $t_{0}$, and the synthesized controller is applied to generate the trajectories. For each sample, a bounded disturbance signal $w(t)$ is randomly generated once at the beginning and held constant over the entire time horizon. This choice of fixed-in-time disturbance is intended to create a more persistent influence on the system, which is shown to be more disruptive than rapidly changing noise. The resulting trajectories are shown in Figure \ref{fig:08_funnel}. To verify funnel invariance, the corresponding Lyapunov function values $V(t,\eta(t))$ are plotted in Figure \ref{fig:08_funnel}. All values remain strictly below the funnel level constant $c_{Q}$, confirming that the invariance condition is satisfied across all sampled trajectories.


\subsection{6-DoF quadrator}

\begin{figure}   \centering   \begin{subfigure}[t]{0.49\linewidth}     \includegraphics[width=\linewidth]{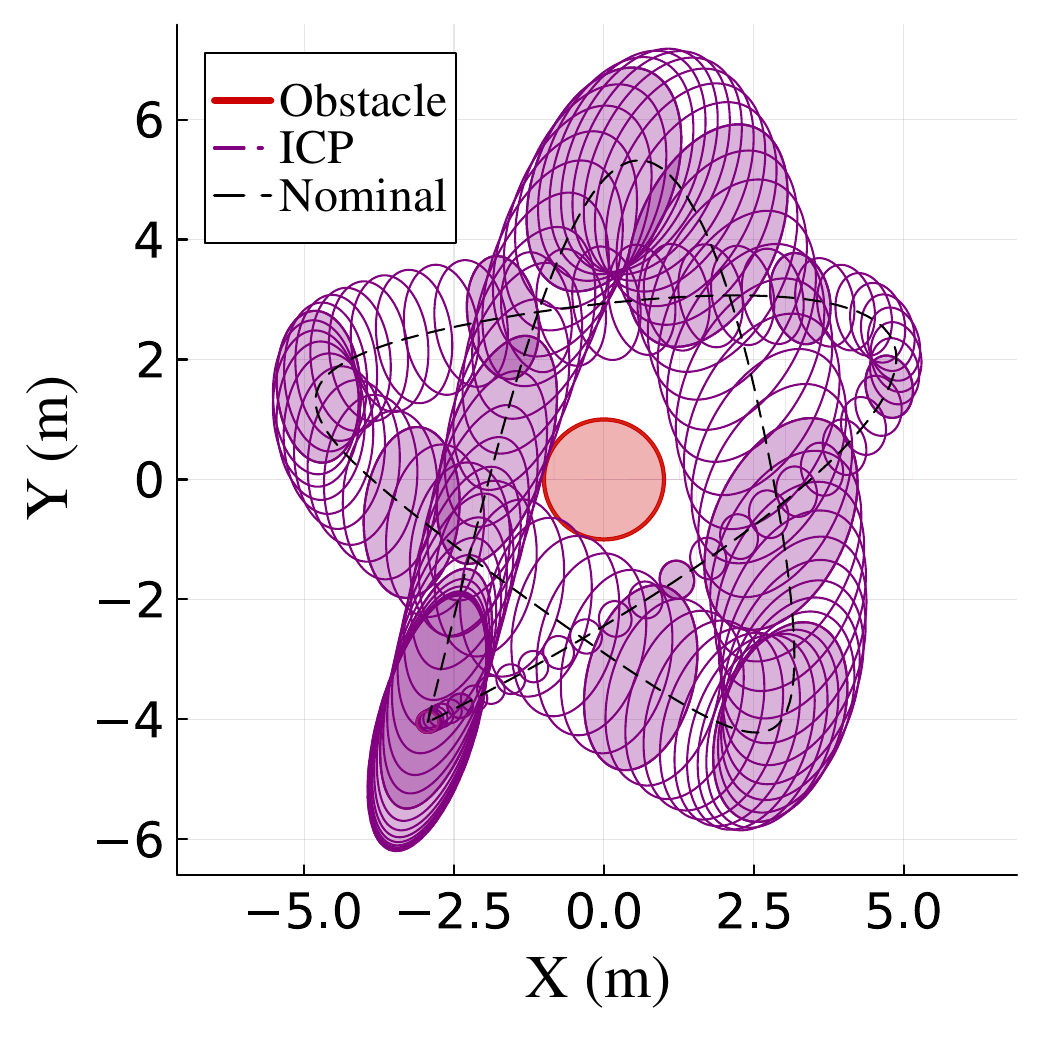}   \end{subfigure}   \hfill   \begin{subfigure}[t]{0.49\linewidth}     \includegraphics[width=\linewidth]{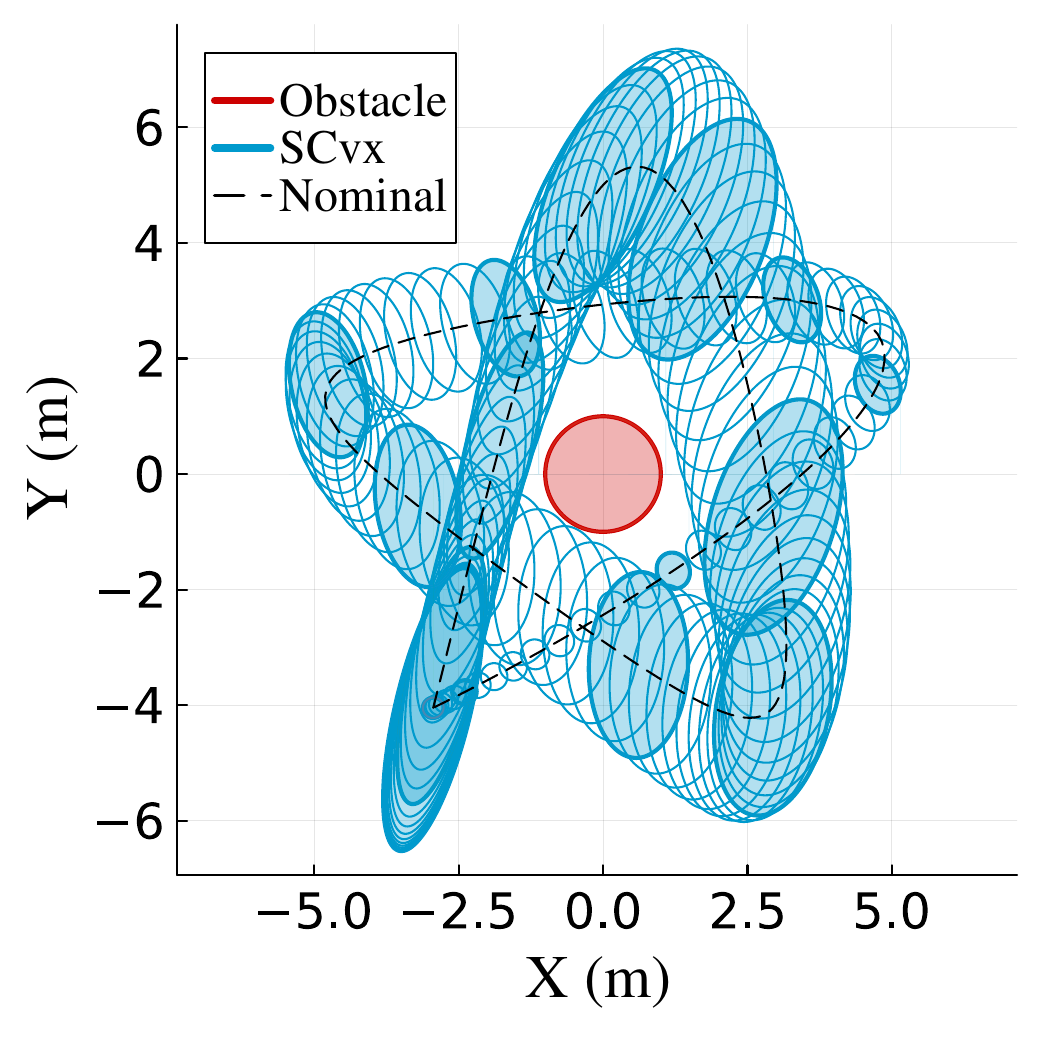}   \end{subfigure}   \caption{State funnels projected onto the $x$ and $y$ plane. Top: result obtained using intermediate checking points. Bottom: result obtained using the SCvx approach. Filled ellipsoids indicate funnels at discrete node points, while unfilled ellipsoids depict intermediate funnels between nodes.}   \label{fig:09_trajectory} \end{figure}

We consider a 6-DoF quadrotor in the East-North-Up (ENU) inertial frame.. The state is $x=(r_{p},v_{p},\Phi,\Omega)$ where $r_{p}\in\mathbb{R}^{3}$ is the position vector, $v_{p}\in\mathbb{R}^{3}$ is the linear velocity, $\Phi\in\mathbb{R}^{3}$ is the Euler angles, and $\Omega\in\mathbb{R}^{3}$ is the body angular velocity. The control input is $u=(F_{z},\tau_{x},\tau_{y},\tau_{z})$ where $F_{z}$ is the thrust along the body $z$-axis, and $(\tau_{x},\tau_{y},\tau_{z})$ are the body-frame torques. The system model is given
\[
\begin{bmatrix}\dot{r_{T}}\\[2pt]
\dot{v_{T}}\\[2pt]
\dot{\Phi}\\[2pt]
\dot{\Omega}
\end{bmatrix}=\begin{bmatrix}v_{T}\\
\frac{1}{m}C_{I/B}(\Phi)F_{B}+g\\
R(\Phi)\Omega\\
J^{-1}\big(\tau-\Omega\times J\Omega\big)
\end{bmatrix},
\]
where $F_{B}=[0,0,F_{z}]^{\top}$, $C_{I/B}(\Phi)$ is the rotation matrix from body to inertia frame, $R(\Phi)$ is the Euler-angle kinematics transformation. Here the following parameters, input constraint set $\mtu$, and the final funnel matrix $Q_{f}$ are used
\begin{align*}
m & =1.325\;(\mathrm{kg}),\;g=[0,0,9.81]^{\top}(\mathrm{m/s^{2}}),\\
J & =\diag\{0.03843,0.02719,0.060528\}\;(\mathrm{kgm^{2}}),\\
\mathcal{U} & =\{u\mid u_{\mathrm{lb}}\le u\le u_{\mathrm{ub}}\},\\
u_{\mathrm{ub}} & =(18,0.1,0.1,0.1),\;u_{\mathrm{lb}}=(0.0,-0.1,-0.1,-0.1)\\
Q_{f} & =\diag\{0.2^{2},0.2^{2},0.2^{2},0.1^{2},0.1^{2},0.1^{2},\\
 & (5^{\circ})^{2},(5^{\circ})^{2},(5^{\circ})^{2},(2^{\circ})^{2},(2^{\circ})^{2},(2^{\circ})^{2}\},
\end{align*}
The state constraint corresponds to obstacle avoidance as depicted in Figure \ref{fig:09_trajectory}.

For the 6-DoF quadrotor dynamics, the first three state equations, corresponding to the position $r_{T}$, are linear in the states and inputs. Therefore, they do not contribute to the nonlinear term $\phi(\cdot)$, and the corresponding rows of the matrix $E$ are zero. Specifically, we have $E=\begin{bmatrix}0_{3\times12}\\
I_{9}
\end{bmatrix}$. We partition the nonlinearities into three channels, each associated with a specific subset of states and inputs: 
\[
q_{[1]}=\begin{bmatrix}\Phi\\
F_{z}
\end{bmatrix},\quad q_{[2]}=\begin{bmatrix}\Phi_{1:2}\\
\Omega_{2:3}
\end{bmatrix},\quad q_{[3]}=\Omega,
\]
where $\Phi_{1:2}$ denotes the first two components of $\Phi$ and $\Omega_{2:3}$ denotes the second and third components of $\Omega$. The functions $\phi_{[1]},\phi_{[2]},\phi_{[3]}\in\mathbb{R}^{3}$ correspond to the nonlinear terms in the dynamics of $v_{T}$, $\Phi$, and $\Omega$, respectively. We use the L-smooth nonlinearity characterization, with the constants $\beta_{1}$, $\beta_{2}$, and $\beta_{3}$ set to $20$, $5$, and $1$, respectively, estimated via sampling. The number of subintervals is set to $N=15$. The nominal trajectory, illustrated in Figure \ref{fig:09_trajectory}, start near $(-3,4)$ follows a star-shaped path, and returns to its starting point. The total time of flight is around 15.78 seconds. 

\begin{figure}
\centering
\begin{subfigure}[t]{0.49\linewidth}
\includegraphics[width=\linewidth]{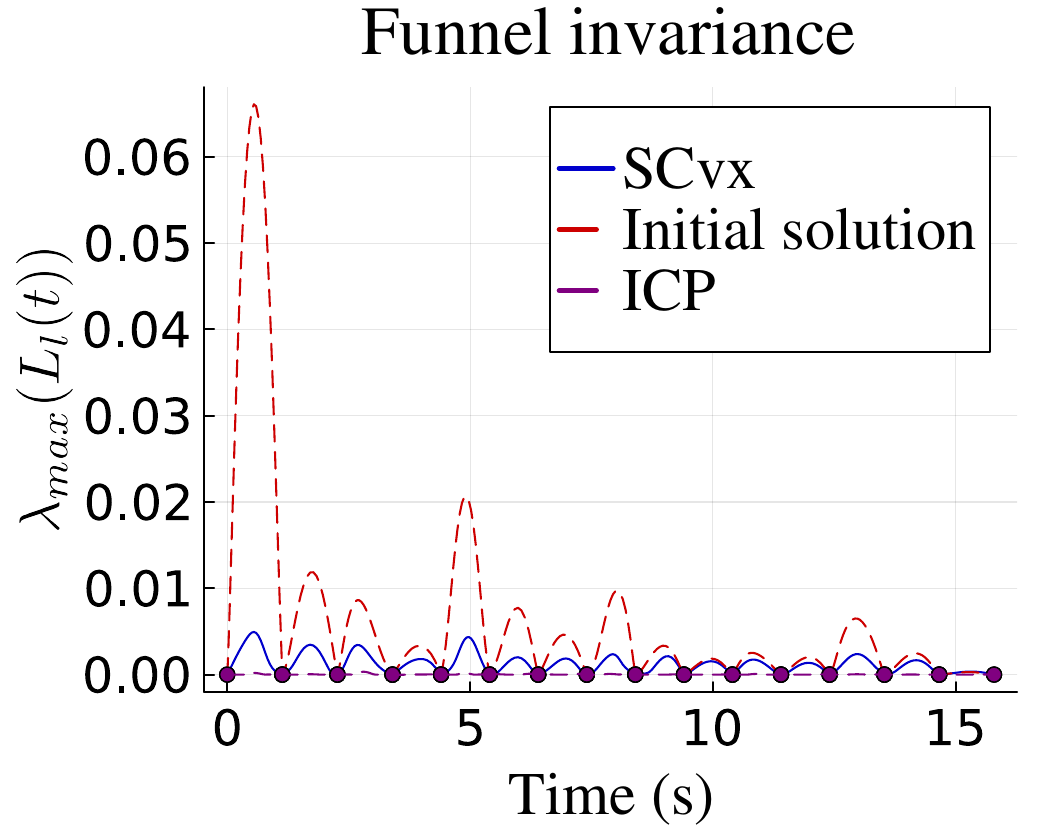}   \end{subfigure}
\begin{subfigure}[t]{0.49\linewidth}   
\includegraphics[width=\linewidth]{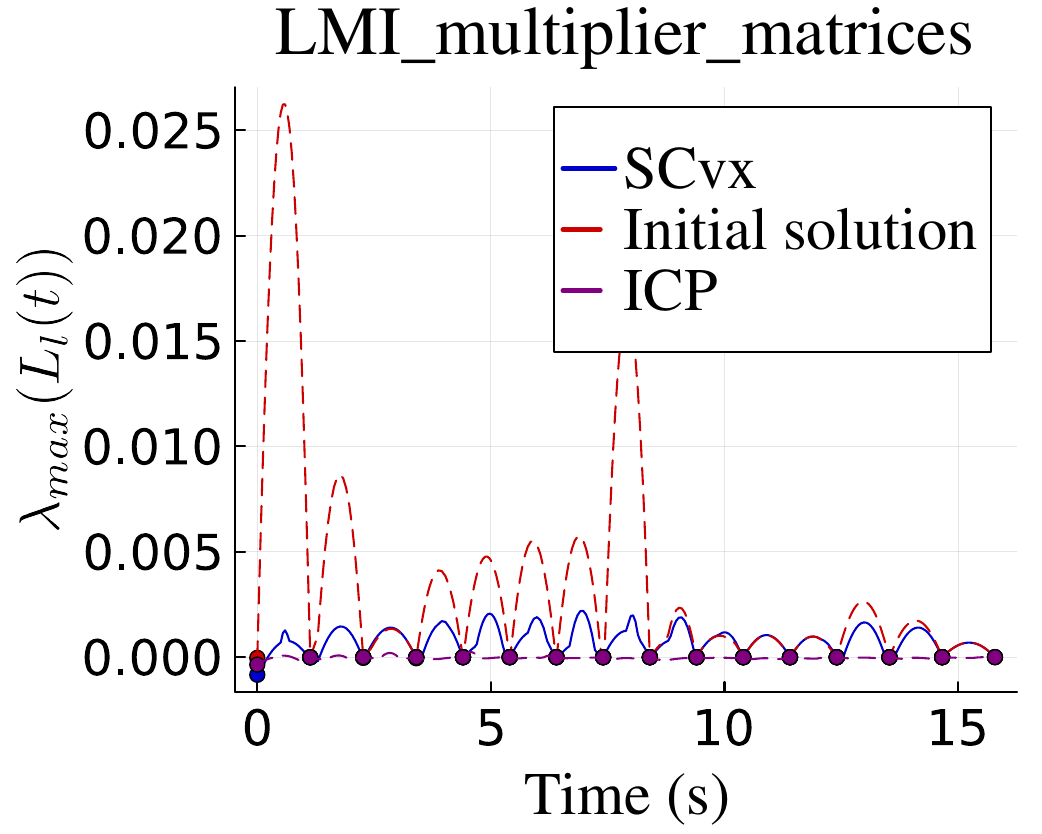}   \end{subfigure}
  \caption{Time evolution of the maximum eigenvalue of each pointwise-in-time LMI constraint $L_l(t)$ associated with the corresponding constraint labeled in each subplot title.}   \label{fig:09_constraint_summary} \end{figure}

We compare two CTCS approaches: (i) introducing intermediate checking points, and (ii) applying SCvx with the subgradient, illustrated in Section \ref{subsec:Intermediate-constraint-checking} and \ref{subsec:Successive-convexification-with}, respectively. Both cases employ the $\fdyn$-type funnel dynamics and have a total 15 time subintervals $N=15$. For the intermediate checking point method, we set $N_{s}=2$. For SCvx, we use the initial guess as the solution of \eqref{eq:funnel_synthesis_discrete} without intermediate checking points. The weights $w_{vc}$ and $w_{tr}$ in \eqref{eq:scvx_cost_function} are set to $2\times10^{5}$ and $2\times10^{3}$, respectively, with $\epsilon=10^{-3}$. The state funnels projected onto the $x$--$y$ coordinates are shown in Figure \ref{fig:09_trajectory}. The resulting costs for the two methods are 76.987 and 75.82, respectively. The total computation time for the former is 68.30, while for the SCvx method, the average computation time per iteration is 23.36 seconds over a total of 50 iterations. The time evolution of the maximum eigenvalues of the pointwise-in-time LMI constraints associated with the invariance condition \eqref{eq:Basic_type_H} and the validity of the multiplier matrices \eqref{eq:inclusion_Lsmooth_LMI} is shown in Figure \ref{fig:09_constraint_summary}.

In both methods, the constraint violation between node points is successfully reduced. For the quadrotor case, the problem dimension is significantly larger than in the unicycle example, both in terms of the system size and the number of discretization nodes $N$, leading to increased computational cost. In this setting, the iterative nature of SCvx can be computationally intensive, making the one-shot SDP approach with intermediate checking points a more efficient alternative.

\section{Conclusion}

This work formulated a continuous-time convex optimization problem for funnel synthesis, with decision variables including time-varying Lyapunov matrices and feedback gains. The formulation combines a DLMI condition to ensure invariance and pointwise-in-time LMI conditions to ensure funnel feasibility. The approach accommodates a broad class of nonlinearities through $\delta$QCs, including sector-bounded and L-smooth dynamics beyond the classical Lipschitz assumption. A solution approach based on numerical optimal control was developed through the notion of funnel dynamics. With the direct-type formulation, it was shown that the positive definiteness of $Q$ can be preserved under constraint \eqref{eq:Z_Z}. To ensure CTCS, two convex methods were proposed: introducing intermediate constraint-checking points and a SCvx approach with subgradients. For the latter, the existence of a valid subgradient for the given constraint was established, enabling its integration into the SCvx framework.

Numerical simulations showed that funnel computation under the L-smooth condition yields less conservative results compared to the Lipschitz condition. Between two funnel dynamics formulations, the Lyapunov-type is not able to preserve the positive definiteness of $Q$ as anticipated, while the direct-type preserved it as proven. Both CTCS methods, intermediate checking points and SCvx with subgradients, successfully reduced constraint violations between node points. For high-dimensional and long-horizon problems such as the quadrotor case, the repeated SDP solving required by SCvx was computationally heavy, making the intermediate-checking-point approach more efficient in practice.



\section*{References}

\bibliographystyle{IEEEtran}
\bibliography{main}

\end{document}